\documentclass[12pt, letter]{article}

\usepackage[left=1in, right=1in, top=1in, bottom=1in]{geometry} 

\usepackage{amsmath}
\usepackage{amsthm}
\usepackage{amssymb}
\usepackage{graphicx}
\usepackage{enumerate}
\usepackage{natbib,bibunits}
\usepackage{enumitem}
\usepackage{setspace}
\usepackage{threeparttable}
\usepackage{xfrac}
\usepackage{braket}
\usepackage{pdflscape}
\usepackage{comment}

\usepackage{tocloft}
\usepackage{titletoc}

\usepackage{xcolor}
\definecolor{darkblue}{rgb}{0.0, 0.0, 0.4} 
\definecolor{mygreen}{RGB}{0,100,0} 

\AtBeginDocument{\hypersetup{citecolor=blue,linkcolor=blue,urlcolor=blue}}

\usepackage{amsfonts}       
\usepackage{dsfont}         
\usepackage{bm}             
\usepackage{pifont}         

\usepackage{caption}        
\usepackage{subcaption}     

\usepackage{rotating}       
\usepackage{booktabs}       
\usepackage{multirow}       
\usepackage{ulem}           
\usepackage{makecell}
\usepackage{rotating}
\useunder{\uline}{\ul}{}    
\usepackage{longtable}      

\usepackage[colorinlistoftodos,prependcaption,textsize=footnotesize]{todonotes}

\def\spacingset#1{\renewcommand{\baselinestretch}{#1}\normalsize}
\allowdisplaybreaks[2]

\newcommand{\var}{\mathbb{V}}

\newcommand{\eps}{\varepsilon}

\newcommand{\1}{\mathds{1}}

\newcommand{\dd}{\mathrm{d}}

\newcommand{\R}{{\mathbb{R}}}
\newcommand{\D}{{\mathbb{D}}}
\newcommand{\Bm}{{\mathcal{B}}}
\newcommand{\Fm}{{\mathcal{F}}}
\newcommand{\Nm}{{\mathcal{N}}}
\newcommand{\Pm}{{\mathcal{P}}}
\newcommand{\Ym}{{\mathcal{Y}}}

\newcommand{\wlB}{{\mathrm{wl}_B}}

\newcommand{\ymdel}{{B}}
\newcommand{\ymdelc}{{B^c}}
\newcommand{\ymdellong}{{B_{\delta}(y_t)}}

\newcommand{\tbu}{{t\vert t-1}} 
\newcommand{\tpu}{{t\vert t}} 

\newcommand{\trace}{\operatorname{tr}}

\newcommand{\score}{{s}}
\newcommand{\scaling}{\mathcal{S}_{t-1}}
\newcommand{\Hessian}{{H}}
\newcommand{\zeromatrix}{O_k}
\newcommand{\unitmatrix}{I_k}

\newcommand{\kl}{\mathsf{KL}} 
\newcommand{\rkl}{{\mathsf{TKL}}} 
\newcommand{\ckl}{{\mathsf{CKL}}} 
\newcommand{\ekl}{\mathsf{EKL}}
\newcommand{\cev}{\mathsf{CEV}} 
\newcommand{\gmm}{\mathsf{EGMM}} 
\newcommand{\mse}{\mathsf{MSE}}

\newcommand{\cmark}{\color{mygreen} \ding{51}}  
\newcommand{\xmark}{\color{red} \ding{55}}  

\usepackage[colorlinks,citecolor=blue]{hyperref}

\renewcommand{\P}{\mathbb{P}}
\newcommand{\E}{\mathbb{E}}

\newtheorem{theorem}{Theorem}

\newtheorem{corollary}{Corollary}
\newtheorem{proposition}{Proposition}
\newtheorem{definition}{Definition}
\newtheorem{example}{Example}
\newtheorem{remark}{Remark}

\newtheorem{assumption}{Assumption}

\newtheorem{assumptionalt}{Assumption}[assumption]

\newenvironment{assumptionp}[1]{
  \renewcommand\theassumptionalt{#1}
  \assumptionalt
}{\endassumptionalt}

\newcommand{\appr}[1]{Appendix~\ref{#1}}

\begin{document}
		
\title{\textbf{Expected Kullback-Leibler-based characterizations of score-driven updates}}

 \author{
    \textsc{Ramon de Punder}\textsuperscript{a,d},\;
    \textsc{Timo Dimitriadis}\textsuperscript{b,e},\; 
     \\ \text{and}
    \textsc{Rutger-Jan Lange}\textsuperscript{c,d}
  }

  \maketitle
  
   \begin{center}
    \small
    \textsuperscript{a}\textit{Department of Quantitative Economics,	University of Amsterdam, The Netherlands} \\
    \textsuperscript{b}\textit{Faculty of Economics and Business, Goethe University Frankfurt, Germany} \\
    \textsuperscript{c}\textit{Econometric Institute,  	Erasmus School of Economics, Rotterdam, Netherlands}
    \\
     \textsuperscript{d}\textit{Tinbergen Institute, The Netherlands}
         \\
     \textsuperscript{e}\textit{Heidelberg Institute for Theoretical Studies, Heidelberg, Germany }
  \end{center}

  \renewcommand{\thefootnote}{\fnsymbol{footnote}}
  \footnotetext[1]{
  E-mail addresses: \href{mailto:r.f.a.depunder@uva.nl}{r.f.a.depunder@uva.nl} (Ramon de Punder), \href{mailto:dimitriadis@econ.uni-frankfurt.de}{dimitriadis@econ.uni-frankfurt.de} (Timo Dimitriadis),  \href{mailto:lange@ese.eur.nl}{lange@ese.eur.nl} (Rutger-Jan Lange).}
  \renewcommand{\thefootnote}{\arabic{footnote}}


\begin{abstract}
	\noindent 
    Score-driven (SD) models are a standard tool in statistics and econometrics, with applications in hundreds of published articles in the past decade. We provide an information-theoretic characterization of SD updates based on reductions in the expected Kullback–Leibler (EKL) divergence relative to the true---but unknown---data-generating density. EKL reductions occur if and only if the expected update direction aligns with the expected score; i.e., their inner product should be positive. This equivalence condition uniquely identifies SD updates (including scaled or clipped variants) as being EKL reducing, even in non-concave, multivariate, and misspecified settings. We further derive explicit bounds on admissible learning rates in terms of score moments, linking SD methods to adaptive optimization techniques. By contrast, alternative performance measures in the literature impose stronger conditions (e.g., concave logarithmic densities) and do not {\it characterize} SD updates: other updating rules may improve these measures, while SD updates need not. Our results provide a rigorous justification for SD models and establish EKL as their natural information-theoretic foundation.
\end{abstract}

\noindent%
\vfill

\pagebreak 

\setlength{\textheight}{8.75in} 
\renewcommand{\baselinestretch}{1.66}
\normalsize
\setlength\abovedisplayskip{6pt}
\setlength\belowdisplayskip{6pt}

\normalem

\begin{bibunit}[chicago]

\section{Introduction}
\label{sec:introduction}

The use of score-driven (SD) models has proliferated over the last decade. They were originally introduced by \cite{creal2013generalized} and \cite{harvey2013dynamic} and known by different names and acronyms; recent literature (e.g., \citealp{artemova2022score1,artemova2022score2}; \citealp{harvey2022score}) has converged on the {terminology} of SD models. These models specify a distribution with time-varying parameters---governing, for example, intensity, location, scale, or shape---whose dynamics are driven by the \emph{score}, the derivative of the log-likelihood with respect to the time-varying parameter vector. 
SD updates can be interpreted as stochastic gradient ascent steps that adjust the parameter after each observation. Unlike in optimization, the parameter does not converge but remains perpetually responsive. Applications are numerous; see \href{www.gasmodel.com}{www.gasmodel.com} for a list of over 400 publications.

Much of this literature assumes that the SD filter coincides with the true data-generating process. 
Yet it remains unclear whether SD updates possess theoretical properties that uniquely characterize them in more general, possibly misspecified settings, and, if so, what these are. In this paper, we resolve this question by showing that, in expectation, sufficiently small parameter adjustments improve the distributional fit if and only if the update direction is driven by the score, precisely the principle on which SD models are based.

To explain this result, let $y_t \in \mathcal{Y}$ with $t \in \mathbb{N}$ denote the realization of $Y_t$ with true density $p_t$ given the information at time $t-1$. The researcher postulates (typically misspecified) parametric densities $f_{t|t-1}(\cdot) \equiv f(\cdot | \vartheta_{t|t-1})$ before updating and $f_{t|t}(\cdot) \equiv f(\cdot | \vartheta_{t|t})$ after updating. 
Here, $\vartheta_{t|t-1}$ and $\vartheta_{t|t} \equiv \vartheta_{t|t}(y_t)$ denote the predicted and updated parameters based on information up to times $t-1$ and $t$, respectively; i.e., the update takes account of the latest observation $y_t$.
The SD update is
\begin{equation}
\label{eq:SD update}
\vartheta_{t|t} \equiv \vartheta_{t|t}(y_t) 
= \vartheta_{t|t-1} + A \mathcal{S}_{t-1} \, \score(y_t,\vartheta_{t|t-1}),
\end{equation}
where $\score(y_t,\vartheta_\tbu) := (\partial/\partial\vartheta) \log f(y_t \vert \vartheta)\vert_{\vartheta_\tbu}$ is the score (i.e., a stochastic gradient). As is standard (e.g., \citealp{creal2013generalized}), the update involves a static learning-rate matrix $A$ and a dynamic scaling matrix $\mathcal{S}_{t-1}$, assumed known at time $t-1$.
We assume their product $A \mathcal{S}_{t-1}$ to be positive definite, which is not automatic even if both $A$ and $\scaling$ are positive definite. The prediction step used to construct $\vartheta_{t+1|t}$ from $\vartheta_{t|t}$ is not essential at present; its discussion is postponed until after Definition~\ref{def:NewtonUpdate}.

To quantify the distributional closeness between $f_\tpu$ and $p_t$, we propose to use the \emph{expected} KL (EKL) divergence, where, in addition to the usual integral, the observation $y \in \Ym$ used in the update $\vartheta_{t|t}(y)$ and thus in $f_\tpu(\cdot) \equiv f(\cdot|\vartheta_{t|t}(y))$ is averaged out using the true density:
\begin{align}
	\label{eqn:EKLDefIntro}
	\ekl(p_t \Vert f_{\tpu})
	:= \int_{\Ym} \int_{\Ym} 
	\log \!\left(\frac{p_t(x)}{f(x \vert \vartheta_{t|t}(y))} \right) 
	p_t(x)\, p_t(y)\, \dd x \,\dd y.
\end{align}
The EKL divergence measure \eqref{eqn:EKLDefIntro} features a double integral over $x$ and $y$, admitting a natural two-sample interpretation: one draw $y$ updates the model density, and an independent redraw $x$ evaluates the updated model’s fidelity on new data. By averaging out the uncertainty in both samples, it provides a natural criterion for evaluating updating rules and combines the realized KL perspective of \cite{blasques2015information} with the expectation-based measures of \cite{Gorgi2023} and \cite{Creal2024GMM} (see related literature below).

As Theorems~\ref{thm:EKLequivalenceIndepCopy}--\ref{thm:EKLequivalenceBounded} show, under various conditions on the score's derivative (the Hessian), any sufficiently small parameter change $\vartheta_\tpu(Y_t) - \vartheta_\tbu$ implies an EKL improvement, 
\begin{align}
    \label{eqn:EKLComparisonIntro}
	\ekl(p_t \Vert f_{\tpu}) < \ekl(p_t \Vert f_{\tbu}),
\end{align}
\emph{if and only if} the expected parameter adjustment is aligned with the expected score, i.e.,
\begin{align}	\label{eqn:ExpectedScoreEquivalence_Intro}
	\E_{p_t} \!\big[ \vartheta_\tpu(Y_t) - \vartheta_\tbu \big]^\top \E_{p_t} \!\big[ \score(X_t,\vartheta_\tbu) \big] > 0,
\end{align}
where $\E_{p_t}[\cdot]$ denotes the expectation under the true density $p_t$. 
The equivalence between~\eqref{eqn:EKLComparisonIntro} and~\eqref{eqn:ExpectedScoreEquivalence_Intro} provides an information-theoretic characterization of a family of updating rules. We refer to any (user-specified) updating rule $\vartheta_{t|t}(Y_t)$ satisfying~\eqref{eqn:ExpectedScoreEquivalence_Intro} as being \emph{score equivalent in expectations} (SEE); we use the plural as~\eqref{eqn:ExpectedScoreEquivalence_Intro} involves two distinct expectations. Since $p_t$ is treated as unknown, we would ideally like~\eqref{eqn:ExpectedScoreEquivalence_Intro} to hold for broad classes of $p_t$.

For SD updates~\eqref{eq:SD update}, in which
$\E_{p_t}[\vartheta_\tpu(Y_t)-\vartheta_\tbu] = A \scaling \E_{p_t}[s(Y_t,\vartheta_\tbu)]$, condition~\eqref{eqn:ExpectedScoreEquivalence_Intro} simplifies to
\begin{equation}
\label{eq:our score equiv}
  \E_{p_t}[s(Y_t,\vartheta_\tbu)]^\top
  A \scaling
  \E_{p_t}[s(Y_t,\vartheta_\tbu)] > 0,
\end{equation}
which holds independently of the true law $p_t$ whenever the expected score is non-zero.
Hence, sufficiently small SD updates with non-zero expectation are EKL reducing. Indeed, SD updates may be the practically most important class of SEE updates. Positive definiteness (and thus symmetry) of $A\mathcal{S}_{t-1}$, although not customary (e.g., \citealp{gasperoni_score-driven_2023,d2024modeling}), is essential in producing gradient-ascent directions and EKL improvement. 
The condition that the expected score at $\vartheta_{\tbu}$ is nonzero can be interpreted as requiring that $\vartheta_{\tbu}$ is not a stationary point of the map $\vartheta \mapsto \E_{p_t}[\log f(Y_t| \vartheta)]$; in particular, no improvement guarantees can be obtained if the prediction $\vartheta_\tbu$ is already optimal.

The broad scope of Theorems~\ref{thm:EKLequivalenceIndepCopy}--\ref{thm:EKLequivalenceBounded} stems from the intrinsic relationship between the EKL criterion~\eqref{eqn:EKLDefIntro} and the driving mechanism in SD models, namely the score. The score is the derivative of the log likelihood, while the EKL criterion can be written in terms of the expected log-likelihood contribution; hence, they are intimately related. This connection allows our characterization results to extend to multivariate settings under mild conditions, accommodating misspecification, non-concave model log-densities, and general (positive-definite) matrix combinations $A \scaling$. As we show in Section~\ref{sec:ComparisonOptimality}, other performance criteria for SD updates are less directly tied to the score and more restricted in scope.

Our improvement guarantee in~\eqref{eqn:EKLComparisonIntro} requires that the parameter adjustment be sufficiently small. To provide practical guidance on the admissible step size, Theorem~\ref{thm:UpperBoundLearningRate} further derives (non-infinitesimal) upper bounds on elements or eigenvalues of the matrix combination $A\scaling$ that still ensure EKL improvements. These bounds depend on the first two population moments of the score, sharpening the constant bound of \citet{Gorgi2023} and motivating the use of adaptive, moment-based learning rates for SD models, in analogy with adaptive methods in the optimization literature. 

\begin{table}[tb]
\centering
\caption{Overview of different performance criteria considered in this paper. \label{tab:OverviewCriteria}}
\begin{footnotesize}
\begin{threeparttable}
\begin{tabular}{ll cc l}
    \toprule
    \multirow{2}{*}{Criterion} & \multirow{2}{*}{Reference}  & Proper & Constructive & \multicolumn{1}{c}{Assumption related to} \\
    & & criterion & ``$\Longleftrightarrow$'' &  \multicolumn{1}{c}{expected Hessian}\\
    \midrule 
    $\ekl$ & this paper: Theorem~\ref{thm:EKLequivalenceIndepCopy} & \cmark &  \cmark & bounded (\ref{ass:HB}) \\
    $\ekl$ & this paper: Theorem~\ref{thm:EKLequivalenceBounded} & \cmark & \cmark & locally bounded (\ref{ass:HLB})
    \\
    \midrule
    $\cev$ & \citet{Gorgi2023} & \cmark & \xmark &  negative~definite  (\ref{ass:HN})\\
    $\mse$ & \citet{Gorgi2023} & \cmark & \xmark &  negative~definite (\ref{ass:HN}) \\
    $\gmm$ & \citet{Creal2024GMM} & \cmark & \xmark & bounded \& negative~definite (\ref{ass:HBNT})   \\
    \midrule
    $\rkl$ & \citet{blasques2015information} & \xmark &    \\
    $\ckl$ & this paper: Appendix~\ref{sec:CKL} & \cmark & \xmark &  \\
    \bottomrule
\end{tabular}
 \vspace{3pt}
{NOTE: As the TKL measure is not a proper divergence measure, the remaining columns are left blank. For EGMM, we additionally need a uniform bound on the third derivative of the log-density.
Table~\ref{tab2} verifies these assumptions for specific model densities. 
}
\end{threeparttable}
\end{footnotesize}
\end{table}

\textbf{Related literature.} We compare the EKL measure with four related performance measures for SD models proposed in three closely related papers;\footnote{ 
These papers focus on providing improvement guarantees for a \emph{single} time step. This differs from alternative approaches based on in-fill asymptotics (e.g.,  \citealp{beutner2023consistency}), which use the Kullback–Leibler divergence to characterize the pseudo-true path in the continuous-time limit.}
see Table~\ref{tab:OverviewCriteria} and Section~\ref{sec:ComparisonOptimality}.

First, \cite{Gorgi2023} show that SD updates move the parameter toward the pseudo-true value, as measured by the first and second moments of the update around the pseudo-true parameter, yielding their conditional expected variation (CEV) and mean squared error (MSE) criteria. Second, \cite{Creal2024GMM} define an expected generalized method-of-moments (EGMM) loss using the expected score as a moment condition and show that SD updates reduce this loss under a scaling that requires a (practically infeasible) expectation under the true density. Third, \citet{blasques2015information} propose a trimmed KL (TKL) measure which, despite known issues (e.g., \citealp{blasques2018information}), remains a standard reference for motivating SD models (e.g., \citealp[p.~1653]{holy2022modeling}, \citealp[p.~1014]{delle2023modeling}, \citealp[p.~1069]{gasperoni_score-driven_2023}, and \citealp[p.~4]{catania2026unobserved}).

For the first two related papers (i.e., \citealp{Gorgi2023} and \citealp{Creal2024GMM}), we show in Sections~\ref{sec:CEV}--\ref{sec:EGMM} that their CEV, MSE, and EGMM performance guarantees for SD updates effectively require log-concave model densities, thereby excluding many fat-tailed distributions (e.g., the Student’s $t$ distribution with time-varying location) that fall within the scope of our EKL criterion. Moreover, neither \cite{Gorgi2023} nor \citet{Creal2024GMM} fully characterize the model classes that improve their performance criteria; here, we generalize their results by establishing characterizations in a multivariate setting (Propositions~\ref{prop:CEVEquiv}--\ref{prop:EGMMEquiv}). 

In the multivariate case, we show that the approaches in \cite{Gorgi2023} and \cite{Creal2024GMM} impose stringent restrictions on the matrix product $A\scaling$, often forcing it to be a scalar multiple of the identity. Moreover, the associated equivalence conditions are of limited use for model construction, as they effectively prescribe adjustment toward the \mbox{(pseudo-)}true parameter.
The relatively greater complexity and narrower scope of these results highlight the appeal of our SEE condition~\eqref{eqn:ExpectedScoreEquivalence_Intro}, which via~\eqref{eq:our score equiv} naturally leads to the general class of (multivariate) SD updates.

For the third related paper (i.e., \citealp{blasques2015information}), we show in Section~\ref{sec:TKLComparison} that the TKL measure is based on a localized, trimmed KL divergence, which is known to be problematic in the scoring-rule literature (e.g., \citealp{diks2011likelihood,gneiting_comparing_2011}). Consequently, the associated equivalence condition is insensitive to the true density $p_t$ and thus unsuitable for gauging closeness to $p_t$. To address this within the localization framework of \cite{blasques2015information}, we replace trimming with censoring, obtaining a censored KL (CKL) measure that preserves the fundamental properties of a KL divergence, but still yields an unhelpful equivalence condition (see Appendices~\ref{sec:CKL} and~\ref{sec:CKLGeneral} for details).

\textbf{Notation.} Vectors are columns and $\|\cdot\|$ denotes the Euclidean norm for vectors and the spectral norm, that is, the operator norm induced by the Euclidean norm, for matrices. For a matrix $A$, $\trace(A)$ is its trace, and $\lambda_{\min}(A), \lambda_{\max}(A)$ its extremal eigenvalues. 
We use $\zeromatrix$ and $\unitmatrix$ for the $k \times k$ zero and identity matrix, respectively.
The Loewner order for symmetric matrices is denoted $A \succeq (\succ)\, B$. Note that $-c\unitmatrix \preceq A \preceq c\unitmatrix$ if and only if $\|A\|\leq c$ for $c>0$.

\section{Expected Kullback-Leibler reducing updates}
\label{sec:EKL}

After introducing notation and preliminaries in Section~\ref{sec:Preliminaries}, we establish our main result in Section~\ref{sec:MainResultEKL}, whose assumptions are relaxed in Section~\ref{sec:LocallyBoundedHessians}.

\subsection{Preliminaries}
\label{sec:Preliminaries}

We consider an outcome space $\Ym \subseteq \mathbb{R}^l$ with $l \in \mathbb{N}$ and a stochastic process $\{Y_t: \Omega \to \Ym\}_{t=1}^T$ on a complete probability space $(\Omega, \Fm, \P)$ with a measurable space $(\Ym^T, \Bm(\Ym^T))$.
We consider the flexible and user-chosen sigma-algebra (or information set) $\Fm_{t-1}$, which is such that  $\sigma(Y_s: s \le t-1) \subseteq \Fm_{t-1} \subseteq \Fm$, but $Y_t$ is not $\Fm_{t-1}$-measurable.
While we typically have $\Fm_{t-1} = \sigma(Y_s : s \le t-1)$ as, e.g., in \citet{Gorgi2023} and \citet{Creal2024GMM}, allowing $\Fm_{t-1}$ to differ from $\sigma(Y_s : s \le t-1)$ is useful for incorporating external covariates or when the data are generated by a state-space model, in which case the true latent state at time~$t$ can be included in $\Fm_{t-1}$ (for details, see Appendix~\ref{app:statespaceDGP}).

The true conditional law of $Y_t$ given $\Fm_{t-1}$ is $P_t \in \Pm \subseteq \Pm_0$, where $\Pm_0$ is the class of absolutely continuous distributions, with density $p_t$.
The subclass $\Pm$ may vary across the results below.
Probabilities and variances with respect to $p_t$ are denoted $\P_{p_t}(\cdot)$ and $\mathbb{V}_{p_t}(\cdot)$, respectively, and we occasionally write $Y_t \sim p_t$ or $p_t \in \Pm$.  
Statements involving random variables are meant to hold $P_t$-almost surely (a.s.) unless stated otherwise. 

The researcher’s (possibly misspecified) predictive density for $Y_t$ is $f_\tbu(\cdot) \equiv f(\cdot|\vartheta_\tbu)$, where $\vartheta_\tbu \in \Theta \subseteq \mathbb{R}^k$ is based on $\mathcal{F}_{t-1}$ and $\Theta$ is open and convex. 
A link function may be embedded in $f(\cdot|\vartheta_\tbu)$, mapping $\vartheta_\tbu$ into some desired domain \citep[pp.\,324--326]{harvey2022score}. 
Throughout, we assume that the support of $p_t$ is a subset of the support of $f(\cdot|\vartheta)$, for all $\vartheta \in \Theta$.
Dependence on static parameters or exogenous covariates known at time $t-1$ is permitted but suppressed for readability. 

Given the realization $y_t$ of $Y_t$, the updated parameter is $\vartheta_\tpu = \phi(y_t,\vartheta_\tbu)$, where $\phi:(\Ym,\Theta)\to\Theta$ is the \emph{updating rule}. We write $\Delta\phi(y_t,\vartheta_\tbu):=\vartheta_\tpu-\vartheta_\tbu$. To guarantee $\vartheta_\tpu\in\Theta$ almost surely, we typically take $\Theta=\mathbb{R}^k$ and employ a suitable link function if needed.
Following \cite{lange2022robust}, we distinguish between the \emph{update} step, which yields $\vartheta_{t| t}$, and the \emph{prediction} step, which yields $\vartheta_{t+1 | t}$. The update $\vartheta_{t| t}$ is designed to improve the fit to the current density $p_t$ from which $y_t$ is drawn, whereas the prediction $\vartheta_{t+1| t}$ is aimed at the next-period density $p_{t+1}$, for which no observations have (yet) been collected.
Ideally, the updated density $f_{\tpu}(\cdot)\equiv f(\cdot|\vartheta_\tpu)$ improves upon $f_\tbu(\cdot)$ by maximizing $\E_{p_t}[\log f(Y_t|\vartheta)]$ over $\vartheta \in \Theta$.
However, this expectation is unobserved and cannot be estimated, as only a single realization $y_t$ of $Y_t$ is available. Stochastic-gradient methods therefore rely on the observed gradient, the \emph{score},
\[
\score(y_t,\vartheta_\tbu):=\left.\frac{\partial}{\partial\vartheta}\log f(y_t | \vartheta)\right|_{\vartheta_\tbu},
\]  
leading to the class of score-driven (SD) filters (e.g., \citealp{creal2013generalized, harvey2013dynamic}).  

\begin{definition}[Score-driven update]
    \label{def:NewtonUpdate} 
    The SD  update is
    \begin{align}
    \label{eqn:SSDUpdate}
    \vartheta_\tpu^\mathrm{SD}  
    = \phi_{\mathrm{SD}}(y_t,\vartheta_\tbu) 
    := \vartheta_\tbu + A \scaling \score(y_t,\vartheta_\tbu),
    \end{align}
    with static learning-rate matrix $A$ and $\mathcal{F}_{t-1}$-measurable and observable scaling matrix $\scaling$, such that $A\scaling \succ \zeromatrix$.  
\end{definition}  

The scaling $\scaling$ is often based on powers of the Fisher information matrix (e.g., \citealp{creal2013generalized}; \citealp{artemova2022score1}), which is then assumed to be non-singular, i.e.,
\begin{equation}
\label{eqn:Fisher}
\scaling= \Big(\int_{\Ym} s(y,\vartheta_\tbu)s(y,\vartheta_\tbu)^\top f(y|\vartheta_\tbu)\,\,\mathrm{d}y\Big)^{-\zeta}, 
\quad \zeta\in\{0,1/2,1\}.
\end{equation}
The scaling is based only on the postulated (not the true) density. For $\zeta=1/2$, we take  the inverse of the (unique) positive-definite square root of the Fisher information, ensuring $\scaling\succ \zeromatrix$. However, positive definiteness of $A\scaling$ is not automatic even if both $A\succ \zeromatrix$ and $\scaling\succ \zeromatrix$ \citep{nicholson1979eigenvalue}. In fact, the interpretation of SD updates as performing steepest ascent in the metric induced by $A\scaling$ requires this matrix to be positive definite (and therefore symmetric). Interestingly, much of the SD literature ignores this point.\footnote{E.g., \citet[Eq.~2]{creal2013generalized}, \citet[Eq.~6--8]{gasperoni_score-driven_2023}, and \citet[Eq.~6]{d2024modeling} do not impose $A\scaling\succ \zeromatrix$. Symmetrized alternatives such as $\scaling A\scaling$ guarantee positive definiteness if $\scaling, A\succ \zeromatrix$, but deviate from the standard SD setup.} In practice, $A\scaling\succ \zeromatrix$ can be ensured if (i) $A$ and $\scaling$ are diagonal with positive entries, or (ii) one of them is positive definite, while the other is a positive multiple of the identity.

After updating, predictions are constructed as $\vartheta_{t+1|t} = \omega + B \vartheta_{t|t}$, for some vector $\omega$ and matrix $B$.
Combining update and prediction steps yields $\vartheta_{t+1|t} = \omega + B(\vartheta_{t|t-1} + A \scaling \score(y_t,\vartheta_\tbu))$, which up to reparameterization is the standard SD formulation. As our focus is on the update, however, we maintain a clean separation between both steps.

\subsection{Main result: Only SD updates are EKL reducing}
\label{sec:MainResultEKL}

We now establish that only SD updates (and their equivalents) guarantee \emph{expected} Kullback--Leibler (EKL) improvement.  
To this end, we examine the EKL measure
$
\E_{p_t}\!\left[\kl\big(p_t \Vert f_\tpu\big)\right],
$
where $f_\tpu(\cdot) \equiv f(\cdot|\vartheta_{\tpu}(Y_t))$ is evaluated at $Y_t \sim p_t$. More explicitly, recalling \eqref{eqn:EKLDefIntro},
\begin{align}
	\label{eqn:EKLDef}
	\ekl(p_t \Vert f_{\tpu}) 
	:= \int_{\Ym} \int_{\Ym} 
	\log \!\left(\frac{p_t(x)}{f(x | \vartheta_{t|t}(y))} \right) 
	p_t(x)\, p_t(y)\, \dd x \,\dd y,
\end{align}
which we assume to be finite throughout.  
Here the dependence of $\vartheta_{t|t}(y) = \phi(y, \vartheta_\tbu)$ on the observation $y$ is made explicit. The assumed independence of the random variables underlying the observation $y$ (driving the update) and the hypothetical redraw $x$ (used to compute the KL divergence) is reflected in the product $p_t(x)\,p_t(y)$. Hence the EKL divergence~\eqref{eqn:EKLDef} is a ``two-sample'' criterion: one draw updates the parameter, and an independent draw evaluates the fidelity to the true distribution.
EKL improvements therefore measure the expected gain under repeated sampling. By incorporating the uncertainty in both draws, the EKL divergence provides a natural criterion for assessing updating rules.

	\begin{figure}[tb]
		\centering
		\begin{subfigure}[b]{0.45\textwidth}
			\centering
              \includegraphics[width=\linewidth]{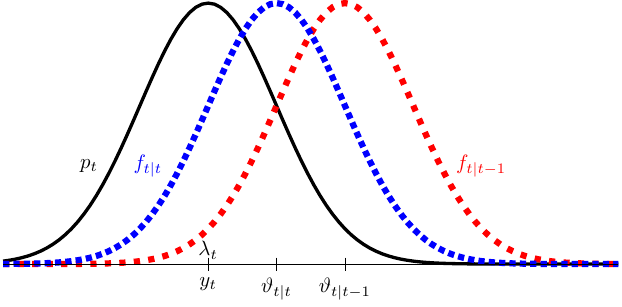}
			\caption{$\lambda_{t} <  \vartheta_\tpu \color{black} <  \vartheta_\tbu$}
			\label{subfig:p<ft+1_mainpaper}
		\end{subfigure}
        \hfill 
		\begin{subfigure}[b]{0.45\textwidth}
			\centering
            \includegraphics[width=\linewidth]{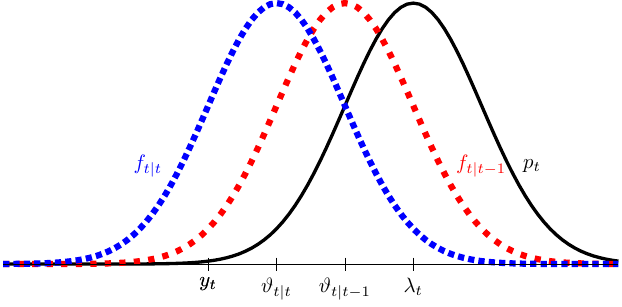}
			\caption{$\lambda_{t} >  \vartheta_\tbu  > \vartheta_\tpu$}
			\label{subfig:p>ft_mainpaper}
		\end{subfigure}
		\caption{SD location model in Example~\ref{exmpl:LKLProblem1} with two hypothetical true densities $p_t$.} 
		\label{fig:IllustrationAlmostSure}
	\end{figure}

\begin{example}
	\label{exmpl:LKLProblem1}
	Consider the true distribution $Y_t \sim p_t = \Nm(\lambda_t,1)$ with unknown time-varying mean $\lambda_t$, and the correctly specified model density $f_\tbu = \Nm(\vartheta_\tbu,1)$. 
    With unit scaling $\scaling=1$ and a (typically small) scalar learning rate $A=\alpha>0$, the SD update $\vartheta_\tpu(y_t) = \vartheta_\tbu + \alpha (y_t - \vartheta_\tbu)$ simply shifts the mean towards the observation $y_t$.  
    Figure~\ref{fig:IllustrationAlmostSure} shows $f_\tbu$ in red, $f_\tpu$ in blue, and $p_t$ in black. 
    In panel~(a), where $y_t$ is typical under $p_t$, the update improves fit; in panel~(b), where $y_t$ is an outlier, it worsens fit. 
    As the setting of panel~(a) is more likely than that of panel~(b), KL improvements may hold \emph{in expectation}.
    Here, the expectation averages out $y_t$ in $\vartheta_{t|t}(y_t)$ with respect to the true law $p_t$, while $x$ represents an independent redraw, also from $p_t$, used to assess distributional closeness in~\eqref{eqn:EKLDef}.
\end{example}

\begin{definition}[EKL difference] 
	\label{def:EKL}
	For any $(\vartheta_\tbu,p_t) \in \Theta \times \Pm$, the \emph{EKL difference} for an updating rule $\phi$ is 
	\begin{align*}
		\Delta^{\ekl}(\phi) 
		\equiv \Delta^{\ekl}(\phi| \vartheta_\tbu,p_t)    
		:=  \mathsf{EKL}(p_{t} \Vert f_\tpu) - \mathsf{EKL}(p_{t} \Vert f_\tbu).  
	\end{align*}
    \sloppy
    We say that an update $\phi$ is \emph{EKL reducing w.r.t.\ the class $\Pm$} if it guarantees an EKL reduction,
    $
    \Delta^{\ekl}(\phi \vert \vartheta_\tbu, p_t) < 0
    $,
    for all $\vartheta_\tbu \in \Theta$ and $p_t \in \Pm$ such that 
    $\E_{p_t}\!\left[ \Delta\phi (Y_t,\vartheta_\tbu) \right]^\top 
	\E_{p_t}\!\left[ \score(X_t,\vartheta_\tbu) \right] \not= 0$.
\end{definition}

The condition $\E_{p_t}\!\left[ \Delta\phi (Y_t,\vartheta_\tbu) \right]^\top 
\E_{p_t}\!\left[ \score(X_t,\vartheta_\tbu) \right] \neq 0$ excludes distributions $p_t$ for which $\vartheta_\tbu$ coincides with the pseudo-true parameter; i.e., a point at which EKL improvements are impossible. The relevance of this condition will become clearer after Theorem~\ref{thm:EKLequivalenceIndepCopy}.

Definition~\ref{def:EKL} involves the difference of two EKL terms.  
Note, however, that $f_\tbu$ does not depend on $Y_t$ and hence
$
\mathsf{EKL}(p_{t}\Vert f_\tbu) =\mathsf{KL}(p_{t} \Vert f_\tbu).
$
In order to analyze the resulting difference, we impose the following assumptions.  
For convenience, the model-implied Hessian is denoted by
\[
\Hessian(x,\vartheta) 
\;:=\; \frac{\partial^2}{\partial\vartheta \,\partial\vartheta^\top} \log f(x | \vartheta).
\]
Our assumptions below ensure that local curvature (in $\vartheta$) of the expected log-likelihood can be bounded; they also justify (the integral form of) exact multivariate mean-value expansions, used in our proofs.

\begin{assumption} 
	\label{ass:thmiffeklindepc}
    \begin{enumerate}[label=(\roman*), itemsep=0.2em, topsep=0.5em]
        \item 
        $\Theta$ is open and convex, and $\vartheta_\tbu \in \Theta$ and $\vartheta_\tpu(Y_t) \in \Theta$.
        
        \item 
        $\log f(x|\vartheta)$ is twice continuously differentiable in $\vartheta$ for all $\vartheta \in \Theta$ and  $x \in \mathcal{Y}$.
        
        \item 
        $\E_{p_t} \big[\score(X_t,\vartheta_\tbu) \big]$ and
        $\E_{p_t} \left[ \Vert \Delta \phi(Y_t, \vartheta_\tbu) \Vert^2\right]$ are finite for all $\vartheta_\tbu \in \Theta$ and  $p_t \in \mathcal{P}$.
    \end{enumerate}
\end{assumption}

\begin{assumptionp}{$\mathcal{HB}$}
    \label{ass:HB}
    There exists a constant $c<\infty$ such that $\sup_{\vartheta \in \Theta} \big\Vert \E_{p_t}[\Hessian (X_t, \vartheta)] \big\Vert \leq c$ for all $p_t \in \Pm$.
\end{assumptionp}

Assumption~\ref{ass:thmiffeklindepc} imposes standard smoothness and moment conditions that are mild and comparable to, or weaker than, those in \citet{blasques2015information}, \citet{Gorgi2023}, and \citet{Creal2024GMM}. 
Our main substantive requirement is the \emph{Hessian boundedness} in expectation in Assumption~\ref{ass:HB}, which can also be written as  $-c \unitmatrix \preceq \E_{p_t}[\Hessian (X_t, \vartheta)] \preceq c \unitmatrix$, uniformly for all $\vartheta \in \Theta$. 
For some model classes (see Section~\ref{sec:Example}), this condition can be shown to hold under mild moment conditions on the class $\Pm$. For this condition to hold independently of $p_t$ (i.e., for all $p_t \in \Pm_0$), a sufficient condition is $\sup_{\vartheta \in \Theta} \Vert \Hessian (x, \vartheta) \Vert < \infty$ for all $x\in\mathcal{Y}$.

For our second main result, we will relax Assumption~\ref{ass:HB} to a localized version (Assumption~\ref{ass:HLB}). However, we note that Assumption~\ref{ass:HB} is \emph{already weaker} than typical conditions in the literature (e.g., \citealp{Gorgi2023};  \citealp{Creal2024GMM}), which additionally demand negative definiteness of the expected Hessian, i.e., 
$\E_{p_t}[\Hessian(X_t,\vartheta)] \prec \zeromatrix$ for all $\vartheta \in \Theta$. This typical condition differs from Assumption~\ref{ass:HB} because it additionally imposes a strict zero upper bound; as a result, it is much harder to achieve and effectively rules out model densities that fail to be log-concave in $\vartheta$ (see Section~\ref{sec:Example} for examples).

To characterize EKL-improving updates, we introduce a localization in the state space via \emph{linearly downscaled} updates. Given an updating rule $\vartheta_\tpu(y) = \phi(y,\vartheta_\tbu)$, we define its downscaled version for a (typically small) value of $\kappa>0$ as
\begin{equation}
\vartheta^\kappa_{t|t}(y) := (1-\kappa)\vartheta_\tbu + \kappa \vartheta_\tpu(y),\quad 0 < \kappa\leq 1.
\end{equation}
Equivalently, 
\begin{equation}
\label{eqn:DefPhilkappa}
\Delta \phi_\kappa(y,\vartheta_\tbu) = \vartheta^\kappa_\tpu(y) - \vartheta_\tbu 
= \kappa \Delta \phi(y,\vartheta_\tbu),
\end{equation}
\sloppy
so adjustments are downscaled by $\kappa$. For SD updates this yields
$
\Delta \phi_{\mathrm{SD},\kappa}(y,\vartheta_\tbu) 
= \kappa A \scaling \score(y,\vartheta_\tbu),
$
where $\kappa$ directly modulates $A \scaling$. For general updating rules, downscaling provides linear control of the step size.

\begin{theorem} 
	\label{thm:EKLequivalenceIndepCopy}
    Consider a class $\Pm$ and let Assumptions~\ref{ass:thmiffeklindepc} and \ref{ass:HB} hold.
    Then, for each $\vartheta_{t\vert t-1}\in \Theta$ and $p_t \in \Pm$ such that $\E_{p_t}\!\left[\Delta \phi(Y_t,\vartheta_\tbu)\right]^\top \E_{p_t}\!\big[\score(X_t,\vartheta_\tbu) \big]\neq 0$, we have:
	\begin{gather*}
		\E_{p_t} \big[ \Delta\phi (Y_t,\vartheta_\tbu) \big]^\top \, \E_{p_t} \big[ \score(X_t,\vartheta_\tbu) \big]  > 0  \\
		\ \iff  \quad 
		\text{ there exists } \bar{\kappa} > 0 \text{ such that for all } \kappa \in (0, \bar{\kappa}] : 
		\ \Delta^{\mathsf{EKL}}(\phi_\kappa| \vartheta_\tbu, p_t) < 0.
	\end{gather*}
\end{theorem}

The proof of Theorem~\ref{thm:EKLequivalenceIndepCopy} applies the path-integral version of the exact (multivariate) mean-value theorem to the EKL criterion at $\vartheta_\tbu$, yielding
\begin{align}
    \label{eqn:EKLProofTechnique}
    \Delta^{\mathsf{EKL}}(\phi_\kappa| \vartheta_\tbu, p_t) 
    = -\kappa \, \E_{p_t}\!\left[ \Delta \phi(Y_t, \vartheta_\tbu) \right]^\top  
      \E_{p_t}\!\left[ \score(X_t, \vartheta_\tbu) \right] 
      + \mathcal{O}(\kappa^2).
\end{align}
Thus, for sufficiently small updates (i.e., $\kappa>0$ small), EKL reductions require the $\mathcal{O}(\kappa)$ term to have the correct sign. 
The admissible update size depends on $p_t$ and $\vartheta_{t\vert t-1}$ and is characterized by an upper bound $\bar{\kappa}\equiv\bar{\kappa}(p_t,\vartheta_{t\vert t-1})$ on the downscaling parameter~$\kappa$.

\sloppy
The condition $\E_{p_t}\left[\Delta \phi(Y_t,\vartheta_\tbu)\right]^\top \E_{p_t} \big[\score(X_t,\vartheta_\tbu) \big] \not= 0$ in Theorem~\ref{thm:EKLequivalenceIndepCopy} assures that the $\mathcal{O}(\kappa)$ term in \eqref{eqn:EKLProofTechnique} does not vanish.
If the inner product of these vectors is zero (e.g., because one of them is zero or they are orthogonal),
then the $\mathcal{O}(\kappa^2)$ term dominates; however, the $\mathcal{O}(\kappa^2)$ term may take either sign, such that no general improvement guarantees are then possible.
We will use a similar approach in Section~\ref{sec:ComparisonOptimality} to analyze the criteria of \citet{Gorgi2023,Creal2024GMM,blasques2015information}.

For a given $\vartheta_\tbu \in \Theta$, Theorem~\ref{thm:EKLequivalenceIndepCopy} establishes EKL improvements for sufficiently small $\kappa > 0$ if and only if the expected update direction and the expected score are aligned, i.e.,
\begin{equation}
    \label{eqn:ExpectedScoreEquivalence}
	\E_{p_t}\!\left[ \Delta\phi (Y_t,\vartheta_\tbu) \right]^\top 
	\E_{p_t}\!\left[ \score(X_t,\vartheta_\tbu) \right] > 0.
\end{equation} 
\sloppy
Updating rules that satisfy this condition for all $\vartheta_\tbu \in \Theta$ and all $p_t \in \Pm$ with the nonzero inner product condition
$\E_{p_t}\!\left[\Delta \phi(Y_t,\vartheta_\tbu)\right]^\top \E_{p_t}\!\big[\score(X_t,\vartheta_\tbu)\big] \neq 0$
(recall Definition~\ref{def:EKL}) are called \emph{score equivalent in expectations} (SEE) with respect to the class $\Pm$.
Because the true density $p_t$ is unknown, we want the class $\mathcal P$ to be as broad as possible. 
The nonzero inner product condition typically excludes any $p_t$ for which $\vartheta_\tbu$ is a stationary point, for then $\E_{p_t}\!\big[\score(X_t,\vartheta_\tbu)\big]=0$.
At such points, the EKL difference is zero to first order in $\kappa$, and the second-order term in~\eqref{eqn:EKLProofTechnique}, whose sign is generally unknown, dominates.

As our EKL criterion is formulated in terms of a double expectation, it is natural that our version~\eqref{eqn:ExpectedScoreEquivalence} of score equivalence also contains two expectations. 
For SD updates, we have $\E_{p_t}[ \Delta \phi_{\mathrm{SD}} (Y_t,\vartheta_\tbu) ] 
= A \scaling \E_{p_t}[ \score(Y_t,\vartheta_\tbu) ]$, such that~\eqref{eqn:ExpectedScoreEquivalence} becomes
\begin{align}
    \label{eqn:SEE_SSDupdates}
    \E_{p_t}[\score(Y_t,\vartheta_\tbu)]^\top\, (A \scaling)\, \E_{p_t}[\score(X_t,\vartheta_\tbu)] \;>\; 0.
\end{align}
For this quadratic form to be positive, we merely require (i) a nonzero expected score (as discussed above), and (ii) $A \scaling\succ O_k$.
Positive definiteness of $A \scaling$ is consistent with the discussion below Definition~\ref{def:NewtonUpdate}.\footnote{Positivity in \eqref{eqn:SEE_SSDupdates} can also hold for non-symmetric $A\scaling$ if (and only if) its \emph{symmetric} part $(A\scaling)+(A\scaling)^\top$ is positive definite. Although non-symmetric versions of $A\scaling$ have been used (e.g., \citealp{gasperoni_score-driven_2023,d2024modeling}), we do not pursue them since they conflict with the SD philosophy that updates move in the direction of the score. Any non-symmetric matrix is an average of its symmetric and anti-symmetric parts: the anti-symmetric part $(A\scaling)-(A\scaling)^\top$, if non-zero, adds an update component $[(A\scaling)-(A\scaling)^\top]\score(y_t,\vartheta_\tbu)$ that is \emph{perpendicular} to the score as
$\score(y_t,\vartheta_\tbu)^\top[(A\scaling)-(A\scaling)^\top]\score(y_t,\vartheta_\tbu)=0$.
Symmetric $A\scaling$ avoids this perpendicular component.}
Theorem~\ref{thm:EKLequivalenceIndepCopy} thus implies  the following result for SD updates.\footnote{The informal arguments in \citet[Eq.~(5)--(7)]{blasques2021finite} motivate a quantity we formalize as the EKL measure, but their simulation-based approach does not provide a theoretical characterization.}

\begin{corollary} 
	\label{cor:EKL_SDupdate}
    Consider a class $\Pm$ and let Assumptions~\ref{ass:thmiffeklindepc} and \ref{ass:HB} hold.
    Then, for any $\vartheta_\tbu \in \Theta$ and $p_t \in \Pm$ with $ \E_{p_t} \big[\score(X_t,\vartheta_\tbu) \big]\neq 0$, there exists a $\bar{\kappa} >0$ such that $\Delta^{\mathsf{EKL}}(\phi_{\mathrm{SD}, \kappa} | \vartheta_\tbu, p_t) < 0$ for all $\kappa \in (0, \bar{\kappa}]$, where $\Delta \phi_{\mathrm{SD}, \kappa} (y_t,\vartheta_\tbu) = \kappa A \scaling \score(y_t,\vartheta_\tbu)$.
\end{corollary}

Corollary~\ref{cor:EKL_SDupdate} shows that SD updates with generic learning-rate and scaling matrices are EKL reducing w.r.t.\ general classes $\Pm$ as long as $A\scaling \succ \zeromatrix$. 
This is in contrast with other methods that tend to require $A \scaling$ to be a scalar multiple of the identity (see Section~\ref{sec:ComparisonOptimality}).
While Corollary~\ref{cor:EKL_SDupdate} invokes Assumption~\ref{ass:HB} for simplicity, its result is valid under the weaker \emph{one-sided} condition, $-c \unitmatrix \preceq \E_{p_t} \big[ \Hessian(X_t,\vartheta) \big]$, uniformly for all $\vartheta \in \Theta$, for some $c<\infty$. 

The SEE condition~\eqref{eqn:ExpectedScoreEquivalence} offers an easily verifiable criterion for whether general updating rules are EKL reducing.
For instance, an update that is usually score-driven but, with some probability less than one, sets $\vartheta_\tpu=\vartheta_\tbu$ remains SEE.  
Similarly, clipped SD updates (Proposition~\ref{prop:CSSD_SEE} below) are SEE if the clipping constant is sufficiently large. 
The Kalman filter, the most widely used algorithm in state-space modeling \citep[Ch.~2]{durbin2012time}, can be expressed as an SD update; hence, it is SEE and its downscaled version is EKL reducing (see Appendix~\ref{sec:KalmanFilter}). Implicit score-driven (ISD, \citealp{lange2022robust}) updates with sufficiently small learning rates are similarly EKL reducing (see Appendix~\ref{app:ISD}).

By contrast, quasi score-driven (QSD, \citealp{blasques2023quasi}) models use two densities: one to compute the score driving the dynamics and another to evaluate the log-likelihood (the model’s ``fit''). Consequently, the SEE condition involves two expectations of distinct scores whose inner product need not be positive. Hence QSD updates with small learning rates are not always EKL reducing (see \appr{app:qsd}).

The SEE condition~\eqref{eqn:ExpectedScoreEquivalence} can be contrasted with the (univariate) \emph{almost sure} score-equivalence condition in \citet{blasques2015information}, i.e., $\operatorname{sign}(\Delta \phi(y_t, \vartheta_\tbu)) = \operatorname{sign}(\score(y_t, \vartheta_\tbu))$ or its multivariate generalization, $\Delta \phi(y_t, \vartheta_\tbu)^\top \score(y_t, \vartheta_\tbu) > 0$. Interestingly, this pointwise condition neither implies nor is implied by the SEE condition~\eqref{eqn:ExpectedScoreEquivalence}, which involves two expectations.

\subsection{Extension to locally bounded Hessians in expectation}
\label{sec:LocallyBoundedHessians}

Assumption~\ref{ass:HB} is violated by several widely used models; for example, 
a Gaussian stochastic volatility model with logarithmic variance (see Section~\ref{sec:Example}). We therefore replace it with the weaker Assumption~\ref{ass:HLB}, which requires the Hessian $H(X_t,\vartheta)$ to be merely \emph{locally bounded} ($\mathcal{HLB}$) in expectation. 

\begin{assumptionp}{$\mathcal{HLB}$}
    \label{ass:HLB}
    For any $\vartheta_\tbu \in \Theta$, there exists a compact ball $\widetilde{\Theta}_r \equiv  \widetilde{\Theta}_r(\vartheta_\tbu) \subset \Theta$ around $\vartheta_\tbu$ with radius $r >0$ 
    such that $\underset{\vartheta \in \widetilde{\Theta}_r(\vartheta_\tbu)}{\sup} \big\Vert \E_{p_t} \big[ \Hessian(X_t,\vartheta) \big] \big\Vert < \infty$ for all $p_t \in \Pm$.
\end{assumptionp}

Assumption~\ref{ass:HLB} is mild and is satisfied by all examples in Section~\ref{sec:Example}: either for any $p_t \in \Pm_0$ or based on mild moment conditions on the class $\Pm$.
For example, for volatility and dependence models based on the Gaussian distribution, the existence of a finite second moment of the data is sufficient for Assumption~\ref{ass:HLB}.
Under Assumption~\ref{ass:HLB}, we obtain a characterization paralleling Theorem~\ref{thm:EKLequivalenceIndepCopy}, now restricted to \emph{bounded} updates.

\begin{theorem} 
	\label{thm:EKLequivalenceBounded}
    Consider a class $\Pm$ and let Assumptions~\ref{ass:thmiffeklindepc} and \ref{ass:HLB} hold.
    Then, for any \emph{bounded update} $\phi$ with $\Vert \Delta\phi (Y_t,\vartheta_\tbu) \Vert \le d(\vartheta_\tbu) < \infty$ and any $\vartheta_\tbu \in \Theta$ and $p_t \in \Pm$ such that $\E_{p_t}\!\left[\Delta \phi(Y_t,\vartheta_\tbu)\right]^\top \E_{p_t}\! \big[\score(X_t,\vartheta_\tbu) \big] \not= 0$, the following equivalence holds:
	\begin{gather*}
		\E_{p_t} \big[ \Delta\phi (Y_t,\vartheta_\tbu) \big]^\top \, \E_{p_t} \big[ \score(X_t,\vartheta_\tbu) \big]  > 0  \\
		\ \iff  \quad 
		\text{ there exists } \bar{\kappa} > 0 \text{ such that for all } \kappa \in (0, \bar{\kappa}] : 
		\ \Delta^{\mathsf{EKL}}(\phi_\kappa| \vartheta_\tbu, p_t) < 0.
	\end{gather*}
\end{theorem}

The focus on bounded updates is consistent with results in optimization and machine learning  \citep{boyd2004convex,nesterov2018lectures,mai2021stability}, which emphasize that for optimization problems with merely locally (but not globally) bounded Hessians, step sizes must be capped.
Applicability of Theorem~\ref{thm:EKLequivalenceBounded} can be guaranteed by a \emph{clipped} SD (CSD) rule, which restricts the update to some distance $d>0$ around $\vartheta_\tbu$,
\begin{align}
    \label{eqn:CSSDupdates}
    \Delta \phi_{\mathrm{CSD}}^d(Y_t,\vartheta_\tbu) 
    := \min \left\{1, \frac{d}{\Vert \Delta\phi_{\mathrm{SD}}(Y_t,\vartheta_\tbu)\Vert} \right\} 
      \Delta\phi_{\mathrm{SD}}(Y_t,\vartheta_\tbu),
\end{align} 
with $1/0:=\infty$, and where the constant $d$ may depend on $\vartheta_\tbu$. Gradient clipping,  standard in stochastic optimization and deep learning \citep{goodfellow2016deep}, improves robustness by capping step sizes. Because clipping is triggered by large scores, it is not obvious that the SEE condition~\eqref{eqn:ExpectedScoreEquivalence} persists. The next proposition shows that, under general conditions on $p_t$, there exists a clipping constant $d>0$ such that CSD updates remain SEE.

\begin{proposition} 
    \label{prop:CSSD_SEE}
    Fix $A\scaling \succ \zeromatrix$ and let Assumption~\ref{ass:thmiffeklindepc} hold. Moreover, let
    (i) $\mu_{p_t} := \E_{p_t} \big[\score(Y_t,\vartheta_\tbu) \big] \neq 0$
    and (ii) $\textnormal{tr}( \Sigma_{p_t})   < \infty$ where
    $\Sigma_{p_t} := \E_{p_t}  \left[\score(Y_t,\vartheta_\tbu)\score(Y_t,\vartheta_\tbu)^\top\right] - \mu_{p_t} \mu_{p_t}^\top$.    If the constant $d>0$ is sufficiently large, such that
    \begin{align}
        \label{eq:probability clipping is active}
        \P_{p_t}\!\left( \big\Vert A \scaling \score(Y_t,\vartheta_\tbu) \big\Vert > d \right)
        \;<\;
       \frac{\|\mu_{p_t}\|_{A\scaling}^2}{\|\mu_{p_t}\|_{A\scaling}^2+\textnormal{tr} \big( A\scaling\,\Sigma_{p_t}\big)},
    \end{align} 
    where $\|\mu_{p_t}\|_{A\scaling}^2:= \mu_{p_t}^\top A\scaling \mu_{p_t}$, then the CSD update~\eqref{eqn:CSSDupdates} satisfies the SEE condition $\E_{p_t} \big[ \Delta \phi_{\mathrm{CSD}}^d (Y_t,\vartheta_\tbu) \big]^\top\, \E_{p_t} \big[ \score(X_t,\vartheta_\tbu) \big]  \; >\; 0$.
\end{proposition}

Inequality~\eqref{eq:probability clipping is active} can be satisfied by choosing the (finite) clipping constant $d$ large enough, with the magnitude dependent on $p_t$. 
For example, if the prediction is inaccurate so that $\|\mu_{p_t}\|_{A {\scaling}}^2$ is large relative to $\mathrm{tr}\!\left(A {\scaling}\Sigma_{p_t}\right)$, condition~\eqref{eq:probability clipping is active} is easily satisfied, as the right-hand side approaches one. Conversely, if the prediction is accurate and $\|\mu_{p_t}\|_{A {\scaling}}^2$ is close to zero, a larger~$d$ may be required to keep clipping sufficiently rare.

Combining Proposition~\ref{prop:CSSD_SEE} with Theorem~\ref{thm:EKLequivalenceBounded}, we find that for any prediction $\vartheta_{\tbu}$ and any true density $p_t$, there exist thresholds $
    \overline{\kappa} \equiv \overline{\kappa}(p_t,\vartheta_{\tbu})>0$ and 
    $\underline{d} \equiv \underline{d}(p_t,\vartheta_{\tbu})<\infty$,
such that for all $\kappa\in(0,\overline{\kappa}]$ and $d \in [\underline{d}, \infty)$, the resulting (C)SD updates improve our EKL criterion. In practice, $\overline{\kappa}$ and $\underline{d}$ are unknown and must be estimated from data (e.g., using empirical gradient moments). Theorem~\ref{thm:UpperBoundLearningRate} in the next section aims to offer practical guidance on how ``large'' SD updates can be.

\begin{remark}
    \label{rem:DiscreteEKL}
    As the proofs reveal, Theorems~\ref{thm:EKLequivalenceIndepCopy} and~\ref{thm:EKLequivalenceBounded} extend to general distributions $P_t \ll \mu$ on a measure space $(\Ym,\mathcal{G},\mu)$. 
    Writing $p_t = \tfrac{\mathrm{d} P_t}{\mathrm{d} \mu}$ as the Radon–Nikodym derivative, all expectations can be interpreted as Lebesgue integrals with respect to $\mu$, which need not be the Lebesgue measure. 
    For discrete distributions, $\mu$ is the counting measure (i.e., integrals reduce to sums), and mixed continuous–discrete cases are also covered. 
    In particular, $p_t$ may be discrete even when $f$ is defined as a density with respect to the Lebesgue measure.
\end{remark}

\begin{remark}
    \label{eqn:PredictionGuarantee}
  We have established guarantees for $\ekl(p_t \Vert f_{\tpu})$ involving the current density $p_t$ based on an observation $Y_t \sim p_t$. We could also seek guarantees for $\ekl(p_{t+1} \Vert f_{\tpu})$, involving the future density $p_{t+1}$, even though $f_{\tpu}$ still depends only on $Y_t \sim p_t$. Adapting the proof of Theorem~\ref{thm:EKLequivalenceIndepCopy} yields the following (necessary and sufficient) condition:
    \begin{align}
        \label{eqn:SEE_SSDupdates_Future}
        \E_{p_t}[\score(Y_t,\vartheta_\tbu)]^\top\, (A \scaling)\, \E_{p_{t+1}}[\score(X_t,\vartheta_\tbu)] \;>\; 0.
    \end{align}
    Because the two expectations are now taken with respect to different laws, this is not a quadratic form; however, \eqref{eqn:SEE_SSDupdates_Future} can hold if \(p_t\) and \(p_{t+1}\) are sufficiently close. For SD models targeting future densities, condition~\eqref{eqn:SEE_SSDupdates_Future} characterizes when this is feasible.
\end{remark}

\section{Upper bounds for the learning-rate matrix}
\label{sec:BoundLearningRate}

So far we asked whether there exists $\bar\kappa>0$ such that any downscaled update $\phi_\kappa$ with $\kappa\in(0,\bar\kappa]$ improves the EKL criterion. 
Here, for a given true density $p_t$, we instead quantify how large the SD matrix combination $\mathcal{A}_{t-1}:=A\scaling$ can be while still guaranteeing EKL improvement. 
If $\mathcal{A}_{t-1}$ is positive definite and sufficiently small, Corollary~\ref{cor:EKL_SDupdate} applies. 
The next result provides explicit upper bounds on the elements (or eigenvalues) of $\mathcal{A}_{t-1}$. 
When $\scaling=\unitmatrix$, these are bounds on the static learning-rate matrix $A$; more generally, they bound the combined matrix $\mathcal{A}_{t-1}$ multiplying the score, beyond the specific case $\mathcal{A}_{t-1}=A\scaling$.

\begin{theorem} 
	\label{thm:UpperBoundLearningRate} 
    Consider a class $\Pm$ and let Assumptions~\ref{ass:thmiffeklindepc} and \ref{ass:HB} hold.
    For a given $\vartheta_{\tbu} \in \Theta$, let $p_t \in \Pm$ be such that $\mu_{p_t} := \E_{p_t} \big[\score(Y_t,\vartheta_\tbu) \big]\neq 0$, with elements $\mu_1,\dots,\mu_k$.
    Moreover, let
    $\Sigma_{p_t} := \E_{p_t}  \left[\score(Y_t,\vartheta_\tbu)\score(Y_t,\vartheta_\tbu)^\top\right] - \mu_{p_t} \mu_{p_t}^\top$ be finite with diagonal elements $\sigma_1^2, \dots,\sigma_k^2$.
    \sloppy 
    Then, the condition $\Delta^{\mathsf{EKL}}(\phi_{\mathrm{SD}} | \vartheta_\tbu, p_t) < 0$ holds either
    \color{black}
    \begin{enumerate}[label=(\Alph*), itemsep=0.2em, topsep=0.5em]
        \item  
        for SD updates with $\mathcal{A}_{t-1} = \alpha_{t-1} \unitmatrix \succ \zeromatrix$ if
        \begin{align}
            \label{eqn:LowerBoundScalar}
            \alpha_{t-1} < 
            \frac{2}{c} \,  \frac{\|\mu_{p_t}\|^2 }{\|\mu_{p_t}\|^2+ \trace(\Sigma_{p_t})}; 
            \qquad \text{or}
        \end{align}
    
        \item     
        for SD updates with a diagonal $\mathcal{A}_{t-1} = \operatorname{diag}(\alpha_{1, t-1}, \dots, \alpha_{k, t-1}) \succ \zeromatrix$ and non-zero elements $\mu_1,\dots,\mu_k$ if
        \begin{align}
            \label{eqn:LowerBoundDiagonal}
            \alpha_{i,t-1} \le \frac{2}{c} \, \frac{\mu_{i}^2}{\mu_{i}^2+\sigma_{i}^2}, \qquad  i=1,\ldots,k,
        \end{align}
        and the inequality~\eqref{eqn:LowerBoundDiagonal} is strict for at least one $i = 1,\dots, k$; or
            \item 
        for SD updates with a general $\mathcal{A}_{t-1} \succ \zeromatrix$ if
        \begin{align}
            \label{eqn:LowerBoundPosDef}
            \frac{\lambda_{\max}(\mathcal{A}_{t-1})^2}{\lambda_{\min}(\mathcal{A}_{t-1})}
            < \frac{2}{c} \, \frac{\Vert \mu_{p_t}\Vert^2}{\Vert \mu_{p_t} \Vert^2 + \trace(\Sigma_{p_t})}.
        \end{align}
    \end{enumerate}
\end{theorem}

The three formulations (A)–(C) of Theorem~\ref{thm:UpperBoundLearningRate} share the factor $2/c$ and the same ``signal-to-noise'' ratio based on its first two moments of the score, analyzed in detail below for the univariate case. In multivariate settings, (A) can be overly conservative, while (B) better adapts to heterogeneity in the cross-section. Its required score moments in each direction can be estimated from the data.
E.g., if $p_t$ conditions on $\Fm_{t-1} = \sigma(Y_s: s \le t-1)$, one can use exponentially weighted moving averages, as is standard in the optimization literature (e.g., \citealp{kingma2014adam}). 
Part (C) allows a generic $\mathcal{A}_{t-1}\succ \zeromatrix$, with EKL improvements guaranteed if $\lambda_{\max}(\mathcal{A}_{t-1})\,\mathrm{cond}(\mathcal{A}_{t-1})$ is upper bounded, where $\mathrm{cond}(\cdot)=\lambda_{\max}(\cdot)/\lambda_{\min}(\cdot)$ is the condition number. Hence EKL improvements can be guaranteed for general learning-rate matrices $\mathcal{A}_{t-1}\succ O_k$ if we control both the largest eigenvalue and the condition number.

Parts (B) and (C) are specific to our framework: they guide admissible choices of the combined matrix
$\mathcal{A}_{t-1}=A\scaling$, allowing both diagonal $\mathcal{A}_{t-1}$ with heterogeneous entries and general positive definite $\mathcal{A}_{t-1}$. 
In contrast, the methods of \cite{Gorgi2023} (Sections~\ref{sec:CEV}--\ref{sec:MSE}) and \cite{Creal2024GMM} (Section~\ref{sec:EGMM}) yield improvements only when $\mathcal{A}_{t-1}$ is proportional to the identity.

To facilitate the comparison with \cite{Gorgi2023}, we take $k=1$ and $\scaling=1$.
Theorem~\ref{thm:UpperBoundLearningRate} guarantees EKL improvements for any SD update with learning rate $\alpha>0$ if
\begin{align}
    \label{eqn:AlphaBarEKL}
    \alpha < \bar{\alpha}_\ekl
	:= \frac{2}{c} \;\frac{\left(\E_{p_t}[\score(Y_t,\vartheta_\tbu)]\right)^2}{\left(\E_{p_t}[\score(Y_t,\vartheta_\tbu)]\right)^2 + \var_{p_t}\!\left(\score(Y_t,\vartheta_\tbu)\right)}.
\end{align}
The bound splits into a static term $2/c$ and a dynamic ``signal-to-noise'' ratio formed from moments of the score. This ratio is always at most one, so $\bar{\alpha}_{\ekl}\le 2/c$. Both moments depend on the prediction $\vartheta_{\tbu}$. If $\vartheta_{\tbu}$ is close to the pseudo-true parameter
\begin{align}
    \label{eqn:PseudoTrueParam}
    \vartheta_t^\ast := \arg \max_{\vartheta \in \Theta} \E_{p_t}\big[\log f(Y_t | \vartheta)\big],
\end{align}
and, for simplicity, we assume $\vartheta_t^\ast \in \Theta$, then standard regularity conditions imply
$\E_{p_t}\!\left[\score(Y_t, \vartheta_t^\ast)\right] = 0$, so that $\bar{\alpha}_{\ekl}$ approaches zero as $\vartheta_{\tbu}\to \vartheta_t^\ast$.
In other words, the more accurate the prediction, the smaller the admissible learning rate.

\citet[Thm.~1]{Gorgi2023} derive a simpler (and larger) bound in the univariate case. If $-c \le \E_{p_t} \big[ \Hessian(X_t, \vartheta) \big] < 0$ for all  $\vartheta \in \Theta$, then any $\alpha < \bar{\alpha}_\cev := 2/c$ ensures improvement under their conditional expected variation (CEV) criterion, 
\allowdisplaybreaks
\begin{alignat}{2}
	\begin{aligned}
		\label{eqn:CEV_definition} 
		\big\vert \vartheta_t^\ast - \E_{p_t}[\vartheta_\tpu(Y_t)] \big\vert &< \big\vert \vartheta_t^\ast -  \vartheta_\tbu \big\vert 
		\qquad  &&\text{if } \vartheta_\tbu \not=  \vartheta_t^\ast, \\
		\E_{p_t}[\vartheta_\tpu(Y_t)] &= \vartheta_t^\ast 
		\qquad  &&\text{if } \vartheta_\tbu =  \vartheta_t^\ast.
	\end{aligned}
\end{alignat}
The CEV criterion depends solely on the distance of $\E_{p_t}[\vartheta_\tpu(Y_t)]$ from $\vartheta_t^\ast$, disregarding further distributional differences between $f_{\tpu}$ and $p_t$. We note that $\bar{\alpha}_\ekl\leq 2/c= \bar{\alpha}_\cev$. Indeed, the EKL criterion demands decreasing learning rates as predictions become more accurate, but $\bar{\alpha}_\cev$ remains fixed at $2/c$. Hence CEV may favor sizable learning rates even when the gradient is  largely dominated by noise. 
The EKL criterion provides the sharper and more intuitive requirement that learning rates should shrink as predictions improve.

\section{Comparison to related performance criteria}
\label{sec:ComparisonOptimality}

This section contrasts our EKL-based characterization of SD updates with four recent alternative performance measures. Sections~\ref{sec:CEV}--\ref{sec:EGMM} show that these approaches require stronger conditions (e.g., concave model log-densities) yet deliver weaker results: (i) their equivalence conditions are not directly useful for model construction, and (ii) non-SD updates may improve these measures, while SD updates are guaranteed to do so only under more restrictive assumptions. The performance measure discussed in Section~\ref{sec:TKLComparison} is \emph{improper} in the statistical sense \citep{gneiting2007strictly}: no conclusions on its basis can be drawn.

\subsection{CEV of \citet{Gorgi2023}}
\label{sec:CEV}

\citet{Gorgi2023} evaluate SD models using the univariate CEV criterion~\eqref{eqn:CEV_definition}. Our multivariate generalization of their CEV criterion can be written as
\begin{align}
	\label{eqn:CEVDef}
	\cev(p_t \Vert f_{\tpu}) 
    :=
    \Big( \E_{p_t} \big[ \vartheta_\tpu (Y_t) \big] -  \vartheta_t^\ast \Big)^\top \Omega_{t-1}  \Big( \E_{p_t} \big[ \vartheta_\tpu (Y_t) \big] -  \vartheta_t^\ast \Big),
\end{align}
where $\Omega_{t-1}\succ \zeromatrix$ is an $\mathcal{F}_{t-1}$-measurable weighting matrix. 
In the univariate case, \cite{Gorgi2023} show that SD updates with sufficiently small learning rates reduce CEV losses, but leave open the question whether this property is unique to SD updates. Proposition~\ref{prop:CEVEquiv} below extends their analysis to the multivariate case, characterizing all CEV-improving updates by examining the CEV difference, defined as
\begin{align*}
	\Delta^{\cev}(\phi) 
	\equiv \Delta^{\cev}(\phi| \vartheta_\tbu,p_t)    
	:= \cev(p_{t} \Vert f_\tpu) - \cev(p_{t} \Vert f_\tbu).
\end{align*}
An update $\phi$ is \emph{CEV reducing w.r.t.\ $\Pm$} if $ \Delta^{\cev}(\phi \vert \vartheta_\tbu, p_t) < 0$ for all $\vartheta_\tbu \in \Theta, p_t \in \Pm$ such that $\E_{p_t} \big[ \Delta\phi (Y_t,\vartheta_\tbu) \big]^\top \Omega_{t-1} \big(\vartheta_t^\ast -  \vartheta_\tbu \big) \not=0$.
This mirrors the exclusion in Definition~\ref{def:EKL} and prohibits $\vartheta_\tbu$ from being equal to $\vartheta_t^\ast$.
The class $\Pm$ may differ from that in Theorem~\ref{thm:EKLequivalenceIndepCopy}.

\begin{proposition}
    \label{prop:CEVEquiv}
    Consider a class $\Pm$ and let Assumption~\ref{ass:thmiffeklindepc} hold.
    Then, for $\phi_\kappa(y, \vartheta_\tbu)$ in~\eqref{eqn:DefPhilkappa}, the following equivalence holds for all $p_t \in \Pm$, $\vartheta_\tbu \in \Theta$:
    \color{black}
	\begin{gather}
		\E_{p_t} \big[ \Delta\phi (Y_t,\vartheta_\tbu) \big]^\top \Omega_{t-1} \big(\vartheta_t^\ast -  \vartheta_\tbu \big) > 0  
           \label{eqn:CEVcondition}
           \\
		\ \iff  \quad 
		\text{ there exists } \bar{\kappa} > 0 \text{ such that for all } \kappa \in (0, \bar{\kappa}] : 
		\ \Delta^{\mathsf{CEV}}(\phi_\kappa| \vartheta_\tbu, p_t) < 0.\notag
	\end{gather}
\end{proposition}

Proposition~\ref{prop:CEVEquiv} shows that an updating rule $\phi$ is CEV reducing if and only if~\eqref{eqn:CEVcondition} holds: on average, $\phi$ must move the parameter toward the pseudo-truth $\vartheta_t^\ast$, with admissible directions governed by $\Omega_{t-1}\succ\zeromatrix$. Since $\vartheta_t^\ast$ is unknown, however, this condition is not directly usable for model construction. Moreover, it is unclear whether~\eqref{eqn:CEVcondition} is satisfied only by SD updates.

Next, we derive minimal conditions under which SD updates
$\Delta \phi_{\mathrm{SD}}(Y_t,\vartheta_\tbu)=A \scaling \score(Y_t,\vartheta_\tbu)$
guarantee CEV improvements. Using $\E_{p_t}\!\left[\score(Y_t,\vartheta_t^\ast)\right]=0$,
which holds under standard regularity conditions (see \citealp[Ass.~1]{Gorgi2023}),
condition~\eqref{eqn:CEVcondition} for CEV reductions becomes
\allowdisplaybreaks
\begin{align}
    &\hspace{-0.5cm}0<\E_{p_t} \left[ \Delta \phi_{\mathrm{SD}}(Y_t,\vartheta_\tbu)\right]^\top \Omega_{t-1} \big( \vartheta_t^\ast-\vartheta_\tbu  \big) \nonumber \\
    &= \E_{p_t} \left[ \score(Y_t,\vartheta_\tbu)^\top - \score(Y_t,\vartheta^\ast_t)^\top \right] (A \scaling) \Omega_{t-1} \big( \vartheta_t^\ast-\vartheta_\tbu  \big) \nonumber \\
    &=\big( \vartheta_t^\ast-\vartheta_\tbu  \big)^\top  \left( \E_{p_t} \left[ - \int_0^1 \Hessian \big( Y_t, \vartheta_t^u \big) \mathrm{d}u \right] (A\scaling) \Omega_{t-1} \right) \big( \vartheta_t^\ast-\vartheta_\tbu  \big),
    \label{eqn:CEVcondition2}    
\end{align}
where we used the integral form of the mean-value theorem with ${\vartheta}_t^u:= u \vartheta_\tbu + (1-u)\vartheta_t^\ast$ for $u\in[0,1]$. To ensure the positivity of the quadratic form~\eqref{eqn:CEVcondition2} for all $\vartheta_t^\ast-\vartheta_\tbu \neq 0$, \cite{Gorgi2023} impose that the central term (in large parentheses) is positive. 
In their setting, this term is a scalar; since $A \scaling \Omega_{t-1} > 0$, it then suffices to assume
$\E_{p_t}\left[-\Hessian\big(Y_t,\vartheta\big)\right] > 0$ for all $\vartheta$.
In the multidimensional case, as considered here, the natural analogue is to require the \emph{symmetric} part of central matrix
$
\E_{p_t}\!\left[-\Hessian\big(Y_t,\vartheta\big)\right]A \scaling \Omega_{t-1}
$
to be positive definite:
\begin{align}
   \E_{p_t}[-\Hessian(Y_t,\vartheta)] (A \scaling) \Omega_{t-1} 
   + \Omega_{t-1} (A \scaling) \E_{p_t}[-\Hessian(Y_t,\vartheta)] \;\succ\; O_k, 
   \quad \forall \vartheta \in \Theta.
   \label{eqn:CEVcondition4}    
\end{align}
However, ensuring that condition~\eqref{eqn:CEVcondition4} holds is difficult, as the symmetrized sum need not be positive definite even if all individual matrices are \citep{nicholson1979eigenvalue}.
The following assumption provides a sufficient condition that guarantees \eqref{eqn:CEVcondition4}.

\begin{assumptionp}{$\mathcal{HN}$}
    \label{ass:HN}
    (i) $\E_{p_t} \big[ \Hessian ( Y_t, \vartheta ) \big] \prec \zeromatrix$ for all $\vartheta\in \Theta$, $p_t \in \Pm$; and
    (ii) $A \scaling \Omega_{t-1} = d_{t-1} \unitmatrix$ for some $\mathcal{F}_{t-1}$-measurable $d_{t-1} > 0$.
\end{assumptionp}

Part~(ii) of Assumption~\ref{ass:HN} limits the generality of $A\scaling$ in the multivariate setting: in the canonical case $\Omega_{t-1}=\unitmatrix$, it effectively enforces a scalar learning rate and scaling function. Part~(i) is also hard to ensure, even in one dimension, unless $\Hessian(y,\vartheta)$ is negative definite for all $y\in\mathcal{Y}$ and $\vartheta\in\Theta$, i.e., the model density is strictly log-concave. By contrast, the EKL-specific Assumption~\ref{ass:HB} only requires boundedness of the expected Hessian, which often follows from mild moment conditions on $p_t$; see Section~\ref{sec:Example}.

Assumption~\ref{ass:HN}~(i) reduces to $\E_{p_t}[\Hessian(Y_t,\vartheta)]<0$ in the scalar case. Relative to \citet[Ass.~3]{Gorgi2023}, this is slightly more general because it does not impose a lower bound. By dropping the lower bound, we already cover the localized theory of \citet[Thm.~2]{Gorgi2023}. An analogue of \eqref{eqn:CEVcondition2} for general scaling matrices is given in Appendix~\ref{sec:ScaledHessianCEVGeneral}.
\citet[Thm.~3]{Gorgi2023} relax the zero upper bound in Assumption~\ref{ass:HN}~(i) further, but only for small deviations $\|\vartheta_\tbu-\vartheta_t^\ast\|$, which cannot be ensured in practice. 

In sum, \eqref{eqn:CEVcondition2} and \eqref{eqn:CEVcondition4} show that CEV reductions are guaranteed under Assumption~\ref{ass:HN} but are otherwise hard to establish. Compared with Assumption~\ref{ass:HB}, which only requires a bounded expected Hessian, Assumption~\ref{ass:HN} imposes a zero upper bound (negative definiteness of the expected Hessian), excluding, for example, \citeauthor{harvey2014filtering}'s (\citeyear{harvey2014filtering}) Student's~$t$ location model (Section~\ref{sec:ExampleUnivariate}) and a bivariate Gaussian location-scale model (Section~\ref{sec:ExampleBivaraiteGaussian}). The difference stems from the proof methods: for the EKL criterion, the Hessian appears in the $\mathcal{O}(\kappa^2)$ term of Theorem~\ref{thm:EKLequivalenceIndepCopy} and only needs to be bounded, whereas for the CEV criterion it enters via \eqref{eqn:CEVcondition2} in the $\mathcal{O}(\kappa)$ term of Proposition~\ref{prop:CEVEquiv}, requiring the stronger sign restriction.

Thus, we conclude that while non-SD rules may yield CEV reductions, generic SD updates with $A\scaling \succ \zeromatrix$ are not necessarily guaranteed to consistently do so.

\subsection{MSE of \citet{Gorgi2023}}
\label{sec:MSE}

\citet[Corol.~1]{Gorgi2023} study updates that reduce the mean squared error (MSE) relative to the pseudo-true parameter $\vartheta_t^\ast$ in~\eqref{eqn:PseudoTrueParam}. A multivariate  extension of their criterion~is
\begin{align}
	\label{eqn:MSEDef}
	\mse(p_t \Vert f_{\tpu}) 
    := \E_{Y_t \sim p_t} \!\left[ 
        \big(\vartheta_\tpu(Y_t) - \vartheta_t^\ast\big)^\top 
        \Omega_{t-1} 
        \big(\vartheta_\tpu(Y_t) - \vartheta_t^\ast\big) 
    \right],
\end{align}
where, as before, $\Omega_{t-1}\succ O_k$ is an $\mathcal{F}_{t-1}$-measurable weighting matrix. 
The dependence of $\mse(p_t \Vert f_{\tpu})$ on $\vartheta_\tpu(Y_t)$ is captured through $f_\tpu$ as $f_\tpu(\cdot) = f(\cdot|\vartheta_\tpu(Y_t))$.
We define 
\begin{align*}
	\Delta^{\mse}(\phi) 
	\equiv \Delta^{\mse}(\phi| \vartheta_\tbu,p_t) 
	:= \mse(p_{t} \Vert f_\tpu) - \mse(p_{t} \Vert f_\tbu),
\end{align*}
and an update $\phi$ is said to be \emph{MSE reducing w.r.t.\ $\Pm$} if $ \Delta^{\mse}(\phi \vert \vartheta_\tbu, p_t) < 0$, for all $\vartheta_\tbu \in \Theta$ and $p_t \in \Pm$ such that $\E_{p_t} \big[ \Delta\phi (Y_t,\vartheta_\tbu) \big]^\top \Omega_{t-1} \big(\vartheta_t^\ast -  \vartheta_\tbu \big) \not=0$, where the last condition prohibits $\vartheta_\tbu$ from being identical to $\vartheta_t^\ast$ (similar to the condition in Definition~\ref{def:EKL}).

\citet[Corol.~1]{Gorgi2023} show that SD updates with sufficiently small learning rates reduce MSE losses, but leave open whether this property is unique to SD updates---a question addressed by our next result, which characterizes all MSE-reducing updates.

\begin{proposition}
    \label{prop:MSEEquiv}
    Consider a class $\Pm$ and let Assumption~\ref{ass:thmiffeklindepc} hold.
    Then, for $\phi_\kappa(y, \vartheta_\tbu)$ in~\eqref{eqn:DefPhilkappa}, the following equivalence holds for all $p_t \in \Pm$, $\vartheta_\tbu \in \Theta$:
    \color{black}
	\begin{gather}
    \label{eqn:MSEcondition}
		\E_{p_t} \big[ \Delta\phi (Y_t,\vartheta_\tbu) \big]^\top \Omega_{t-1} \big(\vartheta_t^\ast -  \vartheta_\tbu \big) > 0  \\
		\ \iff  \quad 
		\text{ there exists } \bar{\kappa} > 0 \text{ such that for all } \kappa \in (0, \bar{\kappa}] : 
		\ \Delta^{\mathsf{MSE}}(\phi_\kappa| \vartheta_\tbu, p_t) < 0.
        \notag
	\end{gather}
\end{proposition}

Somewhat surprisingly, the MSE improvement condition~\eqref{eqn:MSEcondition} coincides with the CEV condition~\eqref{eqn:CEVcondition}, since both criteria share the same $\mathcal{O}(\kappa)$ term (as shown in the proofs). As discussed in Section~\ref{sec:CEV}, this equivalence condition is not directly useful for model construction. Moreover, MSE-improvement guarantees for SD updates require the strong Assumption~\ref{ass:HN}, which imposes a negative-definite expected Hessian and effectively forces $A\scaling$ to be a scalar multiple of $\Omega_{t-1}^{-1}$, a substantial restriction in higher dimensions.

\subsection{EGMM of \citet{Creal2024GMM}}
\label{sec:EGMM}

\citet{Creal2024GMM} propose to update time-varying parameters using a GMM-type influence function based on an expected GMM objective, which resembles an SD update when the moment condition uses the score. In our notation, their proposed updating rule is:\footnote{The last equation on p.~3 of \citet{Creal2024GMM} is missing a minus sign.}
\begin{align}
    \label{eqn:EGMM_SDupdate}
    \Delta \phi_{\mathrm{GMM}}(y_t, \vartheta_\tbu) = A \,
    \big(-\E_{p_t}  \big[\Hessian(X_t,\vartheta_\tbu) \big] \big)^{-1} \,
    \score(y_t, \vartheta_\tbu).
\end{align}
This update is motivated in \cite{Creal2024GMM} as being similar to an SD step with \emph{inverse Fisher scaling}. However, the proposed scaling matrix
\begin{align}
    \label{eqn:EGMM_InfeasibleScaling}
    \scaling^{\mathrm{GMM}} = \Big( -\E_{p_t}\big[\Hessian(X_t,\vartheta_\tbu)\big] \Big)^{-1}
\end{align}
is infeasible in practice due to the dependence on the unknown true density $p_t$ and hence does not follow from the standard Fisher scaling formula~\eqref{eqn:Fisher}.
In general, $\E_{p_t}\!\big[-\Hessian(X_t,\vartheta_\tbu)\big]$ should \emph{not} be interpreted as an information matrix (e.g., \citealp[p.~63]{van1998asymptotic}). 
An information matrix arises only under correct specification; namely, when (i) the model density $f(\cdot| \vartheta_\tbu)$ has the same functional form as the true density $p_t = p(\cdot| \lambda_t)$ for some true parameter $\lambda_t$, and (ii) the parameters coincide, $\vartheta_\tbu = \lambda_t$. Neither condition is guaranteed in~\eqref{eqn:EGMM_InfeasibleScaling}. Consequently, $\E_{p_t}\![-\Hessian(X_t,\vartheta_\tbu)]$ may even be indefinite, so that gradient updates can move in directions that are poorly aligned with (or even opposite to) the score. 
This phenomenon can occur, for instance, in the Student's $t$ location model and the Gaussian location-scale model, discussed in Appendix~\ref{sec:DetailsExamples} and Section~\ref{sec:ExampleBivaraiteGaussian}, respectively.

\citet[Prop.~3]{Creal2024GMM} show that infinitesimal updates of the form~\eqref{eqn:EGMM_SDupdate} improve their objective. 
However, their result relies on the infeasible scaling in \eqref{eqn:EGMM_InfeasibleScaling} and does hence not apply to feasible SD models.
To investigate whether feasible SD models (and only SD models) are EGMM reducing, we consider general downscaled updates $\phi_\kappa$ as defined in~\eqref{eqn:DefPhilkappa}. Specifically, when the score is used as the moment condition, the expected GMM (EGMM) objective in \citet[Eq.~18--19]{Creal2024GMM} takes the form:
\begin{align}
    \begin{aligned} 
    	\label{eqn:GMMCriterion}
    	&\gmm \big(p_t \Vert f_\tpu \big) = \E_{Y_t \sim p_t} \left[  \E_{X_t \sim p_t} \left[\score(X_t, \vartheta_{t\vert t}(Y_t)) \right]^\top  \Omega_{t-1} \E_{X_t \sim p_t} \left[\score(X_t, \vartheta_{t\vert t}(Y_t)) \right] \right],
    \end{aligned}    
\end{align}
where, as before, $\Omega_{t-1}\succ O_k$ is an $\mathcal{F}_{t-1}$-measurable weighting matrix. Proposition~\ref{prop:EGMMEquiv} below extends the analysis in \cite{Creal2024GMM} by fully characterizing which updates $\phi$ guarantee EGMM improvements in the sense that the loss difference $\Delta^{\gmm}(\phi)$ is negative:
\begin{align}
    \label{eqn:EGMM_Improvements}
	\Delta^{\gmm}(\phi) 
	\equiv \Delta^{\gmm}(\phi| \vartheta_\tbu,p_t)
	:= \gmm(p_{t} \Vert f_\tpu) - \gmm(p_{t} \Vert f_\tbu).
\end{align}
An update $\phi$ is \emph{EGMM reducing w.r.t.\ $\Pm$} if $ \Delta^{\gmm}(\phi \vert \vartheta_\tbu, p_t) < 0$ for all $\vartheta_\tbu \in \Theta$ and $p_t \in \Pm$ such that 
$\E_{p_t} \left[\Delta \phi(X_t, \vartheta_{\tbu}) \right]^\top  \E_{p_t} \left[ - \Hessian(X_t, \vartheta_\tbu) \right] \Omega_{t-1} \E_{p_t} \left[\score(X_t, \vartheta_{\tbu}) \right] \not= 0$. 
Here, the last condition is again analogous to (though distinct from) the one in Definition~\ref{def:EKL} and typically prohibits that the prediction $\vartheta_{\tbu}$ coincides with the pseudo-true parameter.

\begin{proposition}
    \label{prop:EGMMEquiv}
    Consider a class $\Pm$, let Assumptions~\ref{ass:thmiffeklindepc} and \ref{ass:HB} hold and let 
    $\sup_{\vartheta \in \Theta} \big\Vert (\partial/\partial \vartheta_j) \, \E_{p_t} \big[ \Hessian(X_t, \vartheta) \big] \big\Vert < \infty$ for all $j=1,\dots,k$.
    Then, for $\phi_\kappa(y, \vartheta_\tbu)$ given in \eqref{eqn:DefPhilkappa}, the following equivalence holds  for each $\vartheta_\tbu \in \Theta$ and $p_t \in \Pm$, for which
    $\E_{p_t} \left[\Delta \phi(X_t, \vartheta_{\tbu}) \right]^\top  \E_{p_t} \left[ - \Hessian(X_t, \vartheta_\tbu) \right] \Omega_{t-1} \E_{p_t} \left[\score(X_t, \vartheta_{\tbu}) \right] \not= 0$:
    \color{black}    
	\begin{gather}
		\E_{p_t} \left[\Delta \phi(X_t, \vartheta_{\tbu}) \right]^\top  \E_{p_t} \left[ - \Hessian(X_t, \vartheta_\tbu) \right] \Omega_{t-1} \E_{p_t} \left[\score(X_t, \vartheta_{\tbu}) \right] > 0 \label{eqn:EGMM equiv condition} 
        \\
		\ \iff  \quad 
		\text{ there exists } \bar{\kappa} > 0 \text{ such that for all } \kappa \in (0, \bar{\kappa}] : 
		\ \Delta^{\mathsf{EGMM}}(\phi_\kappa| \vartheta_\tbu, p_t) < 0.
        \notag
	\end{gather}
\end{proposition}
The assumptions of Proposition~\ref{prop:EGMMEquiv} closely parallel those of Proposition~3 in \cite{Creal2024GMM}. Our proof builds on arguments in their online supplement and incorporates a path-integral version of the exact multivariate mean-value theorem, yielding
\begin{align*}
&\Delta^{\gmm}(\phi_\kappa)\notag
= - 2 \kappa  \E_{p_t} \left[\Delta \phi(X_t, \vartheta_{\tbu}) \right]^\top
\E_{p_t} \!\left[ - \Hessian(X_t, \vartheta_\tbu) \right]
\Omega_{t-1} \E_{p_t} \!\left[\score(X_t, \vartheta_{\tbu}) \right]
+ \mathcal{O}(\kappa^2).
\end{align*}
By analyzing the sign of the $\mathcal{O}(\kappa)$ term, which determines the sign of $\Delta^{\gmm}(\phi_\kappa)$ for sufficiently small $\kappa$, we obtain the equivalence condition~\eqref{eqn:EGMM equiv condition} in Proposition~\ref{prop:EGMMEquiv}. However, \eqref{eqn:EGMM equiv condition} is hard to interpret because it involves a nontrivial matrix product with $\E_{p_t}\!\left[-\Hessian(X_t,\vartheta_\tbu)\right]$, so Proposition~\ref{prop:EGMMEquiv} provides limited guidance for model construction.

To derive practically verifiable conditions, we next examine minimal conditions under which SD updates guarantee EGMM improvements. Substituting $\Delta \phi_{\mathrm{SD}}(Y_t,\vartheta_\tbu) = A \scaling \score(Y_t,\vartheta_\tbu)$ into~\eqref{eqn:EGMM equiv condition}, we get
\begin{equation}
    \label{eqn:EGMMcondition1}
    \E_{p_t} \!\left[\score(X_t, \vartheta_{\tbu}) \right]^\top \Big(A 
    \scaling
    \E_{p_t} \!\left[ - \Hessian(X_t, \vartheta_\tbu) \right]
    \Omega_{t-1} \Big)
    \E_{p_t} \!\left[\score(X_t, \vartheta_{\tbu}) \right] > 0.
\end{equation}
For this quadratic form to be positive for all $\E_{p_t}[\score(X_t,\vartheta_\tbu)] \neq 0$, the symmetric part of the central matrix (in large brackets) must be positive definite for each $\vartheta_\tbu$, i.e., we need
\begin{equation}
    \label{eqn:EGMMcondition2}
    (A \scaling)
    \E_{p_t} \!\left[ - \Hessian(X_t, \vartheta) \right] \Omega_{t-1}
    + \Omega_{t-1} \E_{p_t}\!\left[ - \Hessian(X_t, \vartheta) \right] (A \scaling)
    \succ \zeromatrix, \quad \forall \vartheta \in \Theta.
\end{equation}

\citet[Eq.~A.5]{Creal2024GMM} consider the infeasible scaling~\eqref{eqn:EGMM_InfeasibleScaling}, which is such that the expected Hessian in condition~\eqref{eqn:EGMMcondition2} drops out. For the resulting expression to be positive definite, $A$ should still be a positive scalar multiple of the identity.

However, for general feasible scalings (i.e., not depending on $p_t$), the expected Hessian in condition~\eqref{eqn:EGMMcondition2}  cannot be removed. A sufficient condition to guarantee condition~\eqref{eqn:EGMMcondition2} is to impose $\E_{p_t}\!\left[\Hessian(X_t,\vartheta)\right]\prec\zeromatrix$ and take $A\scaling$ as a positive scalar multiple of $\Omega_{t-1}\succ O_k$. Then \eqref{eqn:EGMMcondition2} reduces to
$2\,\Omega_{t-1}\,\E_{p_t}\!\left[-\Hessian(X_t,\vartheta)\right]\Omega_{t-1}\succ\zeromatrix$, which holds due to symmetry. To guarantee that SD updates yield EGMM improvements, we thus impose the following conditions, combining those of Proposition~\ref{prop:EGMMEquiv} with those ensuring that \eqref{eqn:EGMMcondition2} holds.

\begin{assumptionp}{$\mathcal{HBNT}$}
    \label{ass:HBNT}
    \sloppy
    For all $p_t \in \Pm$, we have
    (i)~$\sup_{\vartheta \in \Theta} \Vert \E_{p_t}[\Hessian (X_t, \vartheta)] \Vert \leq c < \infty$, 
    (ii)~$\E_{p_t} \big[ \Hessian ( X_t, \vartheta ) \big] \prec \zeromatrix$ for all $\vartheta\in \Theta$,
    (iii)~$\sup_{\vartheta \in \Theta} \big\Vert (\partial/\partial \vartheta_j) \, \E_{p_t} \big[ \Hessian(X_t, \vartheta) \big] \big\Vert < \infty$, for all $j=1,\dots,k$, and 
    (iv)~$A \scaling = d_{t-1} \Omega_{t-1}$ for some $\mathcal{F}_{t-1}$-measurable $d_{t-1} > 0$.
\end{assumptionp}

Part (i) of Assumption~\ref{ass:HBNT} requires the Hessian to be \emph{bounded} (denoted $\mathcal{B}$). Part~(ii) imposes the expected Hessian to be \emph{negative definite} (denoted $\mathcal{N}$), while part~(iii) demands the uniform boundedness of the \emph{third} derivatives (denoted $\mathcal{T}$) as in \citet[Prop.~3]{Creal2024GMM}. Part (iv) is critical in the multivariate case to ensure that \eqref{eqn:EGMMcondition2} holds, while in the univariate case it is automatically satisfied.

Unfortunately, even Assumption~\ref{ass:HBNT} can be hard to verify. In particular, the expected Hessian in part~(ii) is typically not negative definite unless $\Hessian(y,\vartheta)$ itself is negative definite for all $y\in\mathcal{Y}$ and $\vartheta\in\Theta$. If $\Hessian(y,\vartheta)$ fails to be negative definite on parts of $\mathcal{Y}$, one would need to restrict $p_t$ from placing too much mass there; since $p_t$ is uncontrolled, in practice this leaves few alternatives to imposing $\Hessian(y,\vartheta)\prec \zeromatrix$ for all $y\in\mathcal{Y},\vartheta\in\Theta$.

Thus, EGMM improvements can be straightforwardly guaranteed only for densities that are log concave. For general densities and generic SD updates (with $A\scaling \succ O_k$, however small), EGMM improvements are not ensured. This underscores that the equivalence condition~\eqref{eqn:EGMM equiv condition} is not directly useful for model design, unlike the SEE condition~\eqref{eqn:ExpectedScoreEquivalence}.

\subsection{TKL of \citet{blasques2015information}}
\label{sec:TKLComparison}

\citet{blasques2015information} were the first to ask what performance guarantees SD updates can deliver. Rather than taking expectations as in Theorem~\ref{thm:EKLequivalenceIndepCopy}, they study a localized KL divergence. In the univariate case ($l=k=1$), they introduce a KL-type measure that restricts $x$ to the vicinity of $y_t$ by \emph{trimming} outcomes:
	\begin{align}
		\label{eq:LKLDefinition}
		\rkl_{B}(p_{t} \Vert f_\tpu) 
		:= \int_{B} \log \left( \frac{p_t(x)}{f(x | \vartheta_\tpu)} \right)p_t(x) \dd x.
	\end{align}
Here, the integration is restricted to the neighborhood
$\ymdel \equiv \ymdellong := \{x\in\Ym:\lvert x-y_t\rvert \le \delta\}$
around the realization $y_t$ for some (small) $\delta>0$, so outcomes outside this ball are trimmed. \citet{blasques2015information} call an updating scheme $\phi$ (locally) ``optimal'' if it necessarily reduces the TKL divergence for any $y_t\in\Ym$, i.e., if
	\begin{align}
		\label{eq:LKLCriterion}
		\Delta^{\rkl}_{\delta}(\phi)
		\equiv \Delta^{\rkl}_{\delta} (\phi|y_t, \vartheta_\tbu,  p_{t}) 
		:= \rkl_B \big(p_{t} \Vert f_\tpu\big) - \rkl_B \big(p_{t} \Vert f_\tbu \big)
		< 0.
	\end{align}
    The criterion $\Delta^{\rkl}_{\delta}(\phi)<0$, which we call \emph{TKL reducing}, requires that updating the model density from $f_{\tbu}$ to $f_{\tpu}$ decreases its discrepancy from the true density $p_t$, at least locally around the observation $y_t$. Roughly, \citet[Proposition~1--2]{blasques2015information} show that SD updates are unique in guaranteeing $\Delta^{\rkl}_{\delta}(\phi)<0$ for all $y_t\in\Ym$. This seems to offer an appealing theoretical property that distinguishes SD updates from other updating rules.

However, while \citet[Propositions~1-2]{blasques2015information} are formally correct, criterion~\eqref{eq:LKLDefinition} relative to which they establish improvement turns out to be vacuous. To illustrate, suppose $f_{\tpu}(x)>f_{\tbu}(x)$ for all $x\in B$. Then a straightforward calculation shows that, for any $p_t$,
    \begin{align*}
        &\hspace{-0.5cm}\rkl_{B}(p_{t} \Vert f_\tpu) - \rkl_{B}(p_{t} \Vert f_\tbu) \\
		&= \int_{\ymdellong} \log \left( \frac{p_t(x)}{f(x | \vartheta_\tpu)} \right)p_t(x) \dd x - \int_{\ymdellong} \log \left( \frac{p_t(x)}{f(x | \vartheta_\tbu)} \right)p_t(x) \dd x \\
        &= \int_{\ymdellong}  \left[ \log \big( f(x | \vartheta_\tbu) \big) - \log \big(f(x | \vartheta_\tpu)\big)  \right] p_t(x) \dd x \quad <\; 0,
    \end{align*}
    where the negativity follows as $f_\tpu(x) > f_\tbu(x)$ for all $x \in B$.\footnote{
        SD models, which update according to the derivative of the model density, are designed to improve the model log-density at $y_t$, such that for smooth densities, they also improve the model density for all $x \in \ymdellong$ for $\delta$ small enough.
        However, since the improper TKL measure is disconnected from $p_t$ such updates cannot be interpreted as genuine improvements toward the true density $p_t$.}
  Because the negativity holds independently of the truth $p_t$, we obtain the puzzling implication
\begin{align}
    \label{eqn:TKL_Problem_Intro}
    \rkl_B(p_t\Vert f_\tpu) < \rkl_B(p_t\Vert f_\tbu) \quad \text{ whenever } \quad 
    f_\tpu(x) > f_\tbu(x) \;\; \text{for all } x \in B.  
\end{align}
This shows that trimming-based localization yields a criterion that is disconnected from the true density and hence uninformative. For instance, even if $f_\tbu=p_t$, the TKL criterion would still favor moving away from the true density $p_t$ whenever $f_\tpu>f_\tbu$ on $B$. This is also evident in Example~\ref{exmpl:LKLProblem1} and panel (b) of Figure~\ref{fig:IllustrationAlmostSure}, where the updated (blue) density is not closer to $p_t$ near $y_t$ than the original (red) density, despite the TKL improvement.

The puzzling implication in \eqref{eqn:TKL_Problem_Intro} (i.e., that improvements are always possible irrespective of $p_t$) has long been known in  the literature on localized scoring rules. Since the KL divergence corresponds to the expected logarithmic score, TKL relates to the trimmed logarithmic score of \citet{amisano2007comparing}. Subsequent work has shown that trimming yields an improper scoring rule, and that localization should instead use censoring to preserve (strict) propriety; see \citet{diks2011likelihood}, \citet{gneiting_comparing_2011}, and \citet{depunder2023localzing}. Following this suggestion of censoring, in Appendix~\ref{sec:CKL} we show that sufficiently small score-driven updates are censored KL (CKL) reducing if and only if $p_t(y_t)>f(y_t| \vartheta_{t| t-1})$. While the dependence on the true density $p_t$ is thus restored, this condition is not verifiable in practice and provides little guidance for model construction. Having shown that the TKL measure is an improper criterion, we do not consider it further.

\section{Examples}
\label{sec:Example}

Here we take several popular model densities and investigate whether the respective guarantees for SD updates being EKL, CEV, MSE or EGMM reducing can be applied.
Section~\ref{sec:ExampleUnivariate} presents results for eleven models with a univariate time-varying parameter and Section~\ref{sec:ExampleBivaraiteGaussian} considers a Gaussian model with a bivariate location-scale time-varying parameter.

To interpret the results below, recall that the performance measures' applicability to SD models primarily differs in their requirements on the expected Hessian. EKL improvements require only (localized) boundedness (Assumptions~\ref{ass:HB}, \ref{ass:HLB}), whereas CEV, MSE, and EGMM improvements (additionally) demand negative definiteness, which is typically possible only if the realized Hessian itself is almost surely negative definite. 
In the multivariate case, Assumptions~\ref{ass:HN} and Assumption~\ref{ass:HBNT} further impose (differing) restrictions on the learning-rate and scaling matrices.

\subsection{Applicability to eleven univariate model densities}
\label{sec:ExampleUnivariate}

Table~\ref{tab2} illustrates the broad applicability of our EKL reduction guarantees of SD updates across eleven model densities with a univariate time-varying parameter. Nine of these models are adapted from \citet{koopman2016predicting}, complemented by two conditional location (local-level) models with Gaussian and Student’s~$t$ distributions. The collection is deliberately diverse, covering models with time-varying parameters for level, scale, dependence, count, intensity, and duration. These densities---or probability functions, in the case of discrete observations (see Remark~\ref{rem:DiscreteEKL})---are postulated by the researcher, while the true data-generating density remains unknown. 
In Appendix~\ref{sec:DetailsExamples}, we provide full details of the postulated models in Table~\ref{tab:DGPs} and discuss three models in detail, related to (i) intensities, (ii) negative binomial counts, and (iii) the Student’s $t$ location model.

For each of the eleven models, Table~\ref{tab2} reports the range of the expected Hessians, together with the moment conditions on the class $\Pm$ of true distribution $p_t$ required for these results.
The ranges are then used to assess whether the corresponding assumptions are satisfied---indicated by check or cross marks---for each performance criterion. 
The required moment conditions are mild, involving only first and second moments, and ensure that the random variables are not zero almost surely under the true distribution.

\begin{table}[t]
\centering
\caption{Applicability of performance guarantees to eleven postulated densities. \label{tab2}}
\begin{footnotesize}
\begin{threeparttable}
\begin{tabular}{l@{\hspace{0.4cm}}l@{\hspace{0.2cm}}c@{\hspace{0.2cm}}c@{\hspace{0.2cm}}c@{\hspace{0.2cm}}c@{\hspace{0.2cm}}c@{\hspace{0.2cm}}c@{\hspace{0.2cm}}c@{\hspace{0.2cm}}c}
\toprule
\multicolumn{2}{l}{\multirow{1}{*}{\bf Postulated model}}
&   \multicolumn{2}{c}{\multirow{1}{*}{\bf Properties of expected Hessian}}
&  \multicolumn{4}{c}{\multirow{1}{*}{\bf Performance criteria}}
\\
\cmidrule(r{5pt}l{5pt}){1-2} 
\cmidrule(r{5pt}l{5pt}){3-4} 
\cmidrule(r{5pt}l{5pt}){5-8}

\multirow{2}{*}{ } &	\multirow{2}{*}{ } &
\multirow{2}{*}{ } & \multirow{2}{*}{ }
& EKL
& EKL
& CEV/MSE
& EGMM
\\

&	 & Moment condition
& Range of
&  Thm.~\ref{thm:EKLequivalenceIndepCopy}
&  Thm.~\ref{thm:EKLequivalenceBounded}
& Pos.~\eqref{eqn:CEVcondition2} 
& Pos.~\eqref{eqn:EGMMcondition1}  
\\

Type &	Distribution & (restriction on $\mathcal{P}$)
& expected Hessian
& \ref{ass:HB}  
& \ref{ass:HLB}
& \ref{ass:HN}  
& \ref{ass:HBNT} 
\\

\cmidrule(r{5pt}l{5pt}){1-2} 
\cmidrule(r{5pt}l{5pt}){3-4} 
\cmidrule(r{5pt}l{5pt}){5-8} 

Count & 	Poisson
&  \multicolumn{1}{c}{$-$}
& \multicolumn{1}{c}{$(-\infty, 0)$}
& \xmark
& \cmark
& \cmark
& \xmark

\\ \addlinespace
Count &	Neg. Binomial
& 	\multicolumn{1}{l}{$\E_{p_t}[Y_t]<\infty$}
& \multicolumn{1}{c}{$\left[- \tfrac{\xi + \mathbb{E}_{p_t}[Y_t]}{4}, 0 \right)$}
& \cmark
& \cmark
& \cmark
& \cmark
\\  \addlinespace

Intensity  &	Exponential
& 	\multicolumn{1}{l}{$0 < \E_{p_t}[Y_t]<\infty$}
& \multicolumn{1}{c}{$(-\infty, 0)$}
& \xmark
& \cmark
& \cmark
& \xmark
\\  \addlinespace

Duration  &	Gamma
& 	\multicolumn{1}{l}{$0 < \E_{p_t}[Y_t]<\infty$}
& \multicolumn{1}{c}{$(-\infty, 0)$}
& \xmark
& \cmark
& \cmark
& \xmark
\\  \addlinespace

Duration  &	Weibull
& 	\multicolumn{1}{l}{$0 < \E_{p_t} [Y_t^\xi] <\infty$}
& \multicolumn{1}{c}{$(-\infty, 0)$}
& \xmark
& \cmark
& \cmark
& \xmark

\\  \addlinespace

Volatility  &	Gaussian
& 	\multicolumn{1}{l}{$0 < \E_{p_t}[Y_t^2]<\infty$}
& \multicolumn{1}{c}{$(-\infty, 0)$}
& \xmark
& \cmark
& \cmark
& \xmark

\\  \addlinespace

Volatility  &	Student's \emph{t}
& \multicolumn{1}{c}{$0<\E_{p_t}[Y_t^2]\qquad\;\;\;$}
& \multicolumn{1}{c}{$[-\tfrac{\nu+1}{8},0)$}
& \cmark 
& \cmark
& \cmark
& \cmark
\\ \addlinespace

Dependence  &	Gaussian
& \multicolumn{1}{l}{$\E_{p_t} [Y^2_{jt}]<\infty$}
& \multicolumn{1}{c}{$\left(-\infty,\tfrac{1}{4}\right]$}
& \xmark
& \cmark
& \xmark
& \xmark
\\ \addlinespace

Dependence  &	Student's \emph{t}
& \multicolumn{1}{c}{$-$}
& \multicolumn{1}{c}{$[-\tfrac{\nu+1}{4},\tfrac{1}{4}]$}
& \cmark
& \cmark
& \xmark
& \xmark
\\ \addlinespace

Level  &	Gaussian 
& \multicolumn{1}{c}{$-$} 
& \multicolumn{1}{c}{$\left[-\sigma^{-2},-\sigma^{-2}\right]$}
& \cmark 
& \cmark 
& \cmark 
& \cmark 
\\ \addlinespace

Level  &	Student's \emph{t} & 
\multicolumn{1}{c}{$-$} 
& \multicolumn{1}{c}{$\left[-\tfrac{\nu+1}{\nu \, \sigma^2},\tfrac{\nu+1}{8\,\nu \, \sigma^2}\right]$}
& \cmark 
& \cmark 
& \xmark
& \xmark\\
\bottomrule
 \end{tabular} 
 \vspace{3pt}
{\textbf{NOTE:} The column \emph{Postulated model} specifies the model type and distribution. The stated range for the expected Hessian, $\E_{p_t}[\Hessian(Y_t,\vartheta)]$, holds under the moment conditions in the column \emph{Moment condition} (with $j\in\{1,2\}$), which restricts the class $\Pm$ of permitted true distributions; these conditions are verified in Table~\ref{tab:DGPs} in Appendix~\ref{sec:DetailsExamples}. 
The lower bounds in the moment conditions exclude $Y_t=0$ a.s.\ and are only required to satisfy Assumption~\ref{ass:HN}.
In \emph{Performance criteria}, line~1 gives the criterion, line~2 the corresponding result/condition, and line~3 the relevant assumption (in addition to Assumption~\ref{ass:thmiffeklindepc}). Check and cross marks in the table indicate which assumptions hold.
}
\end{threeparttable}
\end{footnotesize}
\end{table}

From the ranges and the expected Hessians reported in Table~\ref{tab:DGPs}, we observe that Assumption~\ref{ass:HB} is satisfied for five of the eleven models, thereby guaranteeing EKL improvements for classes $\Pm$ that are only restricted by the given moment conditions.
For the remaining six models, the expected Hessian is unbounded in the time-varying parameter. This issue is addressed by Assumption~\ref{ass:HLB}, which applies to all eleven models. Under this assumption, Theorem~\ref{thm:EKLequivalenceBounded} and Proposition~\ref{prop:CSSD_SEE} ensure that clipped SD updates, with sufficiently large clipping constants, are EKL reducing.

In contrast, the CEV and MSE criteria of \citet{Gorgi2023} rely on a negative, though potentially unbounded, expected Hessian as specified in Assumption~\ref{ass:HN}, which is satisfied by eight of the models. For the three practically relevant cases where this assumption fails, it appears impossible to impose reasonable conditions on $\Pm$ that would restore Assumption~\ref{ass:HN}; see Appendix~\ref{sec:DetailsExamples} for details on the local-level model with a Student’s $t$ distribution discussed in \citet{harvey2014filtering}. Finally, guarantees for the EGMM criterion hold only in three cases where the strongest Assumption~\ref{ass:HBNT} is met.

In sum, Table~\ref{tab2} shows that the EKL criterion, under the mild Assumption~\ref{ass:HLB}, provides the only formal result currently applicable across the full class of models considered.

\subsection{Bivariate dynamic parameter: Gaussian location-scale model}
\label{sec:ExampleBivaraiteGaussian}

We now illustrate in a simple bivariate setting that only the EKL result remains applicable. Consider the Gaussian location-scale model $\mathcal{N}\big(\mu_{\tbu},\,\exp(\lambda_\tbu)\big)$ with time-varying mean $\mu_{\tbu}\in \mathbb{R}$ and variance $\exp(\lambda_{\tbu})>0$, which we collect as a bivariate time-varying parameter $\vartheta_{\tbu} = (\mu_{\tbu}, \lambda_{\tbu})^\top \in \R^2 = \Theta$.
The exponential link guarantees positivity for any $\lambda_\tbu\in \mathbb{R}$. 
For $\vartheta=(\mu, \lambda)^\top \in \mathbb{R}^2$,
a direct calculation for the Hessian matrix yields
\[
\Hessian(x,\vartheta)
=
\begin{pmatrix}
-\,\mathrm{e}^{-\lambda} & (\mu-x)\,\mathrm{e}^{-\lambda}\\[2pt]
(\mu-x)\,\mathrm{e}^{-\lambda} & -\tfrac{1}{2}(\mu-x)^2 \mathrm{e}^{-\lambda}
\end{pmatrix}.
\]
As the top-left element is negative, while the determinant $-\tfrac{1}{2}(\mu-x)^2 \mathrm{e}^{-2\lambda}$ is non-positive, the Hessian is indefinite (unless $x=\mu$, where the Hessian is negative semi-definite).

Next, we consider the expected Hessian. 
Suppose $X_t \mid \Fm_{t-1} \sim p_t \in \Pm$, the class restricted to finite means and variances, denoted $\mu_{p_t}$ and $\sigma_{p_t}^2$. Then, for any $\mu \in \mathbb{R}$, we have $\E_{p_t}[(\mu-X_t)^2] = (\mu_{p_t}-\mu)^2+\sigma_{p_t}^2$
and the expected Hessian becomes
\begin{equation}
\E_{p_t}[\Hessian(X_t,\vartheta)]
=
\begin{pmatrix}
-\,\mathrm{e}^{-\lambda} & (\mu-\mu_{p_t})\,\mathrm{e}^{-\lambda}\\[2pt]
(\mu-\mu_{p_t})\,\mathrm{e}^{-\lambda} & -\tfrac{1}{2}\big((\mu_{p_t}-\mu)^2+\sigma_{p_t}^2\big)\,\mathrm{e}^{-\lambda}
\end{pmatrix}.
\label{eq:expected hessian}
\end{equation}
As this matrix is bounded for $(\mu,\lambda)^\top$ in compact subsets of $\mathbb{R}^2$, Assumption~\ref{ass:HLB} holds. By Theorem~\ref{thm:EKLequivalenceBounded} and Proposition~\ref{prop:CSSD_SEE}, therefore, clipped SD updates with sufficiently large clipping constants and small learning rates are EKL reducing w.r.t.\ $\Pm$.

For the CEV, MSE, and EGMM criteria with $\Omega_{t-1}=\unitmatrix$, Assumptions \ref{ass:HN}~(ii) and~\ref{ass:HBNT}~(iv) force $A\scaling$ to be a scalar multiple of the identity and require a negative definite expected Hessian. In \eqref{eq:expected hessian}, the expected Hessian has negative diagonal entries, but its determinant is
$\frac{1}{2}\big(\sigma_{p_t}^2-(\mu-\mu_{p_t})^2\big)\mathrm{e}^{-2\lambda}$.
When $\mu$ is sufficiently far from $\mu_{p_t}$, this determinant is negative, so the expected Hessian is \emph{indefinite} (i.e., it has both positive and negative eigenvalues). Negative definiteness of the expected Hessian in Assumptions~\ref{ass:HN} and~\ref{ass:HBNT} fails; i.e., no improvement guarantees for CEV, MSE, or EGMM are available.

In sum, even in this simple bivariate case, a Gaussian distribution with time-varying mean and variance (e.g., as in ARMA--GARCH-type specifications), Hessian negative definiteness fails almost surely as well as in expectation.
Improvement guarantees arise only under the EKL criterion, as proposed here, under a localized expected-Hessian boundedness condition that is considerably less restrictive than alternative versions in the literature.

\section{Conclusion}
\label{sec:Conclusion}

We have characterized score-driven (SD) updates via the expected Kullback-Leibler (EKL) measure. We showed that EKL improvements occur if and only if the expected parameter adjustment aligns with the expected score; i.e., their inner product should be positive.  This equivalence continues to hold under merely locally bounded Hessians when updates are clipped. The resulting conditions are interpretable and constructive, directly guiding update design. By deriving explicit upper bounds on admissible learning-rate matrices in terms of score moments, we also connected SD models to adaptive optimization methods.

Relative to recent alternative performance measures, our approach is advantageous in
(i) yielding equivalence conditions directly useful for model construction,
(ii) requiring the mildest Hessian conditions, as illustrated by examples,
(iii) extending naturally to the multivariate case without restrictive choices of weighting, scaling, or learning-rate matrices, and
(iv) providing intuitive, constructive upper bounds for learning-rate matrices.

Together, these findings provide a rigorous justification for SD updates and establish the EKL divergence as their natural information-theoretic foundation.

\vspace{-0.5cm}
\section*{Acknowledgments}
We thank Peter Boswijk, Cees Diks, Dick van Dijk, Simon Donker van Heel, Andrew Harvey, Yi He,
Frank Kleibergen, Roger Laeven, Alessandra Luati, André Lucas, Bram van Os, Andrew Patton, Phyllis Wan 
and Chen Zhou for their valuable comments. We also thank participants at seminars at Heidelberg University and 
Goethe University Frankfurt, and at the 2024 ISF (Dijon), 
the 2024 Bernoulli-IMS World Congress (Bochum), the 2025 QFFE Conference (Marseille), the 2025 IAAE Annual Conference (Turin), and  the 2025 CFE-CMStatistics Conference (London).
T.~Dimitriadis gratefully acknowledges support from the German Research Foundation (DFG) under 
project number~502572912.


\singlespacing
\spacingset{1} 
\setlength{\bibsep}{2pt} 
\putbib[Bibliography-LSPS-v2]
\end{bibunit}

\pagebreak 
\appendix

\begin{center}
\Large \textbf{ Online supplement to:  \\
``Expected Kullback-Leibler-based \\ 
characterizations of score-driven updates''}\\[3pt]
\normalsize
Ramon de Punder, Timo Dimitriadis, and Rutger-Jan Lange
\\
 \today
\end{center}

\noindent This supplementary material complements the main paper with proofs, additional derivations and further technical results. 
Proofs for all results in the main paper are given in Appendix~\ref{sec:ProofsMain}. 
Then, Appendix~\ref{sec:DetailsExamples} provides additional details for the model examples in Section~\ref{sec:Example}.
Appendix~\ref{app:statespace} gives details on the applicability of our results to state-space models and Appendix~\ref{sec:KalmanFilter} illustrates that the Kalman filter can be expressed as a SD model and is hence EKL reducing.
Then, Appendix~\ref{app:ISD} shows EKL reductions for implicit score-driven (ISD) updates and details for Quasi-SD (QSD) models are given in Appendix~\ref{app:qsd}. 
An alternative to guarantee CEV reductions based on a scaled Hessian is discussed in Appendix~\ref{sec:ScaledHessianCEVGeneral}. 
The Appendices~\ref{sec:CKL} and \ref{sec:CKLGeneral} establish equivalence conditions when SD and general updating rules, respectively, are Censored KL (CKL) reducing.
Proofs for the latter two appendix sections are given in Appendix~\ref{sec:ProofsSupplement}.

\setcounter{section}{0} 
\setcounter{subsection}{0}          
\renewcommand{\thesection}{\Alph{section}}
\setcounter{equation}{0}
\renewcommand\theequation{A.\arabic{equation}}

\newtheorem{atheorem}{Theorem}[section]
\newtheorem{alemma}{Lemma}[section]
\newtheorem{acorollary}{Corollary}[section]
\newtheorem{aassumption}{Assumption}[section]
\newtheorem{aproposition}{Proposition}[section]
\newtheorem{aexample}{Example}[section]
\newtheorem{aremark}{Remark}[section]
\newtheorem{adefinition}{Definition}[section]

\renewcommand{\theatheorem}{\thesection.\arabic{atheorem}}
\renewcommand{\thealemma}{\thesection.\arabic{alemma}}
\renewcommand{\theacorollary}{\thesection.\arabic{acorollary}}
\renewcommand{\theaassumption}{\thesection.\arabic{aassumption}}
\renewcommand{\theaexample}{\thesection.\arabic{aexample}}
\renewcommand{\thearemark}{\thesection.\arabic{aremark}}
\renewcommand{\theadefinition}{\thesection.\arabic{adefinition}}

\counterwithin{equation}{section}
\renewcommand{\thetable}{\thesection.\arabic{table}}
\counterwithin{table}{section}
\renewcommand{\thetable}{\thesection.\arabic{table}}
\counterwithin{figure}{section}
\renewcommand{\thefigure}{\thesection.\arabic{figure}}

\normalem




\begin{bibunit}[chicago]

\section{Proofs of the main results}
\label{sec:ProofsMain}

\begin{proof}[Proof of Theorem~\ref{thm:EKLequivalenceIndepCopy}]
	Let $\phi$ be an arbitrary updating rule that updates the time-varying parameter from $\vartheta_\tbu$ to $\vartheta_\tpu = \phi(y,\vartheta_\tbu)$ by using the observation $y$, and which satisfies the conditions of Theorem \ref{thm:EKLequivalenceIndepCopy}.
	For a given $\phi$ and some generic $\kappa > 0$, let the downscaled updating rule $\phi_\kappa$ imply the updated parameter $\vartheta^\kappa_\tpu(y)$, which is such that 
	\begin{align*}
		\Delta \phi_\kappa(y, \vartheta_\tbu) 
		= \phi_\kappa(y,\vartheta_\tbu) - \vartheta_\tbu
		\stackrel{!}{=} \vartheta^\kappa_\tpu(y) - \vartheta_\tbu
		= \kappa \Delta \phi(y,\vartheta_\tbu).
	\end{align*}
    Notice that for any $y \in \mathcal{Y}$ and $\vartheta_\tbu \in \Theta$, the condition $\vartheta^\kappa_\tpu(y) \in \Theta$ is  guaranteed for $\kappa > 0$ small enough by the convexity of $\Theta$.
    
	In the following, we derive the explicit orders how $\Delta^{\mathsf{EKL}}(\phi_\kappa)$ shrinks to zero as $\kappa$ tends to zero, which we use to establish that, for $\kappa$ small enough, $\Delta^{\mathsf{EKL}}(\phi_\kappa)$ approaches zero from below.	For any $x,y \in \mathcal{Y}$, applying the path-integral version of the mean-value theorem to $\log f(x\vert \vartheta^\kappa_\tpu(y))$ at $\vartheta_\tbu$ gives
    \allowdisplaybreaks
	\begin{align*} 
		&\log f(x\vert \vartheta^\kappa_\tpu(y)) -  \log f(x\vert \vartheta_\tbu) 
        = \score(x, \vartheta_\tbu)^\top (\vartheta^\kappa_\tpu(y)-\vartheta_\tbu) \\
        &\hspace{4.5cm}+ (\vartheta^\kappa_\tpu(y) - \vartheta_\tbu)^\top  \left( \int_0^1 (1-u) \Hessian \big(x, \widetilde{\vartheta}_u^\kappa(y) \big) \mathrm{d}u \right)(\vartheta^\kappa_\tpu(y) - \vartheta_\tbu),
	\end{align*}
    where $\widetilde{\vartheta}^\kappa_u(y) = u \vartheta_\tbu + (1-u) \vartheta^\kappa_\tpu(y)$, with $u\in [0,1]$, moves on the line between $\vartheta_\tbu$ and $\vartheta^\kappa_\tpu(y)$.
	Plugging this expansion into the expected KL difference in Definition~\ref{def:EKL} and \eqref{eqn:EKLDef} yields
	\begingroup
	\allowdisplaybreaks
	\begin{align*}
		&-\Delta^{\mathsf{EKL}}(\phi_\kappa| \vartheta_\tbu, p_t) \\
		&\qquad= \E_{X_t,Y_t \sim p_t} \left[\log f(X_t\vert \vartheta^\kappa_\tpu(Y_t))-\log f(X_t\vert \vartheta_\tbu)\right] \\
		&\qquad= \E_{X_t,Y_t \sim p_t} \left[\score(X_t, \vartheta_\tbu)^\top (\vartheta^\kappa_\tpu(Y_t)-\vartheta_\tbu) \right] \\ 
		&\qquad\qquad+  \E_{X_t,Y_t \sim p_t} \left[ (\vartheta^\kappa_\tpu(Y_t)-\vartheta_\tbu)^\top \left( \int_0^1 (1-u)\Hessian \big( X_t, \widetilde{\vartheta}_u^\kappa(Y_t) \big) \mathrm{d}u \right) (\vartheta^\kappa_\tpu(Y_t)-\vartheta_\tbu) \right]\\ 
		&\qquad= \E_{p_t}\left[\Delta \phi_{\kappa}(Y_t,\vartheta_\tbu) \right]^\top \E_{p_t}\left[\score(X_t, \vartheta_\tbu) \right] \\
		&\qquad\qquad+  \E_{X_t,Y_t \sim p_t}\left[ \Delta \phi_{\kappa}(Y_t, \vartheta_\tbu)^\top \left( \int_0^1 (1-u)\Hessian \big( X_t, \widetilde{\vartheta}_u^\kappa(Y_t) \big) \mathrm{d}u \right)   \Delta \phi_{\kappa}(Y_t, \vartheta_\tbu)\right] \\
		&\qquad=\kappa \, \E_{p_t} \big[ \Delta \phi(Y_t, \vartheta_\tbu) \big]^\top  \E_{p_t}  \big[ \score(X_t, \vartheta_\tbu) \big]  \\
        &\qquad\qquad+ \kappa^2 \E_{X_t,Y_t \sim p_t}\left[ \Delta \phi(Y_t, \vartheta_\tbu)^\top \left( \int_0^1 (1-u)\Hessian \big( X_t, \widetilde{\vartheta}_u^\kappa(Y_t) \big) \mathrm{d}u \right)  \Delta \phi(Y_t, \vartheta_\tbu)\right] \\
		&\qquad=\kappa \, \E_{p_t} \big[ \Delta \phi(Y_t, \vartheta_\tbu) \big]^\top  \E_{p_t} \big[ \score(X_t, \vartheta_\tbu) \big] + \mathcal{O}(\kappa^2),
	\end{align*}
	\endgroup
    as $\kappa \downarrow 0$, where we have used that $X_t$ and $Y_t$ are independent (by assumption) and that $\sup_{\vartheta \in {\Theta}} \left\Vert \E_{p_t} \left[ \Hessian \big( X_t, \vartheta \big) \right] \right\Vert \leq c < \infty$ by Assumption~\ref{ass:HB} such that the $\mathcal{O}(\kappa^2)$ term in the penultimate line is bounded:
    \begin{align}
        &\hspace{-0.5cm}\left \vert \E_{X_t,Y_t \sim p_t}\left[ \Delta \phi(Y_t, \vartheta_\tbu)^\top \left( \int_0^1 (1-u)\Hessian \big( X_t, \widetilde{\vartheta}_u^\kappa(Y_t) \big) \mathrm{d}u \right)  \Delta \phi(Y_t, \vartheta_\tbu)\right]\right \vert \nonumber \\ 
        &\leq \E_{Y_t \sim p_t} \left[\left \vert \Delta \phi(Y_t, \vartheta_\tbu)^\top  \left( \int_0^1 (1-u)\E_{X_t \sim p_t} \left[  \Hessian \big( X_t, \widetilde{\vartheta}_u^\kappa(Y_t) \big) \right]  \mathrm{d}u \right)  \Delta \phi(Y_t, \vartheta_\tbu) \right \vert \right] \nonumber\\ 
        &\le \E_{Y_t \sim p_t} \left[ \left\Vert  \Delta \phi(Y_t, \vartheta_\tbu) \right\Vert^2 \left\Vert \int_0^1 (1-u) \E_{X_t \sim p_t} \left[ \Hessian \big( X_t, \widetilde{\vartheta}_u^\kappa(Y_t) \big) \right] \mathrm{d}u  \right\Vert  \right] \nonumber \\
        &\le \E_{Y_t \sim p_t} \left[ \left\Vert  \Delta \phi(Y_t, \vartheta_\tbu) \right\Vert^2 \int_0^1 (1-u) \left\Vert \E_{X_t \sim p_t} \left[ \Hessian \big( X_t, \widetilde{\vartheta}_u^\kappa(Y_t) \big) \right] \right\Vert \mathrm{d}u \right] \nonumber \\
        &\le \E_{Y_t \sim p_t} \left[ \left\Vert  \Delta \phi(Y_t, \vartheta_\tbu) \right\Vert^2 \sup_{\vartheta \in {\Theta}} \left\Vert \E_{p_t} \left[ \Hessian \big( X_t, \vartheta \big) \right] \right\Vert \left(\int_0^1 (1-u) \mathrm{d}u\right)   \right] \nonumber \\
        &\le \frac{c}{2} \, \E_{p_t} \left[ \left\Vert  \Delta \phi(Y_t, \vartheta_\tbu) \right\Vert^2  \right] < \infty,  \label{eqn:NormBound}
    \end{align}    
    using a norm bound for quadratic forms and that $\widetilde{\vartheta}_u^\kappa(Y_t)$ is independent of $X_t$.
   Moreover, $\E_{p_t} \left[ \Vert \Delta \phi(Y_t, \vartheta_\tbu) \Vert^2 \right]  $ is finite by Assumption~\ref{ass:thmiffeklindepc}.

	As $\E_{p_t} \big[ \Delta \phi(Y_t, \vartheta_\tbu) \big]^\top  \E_{p_t} \big[ \score(X_t, \vartheta_\tbu) \big] < \infty$ by assumption, multiplying the following terms yields that
	\begin{align}
		&\Delta^{\mathsf{EKL}}(\phi_\kappa| \vartheta_\tbu, p_t) \E_{p_t} \big[ \Delta \phi(Y_t,\vartheta_\tbu) \big]^\top  \E_{{p_t}} \big[ \score(X_t,\vartheta_\tbu) \big] \label{eqn:EKLProof1} \\
		&\qquad=-\Big(\kappa \,  \E_{p_t}\big[ \Delta \phi(Y_t, \vartheta_\tbu) \big]^\top \E_{p_t} \big[ \score(X_t, \vartheta_\tbu) \big] + \mathcal{O}(\kappa^2) \Big) \nonumber \\
		&\qquad\qquad \times \E_{p_t} \big[ \Delta \phi(Y_t, \vartheta_\tbu) \big]^\top  \E_{p_t} \big[ \score(X_t, \vartheta_\tbu) \big]  \nonumber \\
		&\qquad= -\kappa \Big( \E_{p_t} \big[ \Delta \phi(Y_t, \vartheta_\tbu) \big]^\top  \E_{p_t} \big[ \score(X_t, \vartheta_\tbu) \big] \Big)^2  + \mathcal{O}(\kappa^2). \label{eqn:EKLProof2}
	\end{align}
	For $\kappa$ small enough, the first term in \eqref{eqn:EKLProof2} dominates the $\mathcal{O}(\kappa^2)$ term such that \eqref{eqn:EKLProof2} is negative for $\kappa$ sufficiently small enough.
	This implies that the terms in \eqref{eqn:EKLProof1} have opposing signs (for $\kappa$  small enough), which we formalize in the following:
	
	Starting with the $\Longrightarrow$ direction of the proof, if $\E_{p_t} \big[ \Delta \phi(Y_t,\vartheta_\tbu) \big]^\top  \E_{{p_t}} \big[ \score(X_t,\vartheta_\tbu) \big] > 0$, then we can always find a $\bar \kappa$ small enough such that for all $\kappa \in (0,\bar \kappa]$, \eqref{eqn:EKLProof2} is negative, and hence, $\Delta^{\mathsf{EKL}}(\phi_\kappa| \vartheta_\tbu, p_t)$ must be negative as well for $\kappa \in (0,\bar \kappa]$.
	
	For the $\Longleftarrow$ direction of the proof, we assume that there exists a $\bar \kappa$ such that for all $\kappa \in (0,\bar \kappa]$, $\Delta^{\mathsf{EKL}}(\phi_\kappa| \vartheta_\tbu, p_t)$ is negative. 
	As there always exists a $\tilde \kappa \in (0,\bar \kappa]$ such that \eqref{eqn:EKLProof2} becomes negative, we can conclude that $\E_{p_t} \big[ \Delta \phi(Y_t,\vartheta_\tbu) \big]^\top  \E_{{p_t}} \big[ \score(X_t,\vartheta_\tbu) \big] > 0$, which concludes the proof.
\end{proof}

\begin{proof}[Proof of Theorem~\ref{thm:EKLequivalenceBounded}]
    We follow the proof of Theorem~\ref{thm:EKLequivalenceIndepCopy} to arrive at the expression 
    \begin{align}
        \begin{aligned}
        \label{eqn:ExpansionLocallyBoundedHessian}    
		&-\Delta^{\mathsf{EKL}}(\phi_\kappa| \vartheta_\tbu, p_t) \\
		&\qquad=\kappa \, \E_{p_t} \big[ \Delta \phi(Y_t, \vartheta_\tbu) \big]^\top  \E_{p_t}  \big[ \score(X_t, \vartheta_\tbu) \big]  \\
        &\qquad \qquad+ \kappa^2\E_{X_t,Y_t \sim p_t}\left[ \Delta \phi(Y_t, \vartheta_\tbu)^\top \left( \int_0^1 (1-u)\Hessian \big( X_t, \widetilde{\vartheta}_u^\kappa(Y_t) \big) \mathrm{d}u \right)  \Delta \phi(Y_t, \vartheta_\tbu)\right],
        \end{aligned}
	\end{align}
    where 
    $\widetilde{\vartheta}^\kappa_u(y) = u \vartheta_\tbu + (1-u) \vartheta^\kappa_\tpu(y)$, with $u\in [0,1]$, moves on the line between $\vartheta_\tbu$ and $\vartheta^\kappa_\tpu(y)$.

    Given $\vartheta_\tbu$ and the updating rule $\phi$  with the corresponding $d(\vartheta_\tbu) \ge \Vert \Delta\phi (Y_t,\vartheta_\tbu) \Vert$, we can choose $\kappa \in \big(0, r/d(\vartheta_\tbu)\big)$ such that $\kappa \, d(\vartheta_\tbu) < r$ and hence
    \begin{align*}
        \vartheta^\kappa_\tpu(Y_t)
        =  \vartheta_\tbu +  \kappa \big( \vartheta_\tpu(Y_t) - \vartheta_\tbu \big) \in 
       \widetilde{\Theta}_r(\vartheta_\tbu).
    \end{align*}
    As $\widetilde{\vartheta}^\kappa_u(Y_t) = u \vartheta_\tbu + (1-u) \vartheta^\kappa_\tpu(Y_t)$ is on the line between $\vartheta_\tbu$ and $\vartheta^\kappa_\tpu(Y_t)$, while $\widetilde{\Theta}_r(\vartheta_\tbu)$ is convex, we immediately have that $\widetilde{\vartheta}^\kappa_u(Y_t) \in \widetilde{\Theta}_r(\vartheta_\tbu)$.
    Hence, for the second-order term of \eqref{eqn:ExpansionLocallyBoundedHessian}, we get
    \begin{align*}
        &\hspace{-0.5cm} \left| \E_{X_t,Y_t \sim p_t}\left[ \Delta \phi(Y_t, \vartheta_\tbu)^\top \left( \int_0^1 (1-u)\Hessian \big( X_t, \widetilde{\vartheta}_u^\kappa(Y_t) \big) \mathrm{d}u \right)  \Delta \phi(Y_t, \vartheta_\tbu)\right] \right| \\
        &\le \E_{Y_t \sim p_t} \left[ \left\Vert  \Delta \phi(Y_t, \vartheta_\tbu) \right\Vert^2 \sup_{\vartheta \in \widetilde{\Theta}_r} \left\Vert \E_{X_t \sim p_t} \left[ \Hessian \big( X_t, \vartheta \big) \right] \right\Vert \left(\int_0^1 (1-u) \mathrm{d}u\right)   \right] \\
        &\le \frac{c}{2} \, \E_{p_t} \left[ \left\Vert  \Delta \phi(Y_t, \vartheta_\tbu) \right\Vert^2  \right],
    \end{align*}    
    where we repeated the steps leading to \eqref{eqn:NormBound}, now  using that the norm of the expected Hessian is bounded above by $c$ on the compact set $\widetilde{\Theta}_r(\vartheta_\tbu)$ by Assumption~\ref{ass:HLB}.
    Hence, the last term in \eqref{eqn:ExpansionLocallyBoundedHessian} is $\mathcal{O}(\kappa^2)$ such that following the remaining steps in the proof of Theorem~\ref{thm:EKLequivalenceIndepCopy} gives the desired result.
\end{proof}

\begin{proof}[Proof of Proposition~\ref{prop:CSSD_SEE}]
    For simplicity, we write $\score_t = \score(Y_t,\vartheta_\tbu)$ throughout the proof and define the weighted norm
    $\Vert x\Vert_W^2:=x^\top W x$ for a positive definite matrix $W$ and a vector $x$ of compatible dimensions. As usual, $\|\cdot\|$ denotes the Euclidean norm.
    
    For any event $E$, let $\mathbb{I}\{E\}$ denote the indicator function, which equals one if the event~$E$ occurs and zero otherwise.
    We begin by rewriting 
    
    \begin{align*}
    &\hspace{-0.3cm}\E_{p_t} \big[ \Delta \phi^d_{\mathrm{CSD}} (Y_t,\vartheta_\tbu) \big] \\
    &= \E_{p_t} \left[ \min \left\{1, \frac{d}{\Vert \Delta\phi_{\mathrm{SD}}(Y_t,\vartheta_\tbu)\Vert} \right\} \Delta\phi_{\mathrm{SD}}(Y_t,\vartheta_\tbu) \right] \\
    &= \E_{p_t} \left[ \min \left\{1, \frac{d}{\Vert A \scaling \score_t \Vert} \right\} A \scaling \score_t \right] \\
    &= A \scaling \E_{p_t} \left[ \mathbb{I}\{ \Vert A \scaling \score_t \Vert \le d\}  \score_t
    + \mathbb{I}\{ \Vert A \scaling \score_t \Vert > d\} \frac{d}{\Vert A \scaling \score_t \Vert} \score_t \right] \\
    &= A \scaling \E_{p_t} \left[ (1-\mathbb{I}\{ \Vert A \scaling \score_t \Vert > d\})  \score_t
    + \mathbb{I}\{ \Vert A \scaling \score_t \Vert > d\} \frac{d}{\Vert A \scaling \score_t \Vert} \score_t \right] \\
    &= A \scaling \E_{p_t} [ \score_t]
    - {A \scaling}\,\E_{p_t} \left[ \mathbb{I}\{ \Vert A \scaling \score_t \Vert > d\}
    \left(1-\frac{d}{\Vert A \scaling \score_t \Vert} \right) \score_t  \right] \\
    &= A \scaling \mu_{p_t}
    - {A \scaling}\,\E_{p_t} \left[ \mathbb{I}\{ \Vert A \scaling \score_t \Vert > d\}
    \left(1-\frac{d}{\Vert A \scaling \score_t \Vert} \right) \score_t  \right].
    \end{align*}
    
    Then, taking the inner product of this expression with $\mu_{p_t}$ gives
    \begin{align}
        &\hspace{-0.3cm}\mu_{p_t}^\top \E_{p_t} \big[ \Delta \phi^d_{\mathrm{CSD}} (Y_t,\vartheta_\tbu) \big]\nonumber  \\ 
        &= \mu_{p_t}^\top A \scaling \mu_{p_t}-   
        \mu_{p_t}^\top {A \scaling}\,\E_{p_t} \left[ \mathbb{I}\{ \Vert A \scaling \score_t \Vert > d\}
        \left(1-\frac{d}{\Vert A \scaling \score_t \Vert} \right) \score_t  \right] \nonumber \\
        &= \|\mu_{p_t}\|_{A \scaling}^2
        - \E_{p_t} \left[ \mathbb{I}\{ \Vert A \scaling \score_t \Vert > d\}
        \left(1-\frac{d}{\Vert A \scaling \score_t \Vert} \right)\, \mu_{p_t}^\top {A \scaling}\,\score_t  \right]  \nonumber \\
        &\geq \|\mu_{p_t}\|_{A \scaling}^2
        - \E_{p_t} \left[ \mathbb{I}\{ \Vert A \scaling \score_t \Vert > d\}
        \left(1-\frac{d}{\Vert A \scaling \score_t \Vert} \right)\,
        \big|\mu_{p_t}^\top {A \scaling}\,\score_t\big|  \right]  \nonumber \\
        &\geq \|\mu_{p_t}\|_{A \scaling}^2
        - \E_{p_t} \left[ \mathbb{I}\{ \Vert A \scaling \score_t \Vert > d\}
        \left(1-\frac{d}{\Vert A \scaling \score_t \Vert} \right)\,
        {\|\mu_{p_t}\|_{A\scaling}\, \|\score_t\|_{A\scaling}}  \right]  \nonumber \\
        &\geq \|\mu_{p_t}\|_{A \scaling}^2
        - {\|\mu_{p_t}\|_{A\scaling}}\,\E_{p_t}\!\left[ \mathbb{I}\{ \Vert A \scaling \score_t \Vert > d\}\,
        {\|\score_t\|_{A\scaling}}  \right]  \nonumber \\
        &\geq \|\mu_{p_t}\|_{A \scaling}^2
        - {\|\mu_{p_t}\|_{A\scaling}}\,
        \sqrt{\E_{p_t}\big[\mathbb{I}\{ \Vert A \scaling \score_t \Vert > d\}\big]\;
              \E_{p_t}\big[{\|\score_t\|_{A\scaling}^2}\big]} \nonumber \\
        &= \|\mu_{p_t}\|_{A \scaling}^2
        - {\|\mu_{p_t}\|_{A\scaling}}\,
        \sqrt{\P_{Y_t\sim p_t}(\Vert A \scaling \score_t \Vert > d)\;
              \E_{p_t}\big[{\|\score_t\|_{A\scaling}^2}\big]}
              \nonumber
              \\
        &= \|\mu_{p_t}\|_{A \scaling}^2
        - \|\mu_{p_t}\|_{A\scaling}\,
        \sqrt{\P_{Y_t\sim p_t}(\Vert A \scaling \score_t \Vert > d)\;
              {\Big(\|\mu_{p_t}\|_{A\scaling}^2 + \trace\!\big((A\scaling)\Sigma_{p_t}\big)\Big)}},
       \label{eqn:LastTermCSSDProof}
    \end{align}
    where $\P_{Y_t\sim p_t}$ denotes a probability under the true density $p_t$. Above, the first inequality is algebraic, the second and fourth follow by Cauchy--Schwarz, and the third leaves out the factor
    $\left(1-\frac{d}{\Vert A \scaling \score_t \Vert}\right)$, which is less than one given that $\{\Vert A \scaling \score_t \Vert > d\}$, as imposed by the indicator function.
    
    In deriving the last equality, \eqref{eqn:LastTermCSSDProof}, we use several trace equalities: 
    \begin{align*}
    \E_{p_t}\!\big[\|\score_t\|_{A \scaling}^2\big]
    &= \E_{p_t}\!\big[\score_t^\top (A\scaling)\score_t\big] \\
    &= \E_{p_t}\!\big[\trace\!\big((A\scaling)\score_t\score_t^\top\big)\big] \\
    &= \trace\!\big((A\scaling)\E_{p_t}[\score_t\score_t^\top]\big) \\
    &= \trace\!\big((A\scaling)(\mu_{p_t}\mu_{p_t}^\top+\Sigma_{p_t})\big) \\
    &= \mu_{p_t}^\top (A\scaling)\mu_{p_t} + \trace\!\big((A\scaling)\Sigma_{p_t}\big) \\
    &= \|\mu_{p_t}\|_{A\scaling}^2 + \trace\!\big((A\scaling)\Sigma_{p_t}\big).
    \end{align*}
    Here, we used the
    linearity of $\trace(\cdot)$ and $\E_{p_t}[\score_t\score_t^\top]=\mu_{p_t}\mu_{p_t}^\top+\Sigma_{p_t}$. To write $ \mu_{p_t}^\top (A\scaling)\mu_{p_t}=\|\mu_{p_t}\|_{A\scaling}^2$, we also used the fact that $A \scaling$ is positive definite by assumption.
    
    Finally, the lower bound~\eqref{eqn:LastTermCSSDProof} is positive whenever
    \begin{align*}
    &\|\mu_{p_t}\|_{A \scaling}^2
    >
    \|\mu_{p_t}\|_{A\scaling}\,
    \sqrt{\P_{Y_t\sim p_t}(\Vert A \scaling \score_t \Vert > d)\,
    \Big(\|\mu_{p_t}\|_{A\scaling}^2 + \trace\!\big((A\scaling)\Sigma_{p_t}\big)\Big)} \\
    &\qquad \iff\quad
    \|\mu_{p_t}\|_{A \scaling}
    >
    \sqrt{\P_{Y_t\sim p_t}(\Vert A \scaling \score_t \Vert > d)\,
    \Big(\|\mu_{p_t}\|_{A\scaling}^2 + \trace\!\big((A\scaling)\Sigma_{p_t}\big)\Big)} \\
    &\qquad \iff\quad
    \P_{Y_t\sim p_t}(\Vert A \scaling \score_t \Vert > d)
    <
    \frac{\|\mu_{p_t}\|_{A \scaling}^2}{\|\mu_{p_t}\|_{A\scaling}^2 + \trace\!\big((A\scaling)\Sigma_{p_t}\big)}.
    \end{align*}
    In going from the first to the second inequality, we divided by $\|\mu_{p_t}\|_{A\scaling}>0$, which is permitted since $A\scaling$ is positive definite and $\mu_{p_t}\neq 0$ by assumption. The third inequality is exactly condition \eqref{eq:probability clipping is active}. Hence, if
    \eqref{eq:probability clipping is active} holds, then
    $\mu_{p_t}^\top \E_{p_t} \big[ \Delta \phi^d_{\mathrm{CSD}} (Y_t,\vartheta_\tbu) \big] > 0$, which concludes the proof.
\end{proof}

\begin{proof}[Proof of Theorem~\ref{thm:UpperBoundLearningRate}]
    We consider the SD update with $\Delta \phi_{\mathrm{SD}} (Y_t,\vartheta_\tbu) = \mathcal{A}_{t-1} \score(Y_t,\vartheta_\tbu)$ with a general positive-definite matrix $\mathcal{A}_{t-1}$ and follow the steps in the proof of Theorem~\ref{thm:EKLequivalenceIndepCopy}.
    Specifically, for $\vartheta_\tpu^{\mathcal{A}}(y) := \vartheta_\tbu + \Delta \phi_{\mathrm{SD}} (Y_t,\vartheta_\tbu) = \vartheta_\tbu + \mathcal{A}_{t-1} \score(Y_t,\vartheta_\tbu)$,
    we use the integral version of the mean-value theorem with $\widetilde{\vartheta}^\mathcal{A}_u(y) = u \vartheta_\tbu + (1-u) \vartheta_\tpu^{\mathcal{A}}(y)$ for $u \in [0,1]$, which gives
    \allowdisplaybreaks
	\begin{align*}
		&- \Delta^{\mathsf{EKL}}(\phi_{\mathrm{SD}} | \vartheta_\tbu, p_t) \\
		&\qquad= \E_{X_t,Y_t \sim p_t} \left[\log f(X_t\vert \vartheta^{\mathcal{A}}_\tpu(Y_t))-\log f(X_t\vert \vartheta_\tbu)\right] \\
		&\qquad= \E_{X_t,Y_t \sim p_t} \biggl[ \score(X_t, \vartheta_\tbu)^\top (\vartheta^{\mathcal{A}}_\tpu(Y_t)-\vartheta_\tbu) \\
        &\hspace{2cm} + \bigl( \vartheta^{\mathcal{A}}_\tpu(Y_t)-\vartheta_\tbu \bigr)^\top \left( \int_0^1 (1-u) \Hessian \big( X_t, \widetilde{\vartheta}_u^{\mathcal{A}}(Y_t) \big) \mathrm{d}u \right)  \bigl(\vartheta^{\mathcal{A}}_\tpu(Y_t)-\vartheta_\tbu \bigr) \biggr]\\ 
		&\qquad= \E_{p_t} \bigl[ \score(X_t, \vartheta_\tbu) \bigr]^\top \mathcal{A}_{t-1}  \E_{p_t} \bigl[ \score(Y_t, \vartheta_\tbu) \bigr] \\
		&\hspace{2cm} + \E_{p_t}  \left[  \score(Y_t, \vartheta_\tbu)^\top  \mathcal{A}_{t-1} \left( \int_0^1 (1-u) \Hessian \big( X_t, \widetilde{\vartheta}_u^{\mathcal{A}}(Y_t) \big) \mathrm{d}u \right) \mathcal{A}_{t-1}  \score(Y_t, \vartheta_\tbu) \right] \\ 
		&\qquad\ge \E_{p_t} \bigl[ \score(X_t, \vartheta_\tbu) \bigr]^\top \mathcal{A}_{t-1} \E_{p_t} \bigl[ \score(Y_t, \vartheta_\tbu) \bigr] 
        - \frac{c}{2} \trace \left( \mathcal{A}_{t-1}^2(\Sigma_{p_t} + \mu_{p_t} \mu_{p_t}^\top) \right).
	\end{align*}
	For the inequality in the last line, we have used that $\sup_{\vartheta \in \Theta} \Vert \E_{p_t} \Hessian(x,\vartheta)\Vert \le c<\infty$ from Assumption~\ref{ass:HB} as in \eqref{eqn:NormBound}, and cyclical rotation inside the trace operator to get
	\begin{align*}
        &\left| \E_{X_t, Y_t \sim p_t}  \left[  \score(Y_t, \vartheta_\tbu)^\top  \mathcal{A}_{t-1} \left( \int_0^1 (1-u) \Hessian \big( X_t, \widetilde{\vartheta}_u^{\mathcal{A}}(Y_t) \big) \mathrm{d}u \right) \mathcal{A}_{t-1}  \score(Y_t, \vartheta_\tbu) \right] \right| \\
        &\quad\leq \E_{Y_t \sim p_t}  \left[ \left|  \score(Y_t, \vartheta_\tbu)^\top  \mathcal{A}_{t-1} \left( \int_0^1 (1-u) \E_{X_t \sim p_t} \left[\Hessian \big( X_t, \widetilde{\vartheta}_u^{\mathcal{A}}(Y_t) \big) \right] \mathrm{d}u \right) \mathcal{A}_{t-1}  \score(Y_t, \vartheta_\tbu) \right| \right]  \\
     &\quad\leq   \E_{Y_t \sim p_t}  \left[\score(Y_t, \vartheta_\tbu)^\top  \mathcal{A}_{t-1}^2  \score(Y_t, \vartheta_\tbu)\right] \int_0^1 (1-u) \mathrm{d}u \; \sup_{\vartheta \in \Theta} \left \Vert \E_{X_t \sim p_t} \left[\Hessian \big( X_t, \vartheta \big) \right]  \right \Vert   \\
		&\quad\leq \frac{c}{2} \E_{p_t} \Bigl[ \score(Y_t, \vartheta_\tbu)^\top  \mathcal{A}_{t-1}^2 \score(Y_t, \vartheta_\tbu) \Bigr] \\
        &\quad= \frac{c}{2} \trace \left( \mathcal{A}_{t-1}^2  \E_{p_t} \left[\score(Y_t,\vartheta_\tbu) \score(Y_t,\vartheta_\tbu)^\top \right] \right) \\
        &\quad= \frac{c}{2} \trace \left( \mathcal{A}_{t-1}^2(\Sigma_{p_t} + \mu_{p_t} \mu_{p_t}^\top) \right).
	\end{align*}
    Hence $\Delta^{\mathsf{EKL}}(\phi_{\mathrm{SD}} | \vartheta_\tbu, p_t) < 0 $ if
    \begin{align}
        \label{eqn:GeneralLowerBound}
        \mu_{p_t}^\top \mathcal{A}_{t-1} \mu_{p_t} - \frac{c}{2} \trace \left( \mathcal{A}_{t-1}^2 (\Sigma_{p_t} + \mu_{p_t} \mu_{p_t}^\top) \right) > 0,
    \end{align}
    which we analyze further in the three settings (A)--(C) from the theorem in the following:
    \begin{enumerate}[label=(\Alph*), itemsep=0.2em, topsep=0.5em]
        \item
        If $\mathcal{A}_{t-1}=\alpha_{t-1} \unitmatrix$, the inequality in \eqref{eqn:GeneralLowerBound} becomes
        \begin{align*}
            0 < \alpha_{t-1} \mu_{p_t}^\top \mu_{p_t} - \frac{c\,\alpha_{t-1}^2}{2} \trace \left( \Sigma_{p_t} + \mu_{p_t} \mu_{p_t}^\top \right)
            = \alpha_{t-1} \left( \Vert \mu_{p_t} \Vert^2 - \alpha_{t-1} \frac{c}{2} \left(\Vert \mu_{p_t} \Vert^2 + \trace (\Sigma_{p_t}) \right) \right),
        \end{align*}
        which holds  whenever $\alpha_{t-1}$ satisfies the inequality \eqref{eqn:LowerBoundScalar} stated in part (A) of Theorem~\ref{thm:UpperBoundLearningRate}.

        \item
        For $\mathcal{A}_{t-1} = \operatorname{diag}(\alpha_{1, t-1}, \dots, \alpha_{k, t-1}) \succ \zeromatrix$, the inequality in \eqref{eqn:GeneralLowerBound} becomes
        \begin{align*}
            0 &< \mu_{p_t}^\top \mathcal{A}_{t-1} \mu_{p_t} - \frac{c}{2} \trace \left( \mathcal{A}_{t-1}^2 (\Sigma_{p_t} + \mu_{p_t} \mu_{p_t}^\top) \right)\\
            &= \sum_{i=1}^k \alpha_{i,t-1} \mu_{i}^2 - \frac{c}{2} \sum_{i,j=1}^k (\mathcal{A}_{t-1}^2)_{ij} \big(\Sigma_{p_t} + \mu_{p_t} \mu_{p_t}^\top \big)_{ji} \\ 
            &= \sum_{i=1}^k \alpha_{i,t-1}  \mu_{i}^2 - \frac{c}{2} \sum_{i=1}^k (\mathcal{A}_{t-1}^2)_{ii} \big(\Sigma_{p_t} + \mu_{p_t} \mu_{p_t}^\top \big)_{ii} \\
            &= \sum_{i=1}^k \left(\alpha_{i,t-1} \mu_{i}^2 - \frac{c}{2} \alpha_{i,t-1}^2 \big(\sigma_{i}^2 + \mu_{i}^2 \big) \right),
        \end{align*}
        as $\mathcal{A}_{t-1}^2$ is diagonal.
        The positivity in the term above guaranteed if all summands are weakly positive while at least one summand is strictly positive, which immediately gives the condition  \eqref{eqn:LowerBoundDiagonal} and hence shows part (B).

        \item
        For a general learning-rate matrix $\mathcal{A}_{t-1} \succ O_k$, we have that $$\mu_{p_t}^\top \mathcal{A}_{t-1} \mu_{p_t} \geq \lambda_{\min}(\mathcal{A}_{t-1}) \Vert \mu_{p_t} \Vert^2$$ and by a standard trace inequality for positive semi-definite matrices 
        \begin{align*}
            \frac{c}{2} \trace \left( \mathcal{A}_{t-1}^2(\Sigma_{p_t} + \mu_{p_t} \mu_{p_t}^\top) \right) 
            \le \frac{c}{2} \lambda_{\max}(\mathcal{A}_{t-1})^2 \big(\Vert \mu_{p_t} \Vert^2 + \trace(\Sigma_{p_t}) \big).
        \end{align*}
        Combining these two inequalities, we get that 
        \begin{align*}
            &\mu_{p_t}^\top \mathcal{A}_{t-1} \mu_{p_t} - \frac{c}{2} \trace \left( \mathcal{A}_{t-1}^2 (\Sigma_{p_t} + \mu_{p_t} \mu_{p_t}^\top) \right) \\
            &\qquad
            \ge \lambda_{\min}(\mathcal{A}_{t-1}) \Vert \mu_{p_t} \Vert^2 - \lambda_{\max}(\mathcal{A}_{t-1})^2 \frac{c}{2} \big(\Vert \mu_{p_t} \Vert^2 + \trace(\Sigma_{p_t}) \big).
        \end{align*}
        The right-hand side is strictly positive if 
        \begin{align*}
            \frac{\lambda_{\max}(\mathcal{A}_{t-1})^2}{\lambda_{\min}(\mathcal{A}_{t-1})}
            < \frac{2}{c} \, \frac{\Vert \mu_{p_t}\Vert^2}{\Vert \mu_{p_t} \Vert^2 + \trace(\Sigma_{p_t})},
        \end{align*}
        which is the condition \eqref{eqn:LowerBoundPosDef} and hence shows part (C).    \end{enumerate}
\end{proof}

\begin{proof}[Proof of Proposition~\ref{prop:CEVEquiv}]
    Using $\Delta \phi_\kappa(y, \vartheta_\tbu) = \kappa \Delta \phi(y,\vartheta_\tbu)$, we obtain
    \begin{align*}
        &\Delta^{\cev}(\phi_\kappa| \vartheta_\tbu, p_t) \\
        &=\Big( \E_{p_t}[\vartheta_\tpu^\kappa (Y_t)] -  \vartheta_t^\ast \Big)^\top \Omega_{t-1} \Big(  \E_{p_t}[\vartheta_\tpu^\kappa (Y_t)] -  \vartheta_t^\ast \Big)
        - \Big( \vartheta_\tbu - \vartheta_t^\ast \Big)^\top \Omega_{t-1} \Big( \vartheta_\tbu - \vartheta_t^\ast \Big)  \\
        &=  \E_{p_t} \Big[ \vartheta_\tbu - \vartheta_t^\ast + \kappa \Delta \phi(Y_t,\vartheta_\tbu) \Big]^\top \Omega_{t-1}  \E_{p_t} \Big[  \vartheta_\tbu - \vartheta_t^\ast + \kappa \Delta \phi(Y_t,\vartheta_\tbu) \Big] \\
        &\qquad -  \Big( \vartheta_\tbu - \vartheta_t^\ast \Big)^\top \Omega_{t-1} \Big( \vartheta_\tbu - \vartheta_t^\ast \Big) \\
        &= - 2 \kappa \,  \E_{p_t} \left[ \Delta \phi(Y_t,\vartheta_\tbu) \right]^\top \Omega_{t-1} \big(\vartheta_t^\ast -  \vartheta_\tbu \big)\\
        &\qquad + \kappa^2 \,  \E_{p_t} \big[\Delta \phi(Y_t,\vartheta_\tbu) \big]^\top \Omega_{t-1} \E_{p_t} \big[ \Delta \phi(Y_t,\vartheta_\tbu) \big].
    \end{align*}
    Similar to the proof of Theorem~\ref{thm:EKLequivalenceIndepCopy}, imposing $\E_{p_t} \left[\Delta \phi(Y_t,\vartheta_\tbu) \right] < \infty$ implies that for $\kappa$ small enough, the first order term, $- 2  \kappa  \E_{p_t} \left[ \Delta \phi(Y_t,\vartheta_\tbu) \right]^\top \Omega_{t-1} \big( \vartheta_t^\ast -  \vartheta_\tbu \big)$, determines the sign of $\Delta^{\cev}(\phi_\kappa)$.
    Hence, if $\E_{p_t} \big[ \Delta\phi (Y_t,\vartheta_\tbu) \big]^\top \Omega_{t-1} \big(\vartheta_t^\ast -  \vartheta_\tbu \big)$ is non-zero, the equivalence follows.
    Finally, if $\E_{p_t} \big[ \Delta\phi (Y_t,\vartheta_\tbu) \big]^\top \Omega_{t-1} \big(\vartheta_t^\ast -  \vartheta_\tbu \big) = 0$, the expansion shows that $\Delta^{\cev}(\phi_\kappa| \vartheta_\tbu, p_t) \ge 0$ such that the equivalence also holds in these cases.
\end{proof}

\begin{proof}[Proof of Proposition~\ref{prop:MSEEquiv}]
    Using $\Delta \phi_\kappa(y, \vartheta_\tbu) = \kappa \Delta \phi(y,\vartheta_\tbu)$, we obtain
    \begin{align*}
        &\hspace{-0.3cm}\Delta^{\mse}(\phi_\kappa | \vartheta_\tbu, p_t) \\
        &= \E_{p_t} \left[ \Big( \vartheta_\tpu^\kappa (Y_t) -  \vartheta_t^\ast \Big)^\top \Omega_{t-1} \Big( \vartheta_\tpu^\kappa(Y_t) -  \vartheta_t^\ast \Big) \right] - \E_{p_t}  \left[ \Big( \vartheta_\tbu - \vartheta_t^\ast \Big)^\top \Omega_{t-1} \Big( \vartheta_\tbu - \vartheta_t^\ast \Big) \right]  \\
        &= \E_{p_t}  \left[ \Big( \vartheta_\tbu - \vartheta_t^\ast + \kappa \Delta \phi(Y_t,\vartheta_\tbu) \Big)^\top \Omega_{t-1} \Big( \vartheta_\tbu - \vartheta_t^\ast + \kappa \Delta \phi(Y_t,\vartheta_\tbu) \Big) \right] \\
        &\qquad - \E_{p_t} \left[ \Big( \vartheta_\tbu - \vartheta_t^\ast \Big)^\top \Omega_{t-1} \Big( \vartheta_\tbu - \vartheta_t^\ast \Big) \right]  \\
        &= - 2  \kappa \E_{p_t} \left[ \Delta \phi(Y_t,\vartheta_\tbu) \right]^\top  \Omega_{t-1} \big(\vartheta_t^\ast - \vartheta_\tbu \big) + \kappa^2 \E_{p_t} \left[\Delta \phi(Y_t,\vartheta_\tbu)^\top \Omega_{t-1} \Delta \phi(Y_t,\vartheta_\tbu) \right].
    \end{align*}
    Similar to the proofs of Theorem~\ref{thm:EKLequivalenceIndepCopy} and Proposition~\ref{prop:CEVEquiv}, imposing the moment condition $\E \left[ \Vert \Delta \phi(Y_t,\vartheta_\tbu)\Vert^2 \right] < \infty$ implies that for $\kappa$ small enough, the first order term, 
    \begin{align*}
        - 2  \kappa \E_{p_t}  \left[ \Delta \phi(Y_t,\vartheta_\tbu) \right]^\top \Omega_{t-1} \big(\vartheta_t^\ast - \vartheta_\tbu\big),
    \end{align*}
    determines the sign of $\Delta^{\mse}(\phi_\kappa)$.
    Hence, if $\E_{p_t} \big[ \Delta\phi (Y_t,\vartheta_\tbu) \big]^\top  \Omega_{t-1} \big(\vartheta_t^\ast -  \vartheta_\tbu \big)$ is non-zero, the equivalence follows.
    Finally, if $\E_{p_t} \big[ \Delta\phi (Y_t,\vartheta_\tbu) \big]^\top  \Omega_{t-1} \big(\vartheta_t^\ast -  \vartheta_\tbu \big) = 0$, the expansion shows that $\Delta^{\mse}(\phi_\kappa| \vartheta_\tbu, p_t) \ge 0$ such that the equivalence also holds in these cases.
\end{proof}

\begin{proof}[Proof of Proposition~\ref{prop:EGMMEquiv}]
    This proof closely follows the proof of \citet[Proposition~3]{Creal2024GMM} in their online supplement.
    We use the following shorthand notations:
    \begin{align*}
        g(\vartheta) &:= \E_{p_t} \big[ \score(X_t, \vartheta) \big], \qquad  \qquad
        G(\vartheta) := \E_{p_t} \big[ \Hessian(X_t, \vartheta) \big], \\
        \Delta \phi_\kappa(Y_t) &\equiv \Delta \phi_\kappa(Y_t, \vartheta_\tbu) := \vartheta_\tpu^\kappa(Y_t) - \vartheta_\tbu.
    \end{align*}    
    By the path integral version mean-value expansion, we establish 
    \begin{align*}
        g(\vartheta_\tpu^\kappa(y))
        = g(\vartheta_\tbu) + \int_{0}^1 G\big( \widetilde{\vartheta}^\kappa_u(y)\big) \mathrm{d} u \; \Delta \phi_\kappa(y) 
    \end{align*}
    where $\widetilde{\vartheta}^\kappa_u(y) = u \vartheta_\tbu + (1-u) \vartheta^\kappa_\tpu(y)$ moves on the line between $\vartheta_\tbu$ and $\vartheta^\kappa_\tpu(y)$
    Then,
    \begin{align*}
        &\hspace{-0.3cm}\Delta^{\gmm}(\phi_\kappa) \\
        &= \gmm(p_t \Vert f_\tpu) - \gmm(p_t \Vert f_\tbu)  \\
        &= \E_{p_t} \left[ g(\vartheta_\tpu^\kappa(Y_t))^\top \Omega_{t-1} g(\vartheta_\tpu^\kappa(Y_t)) - g(\vartheta_\tbu)^\top \Omega_{t-1} g(\vartheta_\tbu) \right] \\
        &= \E_{p_t} \biggl[ \Big\{ g(\vartheta_\tbu) + \int_{0}^1 G\big( \widetilde{\vartheta}^\kappa_u(Y_t)\big)\, \mathrm{d} u \Delta \phi_\kappa(Y_t) \Big\}^\top \Omega_{t-1}\Big\{ g(\vartheta_\tbu)  \\
        & \qquad + \int_{0}^1 G\big( \widetilde{\vartheta}^\kappa_u(Y_t)\big) \mathrm{d} u \Delta \phi_\kappa(Y_t) \Big\} - g(\vartheta_\tbu)^\top \Omega_{t-1} g(\vartheta_\tbu) \biggr] \\
        &= \E_{p_t} \Big[ 2 g(\vartheta_\tbu)^\top \Omega_{t-1} \int_{0}^1 G\big( \widetilde{\vartheta}^\kappa_u(Y_t)\big) \mathrm{d} u \,\Delta \phi_\kappa(Y_t) \Big] \\
        &\qquad+ \E_{p_t} \Big[\Delta \phi_\kappa(Y_t)^\top \int_{0}^1 G\big( \widetilde{\vartheta}^\kappa_u(Y_t)\big) \mathrm{d} u\, \Omega_{t-1} \int_{0}^1 G\big( \widetilde{\vartheta}^\kappa_u(Y_t)\big) \mathrm{d} u\, \Delta \phi_\kappa(Y_t) \Big] \\
        &= 2 \kappa \, \E_{p_t} \Big[ g(\vartheta_\tbu)^\top \Omega_{t-1} \int_{0}^1 G\big( \widetilde{\vartheta}^\kappa_u(Y_t)\big) \mathrm{d} u\, \Delta \phi(Y_t) \Big] \\
        &\qquad+ \kappa^2 \, \E_{p_t} \Big[\Delta \phi(Y_t)^\top \int_{0}^1 G\big( \widetilde{\vartheta}^\kappa_u(Y_t)\big) \mathrm{d} u\, \Omega_{t-1} \int_{0}^1 G\big( \widetilde{\vartheta}^\kappa_u(Y_t)\big) \mathrm{d} u\, \Delta \phi(Y_t) \Big] \\
        &= 2 \kappa \, \E_{p_t} \Big[ g(\vartheta_\tbu)^\top \Omega_{t-1} \int_{0}^1 G\big( \widetilde{\vartheta}^\kappa_u(Y_t)\big) \mathrm{d} u \,\Delta \phi(Y_t) \Big] + \mathcal{O}(\kappa^2),
    \end{align*}
    as $\sup_{\vartheta \in \Theta} \Vert  G(\vartheta) \Vert \le c < \infty$
    and $\E_{p_t} \big[ \Vert \Delta \phi(Y_t) \Vert^2 \big] < \infty$ by assumption.

Next consider
     \begin{align*}
        G\big( \widetilde{\vartheta}^\kappa_u(y)\big)&= G\big(\vartheta_\tbu\big) + \int_0^1 \sum_{j=1}^k \left. \frac{\partial G\big(\vartheta\big)}{\partial \vartheta_j} \right|_{\vartheta = \widetilde{\vartheta}_{r,u}^\kappa(y)} \left(  \widetilde{\vartheta}^\kappa_u(y) - \vartheta_\tbu \right)_j \mathrm{d} r \\
        &=G\big(\vartheta_\tbu\big) +\int_0^1   \bar{G}\big(  \widetilde{\vartheta}_{r,u}^\kappa(y)\big) \times_3 \left(  \widetilde{\vartheta}^\kappa_u(y) - \vartheta_\tbu \right) \mathrm{d} r 
    \end{align*}
    where $ \widetilde{\vartheta}_{r,u}^\kappa(y)= r \vartheta_\tbu + (1-r) \widetilde{\vartheta}^\kappa_u(y)$ moves on the line between $\vartheta_\tbu$ and $\widetilde{\vartheta}^\kappa_u(y)$. Moreover,  $\bar{G}(\vartheta)$ denotes the $k\times k \times k$ tensor with elements 
\begin{align*}
    [\bar{G}(\vartheta)]_{i_1,i_2,j}= \frac{\partial [G(\vartheta)]_{i_1,i_2}}{\partial \vartheta_j}, \qquad  i_1, i_2, j=1, \ldots, k,
\end{align*}
which corresponds to the tensor that is assumed to be uniformly bounded in $\Theta$ by \citet[Ass. A.4]{Creal2024GMM}, and $\times_3$ denotes the 3-mode vector product of a tensor \citep{kolda_tensor_2009}.
Then,
     \begin{align*}
        &\Delta^{\gmm}(\phi_\kappa) \\
        &= 2 \kappa \, \E_{p_t} \Big[ g(\vartheta_\tbu)^\top \Omega_{t-1} \int_{0}^1 G\big( \widetilde{\vartheta}^\kappa_u(Y_t)\big) \mathrm{d} u \Delta \phi(Y_t) \Big]  + \mathcal{O}(\kappa^2) \\
        &= 2 \kappa \, \E_{p_t} \Big[ g(\vartheta_\tbu)^\top \Omega_{t-1}G\big(\vartheta_\tbu\big)  \Delta \phi(Y_t) \Big]  + \mathcal{O}(\kappa^2) \\
        &\qquad+ 2 \kappa \, \E_{p_t} \left[ g(\vartheta_\tbu)^\top \Omega_{t-1} \int_0^1 \int_0^1   \bar{G} \big(  \widetilde{\vartheta}_{r,u}^\kappa(y)\big) \times_3  \left(  \widetilde{\vartheta}^\kappa_u(y) - \vartheta_\tbu \right) \mathrm{d} r \mathrm{d} u \Delta \phi(Y_t) \right] \\
         &= 2 \kappa \, g(\vartheta_\tbu)  \Omega_{t-1} G(\vartheta_\tbu) \E_{p_t} \big[ \Delta \phi(Y_t)\big] + \mathcal{O}(\kappa^2),
    \end{align*}
where the order $\mathcal{O}(\kappa^2)$ can be verified as follows. 
By assumption, there exists a $\tilde{c}>0$ such that $\sup_{\vartheta \in \Theta} \big\Vert (\partial/\partial \vartheta_j) G(\vartheta) \big\Vert \leq \tilde{c}<\infty$ for all $j=1,\dots,k$. 
Hence,
\begin{align*}
&\left \vert \E_{p_t} \left[ g(\vartheta_\tbu)^\top \Omega_{t-1} \int_0^1 \int_0^1  \bar{G}\big( \widetilde{\vartheta}_{r,u}^\kappa(Y_t)\big) \times_3 \left(  \widetilde{\vartheta}^\kappa_u(Y_t) - \vartheta_\tbu \right) \mathrm{d} r  \mathrm{d} u \Delta \phi(Y_t) \right]\right \vert \\
   &\qquad \leq  \E_{p_t} \left[  \int_0^1\left\vert g(\vartheta_\tbu)^\top \Omega_{t-1} \int_0^1  \bar{G}\big( \widetilde{\vartheta}_{r,u}^\kappa(Y_t)\big) \times_3 \left(  \widetilde{\vartheta}^\kappa_u(Y_t)- \vartheta_\tbu \right) \mathrm{d} r \Delta  \phi(Y_t) \right \vert  \mathrm{d} u  \right] \\
&\qquad \leq  \int_0^1 \E_{p_t} \left[  \left\vert g(\vartheta_\tbu)^\top \Omega_{t-1} \int_0^1  \bar{G}\big( \widetilde{\vartheta}_{r,u}^\kappa(Y_t)\big)\times_3  \left(  \widetilde{\vartheta}^\kappa_u(Y_t)- \vartheta_\tbu \right) \mathrm{d} r \Delta  \phi(Y_t) \right \vert    \right]\mathrm{d} u \\
&\qquad \leq  \int_0^1 \E_{p_t} \left[  \left\Vert g(\vartheta_\tbu)^\top \Omega_{t-1} \right \Vert  \left\Vert\int_0^1  \bar{G}\big( \widetilde{\vartheta}_{r,u}^\kappa(Y_t)\big) \times_3  \left(  \widetilde{\vartheta}^\kappa_u(Y_t)- \vartheta_\tbu \right) \mathrm{d} r \right \Vert \left\Vert \Delta  \phi(Y_t) \right \Vert    \right]\mathrm{d} u \\
&\qquad \leq \bar{c}_{t\vert t-1} \int_{0}^1  \E_{p_t} \left[\left\Vert\int_0^1    \sum_{j=1}^k \left. \frac{\partial G\big(\vartheta\big)}{\partial \vartheta_j} \right|_{\vartheta = \widetilde{\vartheta}_{r,u}^\kappa(y)} \left(  \widetilde{\vartheta}^\kappa_u(y) - \vartheta_\tbu \right)_j \mathrm{d} r \right \Vert \left\Vert \Delta  \phi(Y_t) \right \Vert    \right]\mathrm{d} u \\
&\qquad \leq \bar{c}_{t\vert t-1} \int_{0}^1  \E_{p_t} \left[\int_0^1 \sum_{j=1}^k \left \Vert \left. \frac{\partial G\big(\vartheta\big)}{\partial \vartheta_j} \right|_{\vartheta = \widetilde{\vartheta}_{r,u}^\kappa(Y_t)} \right \Vert \left \vert\left(  \widetilde{\vartheta}^\kappa_u(Y_t) - \vartheta_\tbu \right)_j \right \vert \mathrm{d} r  \left\Vert \Delta  \phi(Y_t) \right \Vert    \right]\mathrm{d} u \\
&\qquad \leq \frac{\tilde{c}}{2} \bar{c}_{t\vert t-1}  \int_0^1   \E_{p_t} \left[\sum_{j=1}^k \left \vert\left(  \widetilde{\vartheta}^\kappa_u(Y_t) - \vartheta_\tbu \right)_j \right \vert \left\Vert \Delta  \phi(Y_t) \right \Vert    \right]\mathrm{d} u \\
&\qquad \leq \kappa \frac{\tilde{c}   }{4} \bar{c}_{t\vert t-1}  \E_{p_t} \left[\sum_{j=1}^k \left \vert\left(  \Delta  \phi(Y_t) \right)_j \right \vert \left\Vert \Delta  \phi(Y_t) \right \Vert    \right] \\
&\qquad \leq  \kappa \frac{\tilde{c}   \sqrt{k}}{4} \bar{c}_{t\vert t-1} \E_{p_t} \left[\left\Vert \Delta  \phi(Y_t) \right \Vert^2    \right] \\
&\qquad= \mathcal{O}(\kappa),
\end{align*}
as $0 < \bar{c}_{t\vert t-1} := \left\Vert g(\vartheta_\tbu) \right \Vert  \left\Vert\Omega_{t-1} \right \Vert < \infty$ by assumption, and similarly $\E_{p_t} \left[  \left\Vert \Delta  \phi(Y_t) \right \Vert^2    \right]<\infty$. 
The inequality $\sum_{j=1}^k \left \vert\left(  \Delta  \phi(Y_t) \right)_j \right \vert  = \vert \Delta  \phi(Y_t)  \vert^\top\iota_k \leq \Vert \iota_k \Vert \Vert \Delta  \phi(Y_t) \Vert = \sqrt{k} \Vert \Delta  \phi(Y_t) \Vert$ is a direct consequence of Cauchy-Schwarz, where $\iota_k$ denotes a $k$-dimensional vector of ones, and $\vert \cdot \vert$ the elementwise absolute value.

    Similar to the proof of Theorem~\ref{thm:EKLequivalenceIndepCopy}, for $\kappa$ small enough, the first order term,  
    \begin{align*}
        2 \kappa \, g(\vartheta_\tbu)^\top \Omega_{t-1} G(\vartheta_\tbu) \E_{p_t} \big[ \Delta \phi(Y_t)\big],
    \end{align*}
    determines the sign of $\Delta^{\gmm}(\phi_\kappa)$ such that the result follows.
    Further notice that if the product $g(\vartheta_\tbu)^\top \Omega_{t-1} G(\vartheta_\tbu) \E_{p_t} \big[ \Delta \phi(Y_t)\big]$ is zero, the sign of $\Delta^{\gmm}(\phi_\kappa)$ depends on the $\mathcal{O}(\kappa^2)$ terms, whose (combined) sign is unclear.
\end{proof}

\section{Details for the examples in Section~\ref{sec:Example}}
\label{sec:DetailsExamples}

Table~\ref{tab:DGPs} provides additional details for the examples discussed in Section~\ref{sec:Example}.
In the following as well as in Table~\ref{tab:DGPs}, we use the subscript $t$ in $\vartheta_t$ and $\lambda_t$ to clearly distinguish the time-varying parameter from the remaining static parameters.
We now elaborate on three illustrative cases: the exponential intensity model, the negative binomial count model, and the Student’s~$t$ location model.

\begin{sidewaystable}[p]
\caption{\label{tab:DGPs} Full details for models densities in Table~\ref{tab2}.\label{tab3}}
\resizebox{1.0\textwidth}{!}{
\begin{footnotesize}
\begin{threeparttable}
\begin{tabular}{l@{\hspace{0.15cm}}l@{\hspace{0.15cm}}c@{\hspace{0.15cm}}c@{\hspace{0.15cm}}c@{\hspace{0.15cm}}c@{\hspace{0.15cm}}c@{\hspace{0.15cm}}c}
  \toprule
\multicolumn{2}{c}{\bf{Postulated model}} & \bf{Link function}  & \bf{Density} & \bf{Score} & \bf{Hessian}& \multicolumn{2}{c}{\bf{Expected Hessian}}
\\
Type & Distribution &   & $f({y}_t|\vartheta_t)$ & $\score(y_t,\vartheta_t)$ & $\Hessian(y_t,\vartheta_t)$ & Condition on $\Pm$ & Range
\\
  \cmidrule(r{10pt}){1-2}  \cmidrule(lr){3-3} \cmidrule(lr){4-4} \cmidrule(lr){5-5} \cmidrule(lr){6-6} \cmidrule(lr){7-7} \cmidrule(lr){8-8}

  Count 
& Poisson 
& $\lambda_t=\exp(\vartheta_t)$ 
&$\displaystyle\lambda_t^{y_t}\, \exp(-\lambda_t)/y_t! $ 
& $y_t-\lambda_t$ 
& $-\lambda_t$ 
& $-$
& $(-\infty,0)$
\\
Count  
& Neg.~Binomial
& $\lambda_t=\exp(\vartheta_t)$ 
& $\displaystyle \frac{\Gamma(\xi+y_t)\left(\frac{\xi}{\xi+\lambda_t}\right)^\xi\left(\frac{\lambda_t}{\xi+\lambda_t}\right)^{y_t}}{\Gamma(\xi)\Gamma(y_t+1)}$ 
& $\displaystyle y_t-\frac{\lambda_t(\xi+y_t)}{\xi+\lambda_t}$ 
& $\displaystyle- \frac{\xi \lambda_t(\xi+y_t)}{(\xi+\lambda_t)^2}$ 
& $\E_{p_t}Y_t < \infty$
& $\left[- \frac{\xi + \mathbb{E}_{p_t}[Y_t]}{4}, 0 \right)$

\\
Intensity 
& Exponential 
& $\lambda_t=\exp(\vartheta_t)$ 
& $\displaystyle\lambda_t\, \exp(-\lambda_t y_t)$ &  $1-\lambda_t\,y_t$ 
& $-y_t \lambda_t$ 
& $0 < \E_{p_t}Y_t<\infty$
& $(-\infty,0)$

\\
Duration 
& Gamma 
& $\beta_t=\exp(\vartheta_t)$ 
& $\displaystyle
\frac{y_t^{\xi-1}\exp(-y_t/\beta_t)}{\Gamma(\xi)\beta_t^\xi }$ 
& $\displaystyle \frac{y_t}{\beta_t}-\xi$
& $\displaystyle -\frac{y_t}{\beta_t}$ 
& $0 < \E_{p_t}Y_t<\infty$
& $(-\infty,0)$

\\
Duration 
& Weibull  
& $\beta_t=\exp(\vartheta_t)$ 
& $\displaystyle\frac{\xi\, \left(y_t/\beta_t \right)^{\xi-1} }{\beta_t \exp\{(y_t/\beta_t)^\xi\} }
$
& $\displaystyle \xi\left(\frac{y_t}{\beta_t}\right)^\xi-\xi$
& $\displaystyle -\xi^2 \left(\frac{y_t}{\beta_t}\right)^\xi$
& $0 < \E_{p_t}Y_t^\xi<\infty$
& $(-\infty,0)$

\\
Volatility 
& Gaussian 
& $\sigma^2_t=\exp(\vartheta_t)$ 
&$ \displaystyle
\frac{\exp\{-y_t^2/(2\sigma_t^2)\}}{ \{2\pi \sigma_t^2\}^{1/2} }$ & $\displaystyle \frac{y_t^2}{2\sigma_t^2}-\frac{1}{2}$
&$\displaystyle -\frac{y_t^2}{2\sigma_t^2}$ 
& $0 < \E_{p_t}Y_t^2<\infty$
& $(-\infty,0)$

\\
Volatility  
& Student's $t$
& $\sigma^2_t=\exp(\vartheta_t)$ 
&  $\displaystyle\frac{\Gamma\left(\frac{\nu+1}{2}\right)\left(1+\frac{y_t^2}{(\nu-2)\sigma_t^2}\right)^{-\frac{\nu+1}{2}}}{\sqrt{(\nu-2)\pi}\Gamma\left(\nu/2\right)\sigma_t}$ 
&
$\displaystyle \frac{\omega_t \, y_t^2}{2\sigma_t^2}-\frac{1}{2}$
& $\displaystyle -\frac{\nu-2}{\nu+1}\,\frac{\omega_t^2\,y_t^2}{2\sigma_t^2}$
& $0<\E_{p_t}Y_t^2$
& $\displaystyle \left[-\frac{\nu+1}{8},0\right)$ 
\\
&&&& $\displaystyle \omega_t:=\frac{\nu+1}{\nu-2+y_t^2/\sigma_t^2}$ & &

\\
Dependence & Gaussian & $\displaystyle\rho_t=\frac{1-\exp(-\vartheta_t)}{1+\exp(-\vartheta_t)}$ & $\displaystyle\frac{\exp\left\{-\frac{y_{1t}^2+y_{2t}^2-2\rho_t y_{1t}y_{2t} }{2(1-\rho_t^2)} \right\}}{2\pi \sqrt{1-\rho_t^2}}$ & $
\displaystyle \frac{\rho_t}{2}+
\frac{1}{2}\frac{ z_{1t}\,z_{2t}}{1-\rho_t^2} $
& $\displaystyle \frac{1-\rho_t^2}{4}-\frac{1}{4}\frac{z_{1t}^2+z_{2t}^2}{1-\rho_t^2}$
& $\E_{p_t}Y_{jt}^2<\infty$
& $\displaystyle \left(-\infty,\frac{1}{4} \right]$
\\
&& $\displaystyle$ && $z_{1t}:=y_{1t}-\rho_t y_{2t}$
\\
&&&& $z_{2t}:=y_{2t}-\rho_t y_{1t}$
\\
Dependence & Student's $t$ & $\displaystyle\rho_t=\frac{1-\exp(-\vartheta_t)}{1+\exp(-\vartheta_t)}$ & $\displaystyle\frac{\nu \left(1+\frac{y_{1t}^2+y_{2t}^2-2\rho_t y_{1t}y_{2t} }{(\nu-2)(1-\rho_t^2)}\right)^{-\frac{\nu+2}{2}}}{2\pi(\nu-2) \sqrt{1-\rho_t^2}}$ & $\displaystyle \frac{\rho_t}{2}+
\frac{\omega_t}{2}\frac{ z_{1t}\,z_{2t}}{1-\rho_t^2} $
& $\displaystyle \frac{1-\rho_t^2}{4}-\frac{\omega_t}{4}\frac{z_{1t}^2+z_{2t}^2}{1-\rho_t^2}+\frac{1}{2}\frac{\omega_t^2}{\nu+2}\frac{z_{1t}^2 \, z_{2t}^2}{(1-\rho_t^2)^2}$
& $-$
& $\displaystyle \left[-\frac{\nu+1}{4},\frac{1}{4} \right]$
\\
&& $\displaystyle$ && $z_{1t}:=y_{1t}-\rho_t y_{2t}$ &  \multirow{2}{*}{$\displaystyle \omega_t:=\frac{\nu+2}{\nu-2+\frac{y_{1t}^2+y_{2t}^2-2\rho_t y_{1t}y_{2t} }{1-\rho_t^2}}$}
\\
&&&& $z_{2t}:=y_{2t}-\rho_t y_{1t}$
\\
\\
Local level & Gaussian
& $\mu_t=\vartheta_t$ 
& $\displaystyle\frac{1}{\sqrt{2\pi}\sigma}\exp{\left\{-\frac{(y_t - \mu_t)^2}{2\sigma^2}\right\}}$ 
& $\displaystyle \frac{y_t - \mu_t}{\sigma^2}$ 
& $-\frac{1}{\sigma^2}$ 
& $-$
& $\displaystyle \left[-\frac{1}{\sigma^2}  , -\frac{1}{\sigma^2}\right]$
\\
Local level & Student's \emph{t} 
& $\mu_t=\vartheta_t$ 
& $\displaystyle\frac{\Gamma\left(\frac{\nu+1}{2}\right)\left(1+\frac{(y_t-\mu_t)^2}{\nu\,\sigma^2}\right)^{-\frac{\nu+1}{2}}}{\sqrt{\nu\,\pi}\Gamma\left(\frac{\nu}{2}\right)\sigma}$ 
& $\displaystyle \frac{\nu+1}{\nu \, \sigma^2}\frac{y_t - \mu_t}{1+\frac{(y_t - \mu_t)^2}{\nu \,\sigma^2}}$ 
& $\displaystyle \frac{\nu+1}{\nu \, \sigma^2}  \left(1+\frac{(y_t - \mu_t)^2}{\nu \, \sigma^2}\right)^{-2}
\left(\frac{(y_t - \mu_t)^2}{\nu \, \sigma^2} - 1\right)
$
& $-$
& $\displaystyle \left[ - \frac{(\nu+1)}{\nu \, \sigma^2} , \frac{1}{8} \frac{\nu+1}{\nu \, \sigma^2}\right]$
\\
  \bottomrule
\end{tabular}
\begin{tablenotes}
\spacingset{1}
\item {NOTE:} The table contains eleven model densities and link functions, the first nine of which are adapted from \citet{koopman2016predicting}. 
For each model, we report the link function, density, score and Hessian matrix.
The last two columns give the range of the expected Hessian, together with the imposed condition on the admissible truths $p_t \in \Pm$ required for the range.
To facilitate comparison with \cite{koopman2016predicting}, the table largely preserves their parameter notation, at the cost of some internal inconsistency because several symbols are used in the main text for different quantities. 
We distinguish static from time-varying parameters by using a subscript $t$ for the latter.
The restrictions on the static parameters are given by $\xi>0$ and
$\sigma>0$ with $\nu>0$ for the Student's $t$ distribution in the local-level model and $\nu>2$ for all other Student's $t$ distributions.
\end{tablenotes}
\end{threeparttable}
\end{footnotesize}
}
\end{sidewaystable}

First, consider the exponential intensity model, for which the Hessian equals $\Hessian(y_t,\vartheta_t) = -y_t\lambda_t = -y_t\exp(\vartheta_t)$. For $\mathcal{Y}=[0,\infty)$ (a weak restriction for intensities) and $\Theta=\R$, the Hessian is weakly negative, but it is not uniformly bounded in $\vartheta_t$. Imposing $0<\E_{p_t}[Y_t]<\infty$ for all $p_t\in\Pm$ amounts to a mild moment restriction on $\Pm$. It rules out the degenerate case $Y_t=0$ a.s. and guarantees that the expected Hessian satisfies $\E_{p_t}[\Hessian(Y_t,\vartheta_t)] = -\exp(\vartheta_t)\E_{p_t}[Y_t] < 0$. In addition, $\E_{p_t}[\Hessian(Y_t,\vartheta_t)]$ is locally, but not globally, bounded in $\vartheta_t$. Consequently, Assumptions~\ref{ass:HLB} and \ref{ass:HN} hold, whereas Assumptions~\ref{ass:HB} and \ref{ass:HBNT} fail, as reflected by the check and cross marks in Table~\ref{tab2}. If the moment condition were weakened to $0\leq \E_{p_t}[Y_t]<\infty$, then Assumption~\ref{ass:HLB} would still hold, but Assumption~\ref{ass:HN} would fail when $\E_{p_t}[Y_t] = 0$. 
For the EKL measure, only boundedness in expectation is important, whereas the sign is not.

Second, for the negative binomial count model with $\mathcal{Y} = [0, \infty)$ and $\Theta = \R$, the Hessian is $\Hessian(y_t,\vartheta_t) = -\frac{\xi \lambda_t(\xi+y_t)}{(\xi+\lambda_t)^2}$, where $\lambda_t = \exp(\vartheta_t) \in (0,\infty)$ and $\xi > 0$ is a static parameter.
The Hessian is strictly negative such that Assumption~\ref{ass:HN} holds. 
It is furthermore bounded in $\lambda_t$ due to the denominator term, but unbounded in $y_t$.
Imposing the mild condition $\E_{p_t}[Y_t]<\infty$ for all $p_t \in \Pm$, the \emph{expected} Hessian, 
$\E_{p_t} [\Hessian(Y_t,\vartheta_t)] = -\frac{\xi \lambda_t(\xi+\E_{p_t}[Y_t])}{(\xi+\lambda_t)^2}$,
can be bounded between $- \frac{\xi + \mathbb{E}_{p_t}[Y_t]}{4}$ and zero, such that Assumption~\ref{ass:HN} (and hence also Assumption~\ref{ass:HLB}) holds. 
Finally, Assumption~\ref{ass:HBNT} holds as the derivative with respect to $\lambda_t$ of the expected Hessian is $- \frac{\xi ( \xi + \E_{p_t}[Y_t])\lambda_t (\xi - \lambda_t)}{(\xi + \lambda_t)^3}$, which is uniformly bounded in $\lambda_t$.

Third, for the Student's $t$ location model with $\mathcal{Y} = \Theta = \R$, the Hessian is
\begin{align}
    \label{eqn:HessenStudentstLevel}
    \Hessian(y_t,\vartheta_t) = \frac{\nu+1}{\nu \, \sigma^2}  \left(1+\frac{(y_t - \mu_t)^2}{\nu \, \sigma^2}\right)^{-2}
    \left(\frac{(y_t - \mu_t)^2}{\nu \, \sigma^2} - 1\right),
\end{align}
which can be bounded between $- \frac{(\nu+1)}{\nu \, \sigma^2}$ and $\frac{1}{8} \frac{\nu+1}{\nu \, \sigma^2}$, such that Assumptions~\ref{ass:HB} and \ref{ass:HLB} hold for all $p_t \in \Pm_0$; hence, we can choose $\Pm = \Pm_0$.
However, the Hessian in \eqref{eqn:HessenStudentstLevel} is negative only if $\vert y_t - \mu_t \vert < \sqrt{\nu} \sigma$.
This imposes that the observation $y_t$ must be sufficiently close to $\mu_t$, a condition that is impossible to maintain in practice.
Hence, Assumptions~\ref{ass:HN} and \ref{ass:HBNT} generally fail to hold without restrictions on the admissible $p_t \in \Pm$.\footnote{Notably, \citet[App.~B4]{lange2022robust} shows that SD updates based on a (possibly misspecified) Student’s~$t$ location model are MSE reducing by relying on a different proof technique that avoids Hessians as in Assumption~\ref{ass:HN} altogether.}

For the Student's $t$ location model, we ask whether one can impose reasonable and readily interpretable conditions on each $p_t \in \Pm$ that guarantee Assumption~\ref{ass:HN}, and hence ensure that the expected Hessian is negative---at least for a \emph{restricted} class of true densities. If so, the equivalence statements in Propositions~\ref{prop:CEVEquiv}--\ref{prop:EGMMEquiv} would still hold, albeit only for a restricted family of distributions $\Pm$. To this end, recall that for the Hessian in \eqref{eqn:HessenStudentstLevel} to be negative for all $y_t \in \R$ and $\mu_t \in \R$, it is necessary that $\vert y_t-\mu_t\vert < \sqrt{\nu}\,\sigma$.

Although it would suffice for the Hessian to be negative \emph{in expectation}, it is unclear how to derive tractable sufficient conditions on $p_t$. The obstacle is that the factor
$
\left(1+\frac{(y_t-\mu_t)^2}{\nu\,\sigma^2}\right)^{-2}
$
in \eqref{eqn:HessenStudentstLevel} depends on $y_t$ and therefore affects the expectation in a non-constant way. Intuitively, however, any $p_t \in \Pm$ would need to place sufficient probability mass on the event $\vert Y_t-\mu_t\vert < \sqrt{\nu}\,\sigma$. This is difficult to ensure in practice, especially when the researcher's parameter prediction $\mu_t$ is far from the true conditional mean of $Y_t$. Relative to the mild moment restrictions imposed on $\Pm$ in Table~\ref{tab:DGPs}, such a requirement would constitute a severe and practically restrictive condition.

The example of the Student's $t$ location model extends to \emph{any} logarithmic density that is not globally concave in the time-varying parameter: in such cases, it seems infeasible to restore Assumption~\ref{ass:HN} without imposing unduly strong restrictions on the true density $p_t$. By contrast, Assumption~\ref{ass:HLB} typically holds for all $p_t$ or under mild moment conditions.

\section{The conditioning set $\Fm_{t-1}$ and state-space models} 
\label{app:statespace} \label{app:statespaceDGP}

In Section~\ref{sec:Preliminaries}, we define the truth $p_t$ as the density corresponding to the law of $Y_t$ given the flexible information set $\Fm_{t-1}$, which must satisfy $\Fm \supseteq \Fm_{t-1} \supseteq \sigma(Y_{s}; s \le t-1)$
and that $Y_t$ is not measurable with respect to $\Fm_{t-1}$.
This flexibility generalizes upon, e.g., \citet[p.~3]{Gorgi2023} and \citet[p.~6]{Creal2024GMM}, who fix the truth $p_t$ to be conditional on $\sigma(Y_{s}; s \le t-1)$.

Allowing for a user-chosen information set $\Fm_{t-1}$ leads to more transparent interpretations and does not introduce any mathematical complications.
In particular, the more general information set $\Fm_{t-1}$ is permitted to directly include the latent state in state-space models.
Below, we illustrate how the interpretation of the true density $p_t$ differs depending on whether or not the latent state is included in the conditioning set.

Let $\lambda_t$ denote the latent state (that affects $Y_t$), which remains random, even conditional on, for instance, $\sigma(Y_{s}; s \le t-1)$. 
We consider two possibilities for the choice of $\Fm_{t-1}$:
\begin{enumerate}
    \item[(i)] ${\Fm}_{t-1}^{(\mathrm{i})}  = \sigma( \lambda_{s}, Y_{s-1}; s \le t)$
    \item[(ii)] ${\Fm}_{t-1}^{(\mathrm{ii})}  =  \sigma(Y_{s-1}; s \le t)$.
\end{enumerate}
To illustrate, we focus on the linear-state model given by \citet[Eq.~(3.1)]{durbin2012time}, 
\begin{align}
\begin{aligned}
    \label{eqn:LinearSS}
    Y_t &= Z_t \lambda_t +  \varepsilon_t, \qquad 
    \varepsilon_t \stackrel{i.i.d.}{\sim}  \Nm(0,H_t), \\
    \lambda_{t+1} &= T_t \lambda_{t} + R_t \eta_{t}, \qquad 
    \eta_t  \stackrel{i.i.d.}{\sim} \Nm(0,Q_t),
\end{aligned}
\end{align}
with $\lambda_1 \sim \Nm(\bar{a}_1, \bar{P}_1)$, where $\bar{a}_1$ and $\bar{P}_1$ are assumed to be known.
The matrices $Z_t$, $H_t$, $T_t$, $R_t$ and $Q_t$ are also non-stochastic and assumed to be known.
As noted by \citet[p.~82]{durbin2012time}, it is straightforward to show that
 \begin{align*}
 Y_t \mid {\Fm}_{t-1}^{(\mathrm{i})}  \sim p_t(\cdot) = p_y(\cdot \mid \lambda_t),
  \end{align*}
where $p_y(\cdot \mid \lambda_t)$ denotes the density of the conditional law of $Y_t \mid \lambda_t$.
Using \eqref{eqn:LinearSS}, the conditional density further reduces to $p_y(\cdot \mid \lambda_t)=\Nm(Z_t\lambda_t,H_t)$. 
Similarly, the conditional law of $\lambda_{t+1} \mid \sigma({\Fm}_{t-1}^{(\mathrm{i})} ,Y_t)$ reduces to the conditional law of $\lambda_{t+1} \mid \lambda_t$, which density is denoted by  $p_\lambda(\cdot \mid \lambda_{t})$.

Under choice (ii) of the conditioning set, the conditional density of $Y_t \mid {\Fm}_{t-1}^{(\mathrm{ii})}$ is obtained by integrating out the latent state,
 \begin{align}
       Y_t \mid {\Fm}_{t-1}^{(\mathrm{ii})} \sim p_t(\cdot) 
       =  \int \int p_y(\cdot \mid \lambda_t) \; p_\lambda(\lambda_t \mid \lambda_{t-1}) \; q_\lambda(\lambda_{t-1} \mid  {\Fm}_{t-1}^{(\mathrm{ii})}) \; \mathrm{d}\lambda_t \; \mathrm{d}\lambda_{t-1},
    \end{align}
where $q_\lambda( \cdot \mid  {\Fm}_{t-1}^{(\mathrm{ii})})$ denotes the true conditional density of $\lambda_{t-1}$ given ${\Fm}_{t-1}^{(\mathrm{ii})}$.

An analytical expression for $q_\lambda(\lambda_{t-1} \mid  {\Fm}_{t-1}^{(\mathrm{ii})})$ is not immediately available, which may be viewed as a drawback of conditioning on $ {\Fm}_{t-1}^{(\mathrm{ii})}$ alone.
Nevertheless, under \eqref{eqn:LinearSS}, \cite{durbin2012time} show that $q_\lambda(\lambda_{t-1} \mid  {\Fm}_{t-1}^{(\mathrm{ii})})$ is Gaussian, with mean and covariance equal to the filtered state estimates produced by the Kalman filter at time $t-1$.

\section{The Kalman filter for the local-level model}
\label{sec:KalmanFilter}

We now consider a special case of 
the linear state-space model \eqref{eqn:LinearSS}, known as the local-level model \citep[Ch.~2]{durbin2012time}:
\begin{align}
    Y_t &= \lambda_t + \eps_t, \qquad \eps_t \stackrel{i.i.d.}{\sim}  \ \Nm(0,\sigma_\eps^2)\label{eq:obs density} \\
    \lambda_{t+1} &= \lambda_t + \eta_t, \qquad \eta_t  \stackrel{i.i.d.}{\sim} \ \Nm(0,\sigma_\eta^2) \\
    \lambda_1 &\sim \ \Nm (\bar{a}_1, \bar{P}_1),
\end{align}
for which we assume that $\sigma_\varepsilon^2>0$, $\sigma_\eta^2>0$, $\bar{a}_1\in \R$ and $\bar{P}_1>0$ are known constants, and
 $\eps_t$ and $\eta_t$ are mutually independent and independent of $a_1$.

Let $\Fm_{t-1} = \sigma(Y_{s}; s \le t-1)$. Define $\bar{a}_{t\vert t-1}$ and $\bar{P}_{t\vert t-1}$ to be such that $\lambda_t \mid \Fm_{t-1} $ is distributed as $ \Nm(\bar{a}_{t\vert t-1}, \bar{P}_{t\vert t-1})$ and $\bar{P}_{t\vert t}$ such that $\lambda_t \mid \Fm_{t} \sim \Nm(\bar{a}_{t\vert t}, \bar{P}_{t\vert t})$. The associated recursions in \citet[Eq.~2.11--14]{durbin2012time} read 
\begin{align}
    \bar{a}_{t\vert t} &= \bar{a}_{t\vert t-1} + \frac{\bar{P}_{t\vert t-1}}{\bar{P}_{t\vert t-1}+\sigma_\varepsilon^2} (y_t - \bar{a}_{t\vert t-1} )\label{eq:KF1}\\
    \bar{P}_{t\vert t} &= \frac{\bar{P}_{t\vert t-1} \sigma_\varepsilon^2}{\bar{P}_{t\vert t-1} + \sigma_\varepsilon^2}\label{eq:KF2}\\
    \bar{a}_{t+1\vert t} &= \bar{a}_{t\vert t}  \\
    \bar{P}_{t+1\vert t} &= \bar{P}_{t\vert t} + \sigma_\eta^2= \frac{\bar{P}_{t\vert t-1} \sigma_\varepsilon^2}{\bar{P}_{t\vert t-1} + \sigma_\varepsilon^2} + \sigma_\eta^2,
\end{align}
known as the Kalman filter. 

Here we show that Kalman's level update can also be viewed as a score-driven update. Let the postulated density be Gaussian with unknown location and variance $\sigma_\varepsilon^2$; this matches the true observation density~\eqref{eq:obs density}. The score with respect to the unknown location, evaluated at the prediction $\bar{a}_{t|t-1}$, reads $(y_t-\bar{a}_{t|t-1})/\sigma_\varepsilon^2$. Further, if we take $A=1$ and
 $\scaling=\bar{P}_{t\vert t}$, which is $\mathcal{F}_{t-1}$-measurable as can be seen from~\eqref{eq:KF2}, then we obtain the scaled score-driven update as
\begin{equation}
    \bar{a}_{t\vert t} = \bar{a}_{t\vert t-1} + \underbrace{\bar{P}_{t\vert t}}_{=\scaling}\,\underbrace{\frac{y_t - \bar{a}_{t\vert t-1} }{\sigma_\varepsilon^2}}_{=\text{score}} =
    \bar{a}_{t\vert t-1}  + \frac{\bar{P}_{t\vert t-1} \sigma_\varepsilon^2}{\bar{P}_{t\vert t-1} + \sigma_\varepsilon^2} \frac{y_t - \bar{a}_{t\vert t-1} }{\sigma_\varepsilon^2},
\end{equation}
where the second equality uses~\eqref{eq:KF2}. The final result matches~\eqref{eq:KF1}, such that the level update of the Kalman filter is a score-driven update with a specific scaling matrix. 
As the Kalman filter level update moves the parameter in the direction of the score, its downscaled version is EKL reducing. The only exception is when $\bar{a}_{t|t-1}$ coincides with the pseudo-true location at time $t$, in which case no further improvement is possible.

\section{EKL reduction for implicit score-driven updates}
\label{app:ISD}

Here we show that implicit score-driven (ISD; \citealp{lange2022robust}) updates reduce the expected Kullback-Leibler (EKL) divergence relative to the true density $p_t$ whenever the localization parameter $\kappa>0$ is sufficiently small. A key step in the proof is that, to first order in $\kappa$, ISD updates coincide with standard score-driven (SD) updates. This first-order equivalence, however, does not place ISD updates within the framework of the main text, as Theorems~\ref{thm:EKLequivalenceIndepCopy} and \ref{thm:EKLequivalenceBounded} apply to updating rules that are \emph{linearly} downscaled by $\kappa$, whereas the localization in ISD updates is \emph{nonlinear} in $\kappa$. Hence, even for small $\kappa$, a separate argument is required to establish EKL reductions for ISD updates, which we provide here.

For expository clarity, we work in a Euclidean parameter space (thus excluding parameter-space boundaries) and assume a static learning-rate matrix $A \succ \zeromatrix$. A time-varying version can be accommodated under mild measurability conditions. While the result below is stated for sufficiently small values of the localization parameter $\kappa$, for completeness we note that ISD updates can differ substantially from SD updates for non-localized updates (i.e., if $\kappa$ is fixed at unity), as shown in \cite{lange2022robust} and \cite{donker2025stability}.

We use the following assumption.

\begin{aassumption}\label{ass:ISD}
    \begin{enumerate}[nosep]
   \item[\textnormal{(i)}] 
   \emph{(Euclidean parameter space)}  $
   \vartheta\in\Theta=\mathbb{R}^k$.
\item[\textnormal{(ii)}] \emph{(Smoothness)} For all $y\in\mathcal Y$, the map $\vartheta\mapsto \log f(y|\vartheta)$ is twice continuously differentiable with score $\score(y,\vartheta)$ and Hessian $\Hessian(y,\vartheta)$.
\item[\textnormal{(iii)}] \emph{(One-sided Hessian bound)} There exists $c<\infty$ such that $\Hessian(y,\vartheta)\ \preceq\ c\,\unitmatrix$ for all $y\in\mathcal Y$ and all $\vartheta\in\mathbb R^k$.
\item[\textnormal{(iv)}] \emph{(Bounded expected Hessian in operator norm)} There exists $C<\infty$ such that \\
$
\sup_{\vartheta\in \mathbb{R}^k}\Big\|\E_{Y_t\sim p_t}\big[\Hessian(Y_t,\vartheta)\big]\Big\|\ \le\ C,
$
where $\|\cdot\|$ denotes the operator norm.
\item[\textnormal{(v)}] \emph{(Square-integrability of the score at $\vartheta_\tbu$)} $
\E_{Y_t\sim p_t}\big[\ \|\score(Y_t,\vartheta_\tbu)\|^2\ \big]\ <\ \infty.
$
\item[\textnormal{(vi)}] \emph{(Independence and non-zero expected score)} $X_t$ and $Y_t$ are independent and identically distributed with marginal $p_t$, and $
\mu_{p_t}:=\E_{Y_t\sim p_t}\big[\score(Y_t,\vartheta_\tbu)\big]\ \neq\ 0$.
\end{enumerate}
\end{aassumption}

We then obtain the following result. 

\begin{atheorem}[ISD updates are EKL reducing for sufficiently small $\kappa$]
\label{prop:ISD_EKL_smallkappa_noUB_explained}
Let Assumption~\ref{ass:ISD} hold. Fix the prediction $\vartheta_\tbu\in\mathbb R^k$ and learning-rate matrix $A\succ\zeromatrix$. For a localization parameter $0<\kappa\leq 1$ and an observation $y\in\mathcal Y$, define the \emph{localized ISD update} as
\begin{equation}
\label{eq:ISD_prox_def_lowner}
\vartheta^\mathrm{ISD}_\tpu(y)
=\arg\max_{\vartheta\in\mathbb R^k}\Big\{
\kappa \log f(y|\vartheta)
-\tfrac12(\vartheta-\vartheta_\tbu)^\top A^{-1}(\vartheta-\vartheta_\tbu)
\Big\},
\end{equation}
where the usual ISD update takes $\kappa=1$. Then there exists $\bar\kappa>0$ such that for all $\kappa\in(0,\bar\kappa)$:
\begin{enumerate}[nosep]
\item[\textnormal{(a)}] the maximizer in \eqref{eq:ISD_prox_def_lowner} exists and is unique for every $y\in\mathcal Y$;
\item[\textnormal{(b)}] for every fixed $y\in\mathcal Y$,
$\vartheta^\mathrm{ISD}_\tpu(y)
=\vartheta_\tbu+\kappa A\,\score(y,\vartheta_\tbu)+o(\kappa)$;
\item[\textnormal{(c)}] writing $\phi_\mathrm{ISD}(y,\vartheta_{t|t-1}):=\vartheta^\mathrm{ISD}_\tpu(y)$, we have $-\Delta^{\mathsf{EKL}}(\phi_\mathrm{ISD}| \vartheta_\tbu,p_t)=\kappa\,\mu_{p_t}^\top A\mu_{p_t}+o(\kappa)$, where $\mu_{p_t}^\top A \mu_{p_t}>0$; hence, ISD updates are EKL reducing for sufficiently small $\kappa>0$. 
\end{enumerate}
\end{atheorem}

\begin{proof}
\textbf{Proof of part (a).}
Fix $y\in\mathcal Y$ and define the objective
\begin{equation}
Q_\kappa(\vartheta;y)
:=\kappa \log f(y|\vartheta)
-\tfrac12(\vartheta-\vartheta_\tbu)^\top A^{-1}(\vartheta-\vartheta_\tbu),\label{eq: Q def}
\end{equation}
so that  $\vartheta_{t|t}^\mathrm{ISD}(y)$ is the maximizer of this function, i.e., $\vartheta_{t|t}^\mathrm{ISD}(y):=\arg\max_{\vartheta\in\mathbb R^k} Q_\kappa(\vartheta;y)$, where the maximization is over $\Theta=\mathbb{R}^k$ by Assumption~\ref{ass:ISD}(i). For sufficiently small $\kappa>0$, we prove that the maximizer exists and is unique. To show existence, we write $\vartheta=\vartheta_\tbu+\Delta$. By the smoothness of Assumption~\ref{ass:ISD}(ii), the function $\vartheta\mapsto \log f(y|\vartheta)$ is twice continuously differentiable.
Together with the one-sided Hessian bound in Assumption~\ref{ass:ISD}(iii), a second-order Taylor inequality gives, for every $\Delta\in\mathbb R^k$,
\begin{equation}
\label{eq:Taylor_upper_noUB}
\log f(y|\vartheta_\tbu+\Delta)
\le \log f(y|\vartheta_\tbu)+\score(y,\vartheta_\tbu)^\top\Delta+\tfrac{c}{2}\,\|\Delta\|^2,
\end{equation}
because the second-order term can be bounded above by $\frac12\Delta^\top(cI)\Delta=\frac{c}{2}\|\Delta\|^2$ when $\Hessian(y,\cdot)\preceq c \unitmatrix$.
Substituting \eqref{eq:Taylor_upper_noUB} into $Q_\kappa(\vartheta_\tbu+\Delta;y)$ yields
\begin{align}
Q_\kappa(\vartheta_\tbu+\Delta;y)&=\kappa \log f(y|\vartheta_\tbu+\Delta)
-\tfrac12\Delta^\top A^{-1}\Delta\notag
\\
&\le \kappa\log f(y|\vartheta_\tbu)
+\kappa\,\score(y,\vartheta_\tbu)^\top\Delta
-\tfrac12\,\Delta^\top\big(A^{-1}-\kappa c\,\unitmatrix\big)\Delta.
\label{eq: quadratic term}
\end{align}
Next, we choose $\kappa>0$ such that 
\begin{equation}
    A^{-1}-\kappa c\,\unitmatrix\succ\zeromatrix. \label{eq:kappa matrix restriction}
\end{equation}
If $c\le 0$, then \eqref{eq:kappa matrix restriction} automatically holds without any restriction on $\kappa>0$, as $A^{-1}$ is positive definite by assumption. If $c>0$, we can make~\eqref{eq:kappa matrix restriction} hold by making $\kappa$ sufficiently small. Specifically, for $c>0$ we will impose
\begin{equation}
 \label{eq:kappa stricter range} 
0<\kappa\ \leq \ \frac{1}{2}\frac{1}{c\,\lambda_{\max}(A)},
\end{equation} 
which ensures that 
\begin{equation}
\label{eq: kappa lower bound on eigenvalue}
\lambda_{\min}\Big( A^{-1} -\kappa c \unitmatrix\Big) \quad \geq \quad \lambda_{\min}(A^{-1})-c \,\lambda_{\max} (\kappa \unitmatrix) \quad \geq  \quad \frac{1}{2\lambda_{\max}(A)} ,
\end{equation}
where, in the last equality, we used $\lambda_{\min}(A^{-1})=1/\lambda_{\max}(A)$ as well as \eqref{eq:kappa stricter range}. Hence, for $\kappa$ satisfying~\eqref{eq:kappa stricter range}, the matrix $A^{-1} -\kappa c \unitmatrix$ is positive definite, where the smallest eigenvalue is lower bounded by $1/(2\lambda_{\max}(A))$. Importantly, this lower bound is independent of $\kappa$.

For $\kappa$ satisfying~\eqref{eq:kappa stricter range}, the quadratic form $\Delta^\top(A^{-1}-\kappa cI)\Delta$ in~\eqref{eq: quadratic term} is strictly positive for $\Delta\neq 0$,
so the last term on the right-hand side of~\eqref{eq: quadratic term} is strictly negative and dominates the linear term as $\|\Delta\|\to\infty$.
Therefore the right-hand side tends to $-\infty$ as $\|\Delta\|\to\infty$, and hence
$Q_\kappa(\vartheta;y)\to -\infty$ as $\|\vartheta\|\to\infty$.
Since $Q_\kappa(\cdot;y)$ is continuous on $\mathbb R^k$ by Assumption~\ref{ass:ISD}(ii), it attains its maximum; hence, we have shown existence.

To prove uniqueness of the maximizer, nearly the same argument suffices. By Assumption~\ref{ass:ISD}(ii), the Hessian of $Q_\kappa(\cdot;y)$ exists and equals
\begin{equation}
\label{eq: Hessian of Q}
\frac{ \partial^2 Q_\kappa(\vartheta;y)}{\partial \vartheta\partial\vartheta^\top}=\kappa\,\Hessian(y,\vartheta)-A^{-1}\preceq\ \kappa c\,\unitmatrix-A^{-1}
=-(A^{-1}-\kappa c\,\unitmatrix)\; \prec\; \zeromatrix,
\end{equation}
where the first inequality follows from the upper-bounded Hessian in Assumption~\ref{ass:ISD}(iii), and the last inequality follows from~\eqref{eq:kappa matrix restriction}. This implies that the Hessian of $Q_\kappa(\cdot;y)$ can be made strictly negative definite: either for all $\kappa>0$ if $c\leq 0$, or for sufficiently small $\kappa>0$ if $c>0$, in which case $\kappa$ satisfying~\eqref{eq:kappa stricter range} is sufficient. It follows that for sufficiently small $\kappa>0$, the Hessian of $Q_\kappa(\cdot;y)$ is negative definite for all $\vartheta\in \mathbb{R}^k$; hence, for every $y\in \Ym$, the map $\vartheta\mapsto Q_\kappa(\vartheta;y)$ is strictly concave. A strictly concave function on $\mathbb{R}^k$ has at most one maximizer; i.e., the maximizer is unique. This proves part~(a).

\medskip
\textbf{Proof of part (b).}
Fix $y\in\mathcal Y$, let $\vartheta^\mathrm{ISD}_\tpu(y)$ be the unique maximizer from part (a), and denote the \emph{difference} between the prediction and the update, indexed by $\kappa>0$, as
\[
\Delta_\kappa(y):=\vartheta^\mathrm{ISD}_\tpu(y)-\vartheta_\tbu.
\]
By part (a) and the optimality of $\vartheta^\mathrm{ISD}_\tpu(y)$, the first-order condition of the objective $Q_\kappa(\vartheta;y)$ evaluated at the maximizer $\vartheta^\mathrm{ISD}_\tpu(y)$ is
\begin{equation}
\label{eq:FOC_ISD_noUB_explained}
0=\kappa\,\score(y,\vartheta^\mathrm{ISD}_\tpu(y))-A^{-1}\Delta_\kappa(y),
\quad\text{equivalently}\quad
\Delta_\kappa(y)=\kappa A\,\score(y,\vartheta^\mathrm{ISD}_\tpu(y)).
\end{equation}
We first show that $\Delta_\kappa(y)\to 0$ as $\kappa\downarrow 0$.
To this end, note that by the upper-bounded Hessian in Assumption~\ref{ass:ISD}(iii), we have~\eqref{eq: Hessian of Q}. To analyze this inequality further, define the positive constant
\begin{equation}
\label{eq:m constant}
m_\kappa:=\lambda_{\min}(A^{-1}-\kappa c\,\unitmatrix)>0,
\end{equation}
where positivity follows automatically if $c\leq 0$, while if $c>0$ it follows by taking $\kappa$ to satisfy~\eqref{eq:kappa stricter range}. If we impose~\eqref{eq:kappa stricter range}, then it follows that
\begin{equation}
    m_\kappa\;=\;\lambda_{\min}(A^{-1})-\kappa c\;\geq\; 
\frac{1}{2\lambda_{\max}(A)}>0.\label{eq: m_kappa minimum}
\end{equation}
Summing up, \eqref{eq: Hessian of Q} implies
\begin{equation}
\label{eq: give condition}
\frac{ \partial^2 Q_\kappa(\vartheta;y)}{\partial \vartheta\partial\vartheta^\top}
 \preceq -m_\kappa \unitmatrix, \quad \text{ where }\quad m_\kappa\geq \frac{1}{2\lambda_{\max}(A)}>0,
\end{equation}
for all $\kappa>0$ satisfying~\eqref{eq:kappa stricter range}. For such $\kappa>0$, the function $\vartheta\mapsto Q_\kappa(\vartheta;y)$ is $m_\kappa$-\emph{strongly concave} with parameter of strong concavity given by $m_\kappa>(2\lambda_{\max}(A))^{-1}$. 

Strong concavity (by the inequalities in~\eqref{eq: give condition}) implies that the gradient difference of $Q_\kappa(\vartheta;y)$, evaluated at two points $\vartheta_1$ and $\vartheta_2$, satisfies the following standard inequality:
\begin{equation}
\label{eq:strong_concavity_gradient}
\left(\frac{\partial Q_\kappa(\vartheta_1;y)}{\partial\vartheta_1}-\frac{\partial Q_\kappa(\vartheta_2;y)}{\partial\vartheta_2}\right)^\top(\vartheta_1-\vartheta_2)
\le -m_\kappa\,\|\vartheta_1-\vartheta_2\|^2,\quad \forall \vartheta_1,\vartheta_2\in \Theta,y\in \mathcal{Y}.
\end{equation}
Next, we take \eqref{eq:strong_concavity_gradient} with $\vartheta_1=\vartheta^\mathrm{ISD}_\tpu(y)$ and $\vartheta_2=\vartheta_\tbu$, such that $\vartheta_1-\vartheta_2=\Delta_{\kappa}(y)$.
On the left-hand side of~\eqref{eq:strong_concavity_gradient}, we can use
\begin{equation}
   \left.\frac{\partial Q_\kappa(\vartheta_1;y)}{\partial\vartheta_1}\right|_{\vartheta_1=\vartheta^\mathrm{ISD}_\tpu(y)}=0,\qquad   \left.\frac{\partial Q_\kappa(\vartheta_2;y)}{\partial\vartheta_2}\right|_{\vartheta_2=\vartheta_\tbu}=\kappa\,\score(y,\vartheta_\tbu),
\end{equation}
where the equation on the left is simply the first-order condition for the optimality of $\vartheta^\mathrm{ISD}_{t|t}(y)$, while the equation on the right follows from the definition of $Q_\kappa$ in~\eqref{eq: Q def}. For these specific choices of $\vartheta_1$ and $\vartheta_2$, inequality \eqref{eq:strong_concavity_gradient} becomes
\[
\big(0-\kappa\,\score(y,\vartheta_\tbu)\big)^\top\Delta_\kappa(y)
\le -m_\kappa\,\|\Delta_\kappa(y)\|^2.
\]
Multiplying by negative one, we have
\[
m_\kappa\,\|\Delta_\kappa(y)\|^2
\;\le\; \kappa\,\score(y,\vartheta_\tbu)^\top\Delta_\kappa(y)
\;\le\; \kappa\,\|\score(y,\vartheta_\tbu)\|\,\|\Delta_\kappa(y)\|,
\]
where the last inequality follows by applying Cauchy--Schwarz. The inequality is trivial if $\|\Delta_\kappa(y)\|= 0$; hence, we may assume $\|\Delta_\kappa(y)\|\neq 0$ and divide both sides by $\|\Delta_\kappa(y)\|$ to obtain
\begin{equation}
\label{eq:Delta_bound_by_score_noUB_explained}
\frac{\|\Delta_\kappa(y)\|}{\kappa}\ \le\ \frac{1}{m_\kappa}\,\|\score(y,\vartheta_\tbu)\|.
\end{equation}
For fixed $y$, the right-hand side is finite, and from~\eqref{eq:m constant} we see that $m_\kappa\geq 1/(2\lambda_{\max}(A))>0$ as $\kappa\downarrow 0$,
so \eqref{eq:Delta_bound_by_score_noUB_explained} implies $\|\Delta_\kappa(y)\|\to 0$ as $\kappa \downarrow 0$. Indeed, $\|\Delta_\kappa(y)\|=\mathcal{O}(\kappa)$. 

Next, we use this fact inside the first-order condition \eqref{eq:FOC_ISD_noUB_explained}. Since $\vartheta^\mathrm{ISD}_\tpu(y)=\vartheta_\tbu+\Delta_\kappa(y)\to\vartheta_\tbu$
and because the score is continuous in $\vartheta$, the smoothness of Assumption~\ref{ass:ISD}(ii) implies that
$\score(y,\vartheta^\mathrm{ISD}_\tpu(y))\to \score(y,\vartheta_\tbu)$ as $\kappa\downarrow 0$.
By \eqref{eq:FOC_ISD_noUB_explained}, this means
\[
\frac{\Delta_\kappa(y)}{\kappa}
=A\,\score(y,\vartheta^\mathrm{ISD}_\tpu(y))
\to A\,\score(y,\vartheta_\tbu),
\qquad(\kappa\downarrow 0),
\]
which is equivalent to
\[
\Delta_\kappa(y)=\kappa A\,\score(y,\vartheta_\tbu)+o(\kappa),\qquad(\kappa\downarrow 0),
\]
which proves part (b).
\medskip

\textbf{Proof of part (c).}
Fix $x,y\in\mathcal Y$ and recall $\Delta_\kappa(y):=\vartheta^\mathrm{ISD}_{t|t}(y)-\vartheta_{t|t-1}$. By the smoothness of Assumption~\ref{ass:ISD}(ii), an exact integral version of the second-order Taylor expansion gives
\begin{align}
\log f(x|\vartheta_\tpu(y))-\log f(x|\vartheta_\tbu)&=\log f(x|\vartheta_\tbu+\Delta_\kappa(y))-\log f(x|\vartheta_\tbu)\notag
\\
&=\score(x,\vartheta_\tbu)^\top \Delta_\kappa(y) + R_\kappa(x,y),
\label{eq:path_integral_loglik_noUB_explained}
\end{align}
where, as is standard, the remainder term $R_\kappa(x,y)$ contains an integral over the Hessian as follows:
\begin{equation}
R_\kappa(x,y)\;: =\;\Delta_\kappa(y)^\top \left( \int_0^1 (1-u) \,\Hessian(x,\vartheta_\tbu+u\Delta_\kappa(y))\,\mathrm du\right) \,\Delta_\kappa(y),
\label{eq: remainder}
\end{equation}
Here, the factor $(1-u)$ under the integral is obtained by a standard trick involving changing the order of the double integral in the usual second-order expansion (e.g., \citealp[p.~769, Prop.~A23]{bertsekas2016nonlinear}). 

Next, we evaluate~\eqref{eq:path_integral_loglik_noUB_explained} at the random variables $X_t$ and $Y_t$ instead of $x$ and $y$ and take a (joint) expectation over the independent draws $X_t,Y_t\sim p_t$ to obtain the negative EKL difference on the left-hand side, i.e.,
\begin{align}
-\Delta^{\mathsf{EKL}}(\phi_\mathrm{ISD}\mid \vartheta_\tbu,p_t) &= \E_{X_t,Y_t\sim p_t}\big[\log f(X_t|\vartheta_\tpu(Y_t))-\log f(X_t|\vartheta_\tbu)\big]\notag
\\
&=\E_{X_t,Y_t\sim p_t}\big[\score(X_t,\vartheta_\tbu)^\top \Delta_\kappa(Y_t)\big]+\E_{X_t,Y_t\sim p_t}\big[R_\kappa(X_t,Y_t)\big],
\label{eq:EKL_decomp_noUB_explained}
\end{align}
which we hope to show is positive, at least to leading order in $\kappa$. By independence of $X_t$ and $Y_t$ as ensured by Assumption~\ref{ass:ISD}(vi), the expectation in the first term of~\eqref{eq:EKL_decomp_noUB_explained} can be split up as follows:
\[
\E_{X_t,Y_t\sim p_t}\big[\score(X_t,\vartheta_\tbu)^\top \Delta_\kappa(Y_t)\big]
=\E_{X_t\sim p_t}[\score(X_t,\vartheta_\tbu)]^\top\,\E_{Y_t\sim p_t}[\Delta_\kappa(Y_t)]
=\mu_{p_t}^\top\,\E_{Y_t\sim p_t}[\Delta_\kappa(Y_t)].
\]
Next, we claim that $\E_{Y_t\sim p_t}[\Delta_\kappa(Y_t)]/\kappa\to A\mu_{p_t}$.
Indeed, part (b) gives pointwise convergence $\Delta_\kappa(y)/\kappa\to A\,\score(y,\vartheta_\tbu)$ for each fixed $y$. Here, however, we require an extra argument to be able to exchange the limit and the expectation. In fact, this is permitted due to~\eqref{eq:Delta_bound_by_score_noUB_explained}, which combined with~\eqref{eq: m_kappa minimum} gives
\[
\Big\|\frac{\Delta_\kappa(Y_t)}{\kappa}\Big\|
\le \frac{1}{m_\kappa}\,\|\score(Y_t,\vartheta_\tbu)\| \le 2\lambda_{\max}(A)\,\|\score(Y_t,\vartheta_\tbu)\|,
\]
where the validity of~\eqref{eq: m_kappa minimum} is ensured for all $\kappa>0$ satisfying~\eqref{eq:kappa stricter range}. Now, integrability of the right-hand side is ensured by Assumption~\ref{ass:ISD}(v), such that its expectation exists. This means that we can apply a dominated-convergence argument as $\kappa\to 0$ to obtain
\[
\frac{\E_{Y_t\sim p_t}[\Delta_\kappa(Y_t)]}{\kappa}\ \to\ A\,\E_{Y_t\sim p_t}[\score(Y_t,\vartheta_\tbu)]=A\mu_{p_t},
\] 
as claimed. Therefore, the leading-order term (in $\kappa$) of~\eqref{eq:EKL_decomp_noUB_explained} becomes
\begin{equation}
\label{eq:main_term_asymp_noUB_explained}
\E_{X_t,Y_t\sim p_t}\big[\score(X_t,\vartheta_\tbu)^\top \Delta_\kappa(Y_t)\big]
=\kappa\,\mu_{p_t}^\top A\mu_{p_t}+o(\kappa),\qquad(\kappa\downarrow 0),
\end{equation}
as claimed. To complete the proof, we still need to use Assumption~\ref{ass:ISD}(iv) to show that the remainder term $\E_{X_t,Y_t\sim p_t}[R_\kappa(X_t,Y_t)]$ in~\eqref{eq:EKL_decomp_noUB_explained} is $o(\kappa)$. To this end, we consider
\begin{align*}
& \Big| \E_{X_t,Y_t\sim p_t} \big[R_\kappa(X_t,Y_t)\big] \Big|
\\
&= \left|  \E_{X_t,Y_t\sim p_t}\left[ \Delta_\kappa(Y_t)^\top  \left( \int_0^1 (1-u)\,\Hessian(X_t,\vartheta_\tbu+u\Delta_\kappa(Y_t))\,\mathrm du\right)\,\Delta_\kappa(Y_t)\, \right] \right|
 \text{\footnotesize by~\eqref{eq: remainder}}\\
 &=  \left| \E_{Y_t\sim p_t}\left[ \Delta_\kappa(Y_t)^\top  \left( \int_0^1 (1-u)\,\E_{X_t\sim p_t}\Big[\Hessian(X_t,\vartheta_\tbu+u\Delta_\kappa(Y_t))\Big]\,\mathrm du\right)\,\Delta_\kappa(Y_t)\, \right] \right|
 \\
 & 
\qquad \text{\footnotesize by independence as in Assumption~\ref{ass:ISD}(vi)}
 \\
   &\leq   \E_{Y_t\sim p_t}\left[  \left\|  \int_0^1 (1-u)\,\E_{X_t\sim p_t}\Big[\Hessian(X_t,\vartheta_\tbu+u\Delta_\kappa(Y_t))\Big]\,\mathrm du \right\| \;\| \Delta_\kappa(Y_t)\|^2\, \right] 
   \\
      &\leq   \E_{Y_t\sim p_t}\left[   \int_0^1 (1-u)\, \left\|  \E_{X_t\sim p_t}\Big[\Hessian(X_t,\vartheta_\tbu+u\Delta_\kappa(Y_t))\Big] \right\|\mathrm du \;\| \Delta_\kappa(Y_t)\|^2\, \right] 
   \\ 
  &\leq    \E_{Y_t\sim p_t}\left[\left\{ \sup_{\vartheta\in \Theta} \Big \|\E_{X_t\sim p_t}[\Hessian(X_t,\vartheta)]\Big \| \right\}  \;  \int_0^1 (1-u)\,\mathrm du\;  \| \Delta_\kappa(Y_t)\|^2\, \right]
   \\
   &\leq \frac{C}{2}\,\E_{Y_t\sim p_t}[\|\Delta_\kappa(Y_t)\|^2]\quad \text{\footnotesize by Assumption~\ref{ass:ISD}(iv)}
\\
& \leq  \frac{C}{2} \frac{\kappa^2}{m_\kappa^2}\,\E_{Y_t\sim p_t}[ \|\score(Y_t,\vartheta_\tbu)\|^2]
 \quad \text{\footnotesize by~\eqref{eq:Delta_bound_by_score_noUB_explained}}
 \\
 &\leq  2 C\, \lambda_{\max}(A)^2  \, \kappa^2\,\E_{Y_t\sim p_t}[ \|\score(Y_t,\vartheta_\tbu)\|^2]  \quad \text{\footnotesize by~\eqref{eq: m_kappa minimum}}
\\
&=  \mathcal{O}(\kappa^2)\qquad \text{\footnotesize by Assumption~\ref{ass:ISD}(v)} 
\end{align*}
 In the last line, the square integrability of Assumption~\ref{ass:ISD}(v) allows us to conclude that $| \E_{X_t,Y_t\sim p_t} [R_\kappa(X_t,Y_t)]|=\mathcal{O}(\kappa^2)$, which is obviously $o(\kappa)$. Using this in  \eqref{eq:EKL_decomp_noUB_explained} and applying~\eqref{eq:main_term_asymp_noUB_explained}, we conclude
\[
-\Delta^{\mathsf{EKL}}(\phi_\mathrm{ISD}\mid \vartheta_\tbu,p_t)
=\kappa\,\mu_{p_t}^\top A\mu_{p_t}+o(\kappa),\qquad(\kappa\downarrow 0).
\]
Since $A\succ\zeromatrix$ and $\mu_{p_t}\neq 0$ by Assumption~\ref{ass:ISD}(vi), we have $\mu_{p_t}^\top A\mu_{p_t}>0$; hence, the expression is positive for all sufficiently small $\kappa$. This concludes the proof of part (c).
\end{proof}

\section{Quasi score-driven updates}
\label{app:qsd}

\cite{blasques2023quasi} generalize the class of score-driven (SD) updates to the class of so-called Quasi-SD (QSD) updates by allowing the update to be based on the score of a postulated density $\tilde{f}(\cdot \vert \vartheta_\tbu)$ that possibly differs from the model density ${f}(\cdot \vert \vartheta_\tbu)$. 
While both, $f$ and $\tilde f$ are assumed to be driven by the same time-varying parameter $\vartheta_\tbu$, neither of these distributions is assumed to coincide with the truth $p_t \in \Pm$ for some class $\Pm$.
If $\tilde{f} \equiv f$, the class of SD updates is obtained, which is  EKL reducing by Corollary~\ref{cor:EKL_SDupdate}. 
Formally, QSD updates are defined as 
\begin{align*}
	\phi_{\mathrm{QSD}}(y_t, \vartheta_\tbu) := \vartheta_\tbu + A \scaling  \widetilde{\score}(y_t, \vartheta_\tbu), \qquad 
    \widetilde{\score}(y_t,\vartheta_\tbu) := \left.\frac{\partial \log \tilde{f}(y_t \vert \vartheta)}{\partial \vartheta}\right\vert_{\vartheta_\tbu}.
\end{align*}
As we have characterized all updating rules that are guaranteed to be  EKL reducing in Theorem~\ref{thm:EKLequivalenceIndepCopy}, we can now exploit this result to immediately conclude that  $\phi_{\mathrm{QSD}}$ is EKL reducing w.r.t.~$\Pm$ if, and only if,
\begin{align} 
	\label{eq:EKLQSDCondition}
	\E_{p_t} \big[ \widetilde{\score}(Y_t, \vartheta_\tbu) \big]^\top \, (A \scaling) \, \E_{p_t} \big[ \score(X_t,\vartheta_\tbu) \big] >0,
\end{align}
for all $\vartheta_\tbu \in \Theta$ and $p_t \in \Pm$ such that 
$\E_{p_t}[\widetilde{\score}(Y_t,\vartheta_\tbu)]^\top\, (A \scaling)\, \E_{p_t}[\score(X_t,\vartheta_\tbu)] \;\not=\; 0$.
Hence, QSD updates need not be EKL reducing; this depends, among other factors, on the similarity of $\score(Y_t,\vartheta_\tbu)$ and $\widetilde{\score}(X_t,\vartheta_\tbu)$ and on the richness of the class $\Pm$ such that  \eqref{eq:EKLQSDCondition} can be guaranteed to hold for all $p_t \in \Pm$.

\begin{aexample}[GARCH-$t$] 
	\cite{blasques2023quasi} show that the univariate GARCH-$t$ filter of \citet{bollerslev1987conditionally} is an example of a QSD update where $\tilde f(\cdot \vert \vartheta_{\tbu})$ is a Student's $t$ density with $\nu \in \mathbb{R}$, $\nu > 2$ degrees of freedom and time-varying variance $\vartheta_\tbu$, and $f(\cdot \vert \vartheta_{\tbu})$ is Gaussian density with time-varying variance $\vartheta_\tbu$. 
	To formally obtain equivalence with the \hbox{GARCH-$t$} filter, the scaling factor $\scaling = -\left(\E_{\tilde f_{t\vert t-1}}\left[\tilde{H}(Y_t, \vartheta_{t\vert t-1}) \right]\right)^{-1}$ is required, where $\tilde{H}(y_t, \vartheta_{t\vert t-1}):=\frac{\partial^2\tilde{f}(y_t \vert \vartheta) }{\partial \vartheta^2}\big\vert_{\vartheta_\tbu}$. 
    
    In the univariate case, $A = \alpha > 0$ and as $\scaling$ is $\mathcal{F}_{t-1}$-measurable and positive, $A \scaling$ is positive and we can ignore them in \eqref{eq:EKLQSDCondition}.
	Then, straight-forward calculations show that \eqref{eq:EKLQSDCondition} is satisfied if and only if
	\begin{align*}
		\E_{p_t} \left[Y_t^2-\vartheta_\tbu \right] \E_{p_t} \left[   \frac{\nu+1}{\nu-2 +X_t^2/\vartheta_\tbu}X_t^2-\vartheta_\tbu\right] >0.
	\end{align*}
	While this holds trivially for the SD case of $\nu \to \infty$, it may or may not hold for any fixed $\nu > 2$, depending on the class of true distributions $p_t \in \Pm$.
    \color{black}
\end{aexample}

\setcounter{equation}{0}
\section{CEV-reduction guarantees for scaled Hessians}
\label{sec:ScaledHessianCEVGeneral}

Equation~\eqref{eqn:CEVcondition2} in Section~\ref{sec:CEV} of the main paper illustrates that negative definite expected Hessian matrices---as formalized through Assumption~\ref{ass:HN}~(i)---are required to guarantee CEV (and MSE) reductions.
The closely related Assumption~3 of \citet{Gorgi2023} however concerns a \emph{scaled} Hessian, which differs slightly from the condition required for the derivation in \eqref{eqn:CEVcondition2}.
We hence provide an alternative derivation based on scaled Hessians that aligns more closely with \citet{Gorgi2023}.

To this end, we consider scaling matrices that depend on the time-varying parameter, $\scaling = \scaling(\vartheta)$, for which we obtain
\begin{align}
    &\hspace{-0.5cm}\E_{p_t} \left[ \Delta \phi_{\mathrm{SSD}}(Y_t,\vartheta_\tbu)\right]^\top \Omega_{t-1} \big( \vartheta_t^\ast-\vartheta_\tbu  \big) \nonumber \\
    &= \E_{p_t} \left[ \score(Y_t,\vartheta_\tbu)^\top \scaling(\vartheta_\tbu)^\top - \score(Y_t,\vartheta^\ast_t)^\top \scaling(\vartheta^\ast_t)^\top \right] A^\top \Omega_{t-1} \big( \vartheta_t^\ast-\vartheta_\tbu  \big) \nonumber \\
    &=\big( \vartheta_t^\ast-\vartheta_\tbu  \big)^\top   \E_{p_t} \left[ - \int_0^1 \Hessian_S \big( Y_t, {\vartheta}_t^u \big) \mathrm{d}u \right] A^\top \Omega_{t-1} \big( \vartheta_t^\ast-\vartheta_\tbu  \big),
    \label{eqn:CEVconditionScaled}    
\end{align}
where $\Hessian_S(y, \vartheta) = (\partial / \partial \vartheta) \big\{ \score(y,\vartheta)^\top \scaling(\vartheta)^\top \big\}$ is the \emph{scaled Hessian}, and ${\vartheta}_t^u := u \vartheta_\tbu + (1-u)\vartheta_t^\ast$ for $u\in[0,1]$. 
If $A^\top \Omega_{t-1} = d_{t-1}$ for some scalar, positive and $\mathcal{F}_{t-1}$-measurable $d_{t-1}$, the inner product in \eqref{eqn:CEVconditionScaled} guarantees CEV reductions for SD models given negative definiteness of $\E_{p_t} \big[ \Hessian_S(Y_t, \vartheta) \big]$ for all $\vartheta \in \Theta$, $p_t \in \Pm$, which is closely related to Assumption~3 of \citet{Gorgi2023}.
This generalization however still requires a specific multivariate setting ($A^\top \Omega_{t-1} = d_{t-1} \unitmatrix$) as well as negative definiteness of the (scaled) Hessian, which is stronger than merely (local) boundedness in the EKL-specific Assumptions~\ref{ass:HB} and \ref{ass:HLB}.

\setcounter{equation}{0}	
\section{Censored KL reductions for SD models}
\label{sec:CKL}
    We start by illustrating the puzzling conclusion from \eqref{eqn:TKL_Problem_Intro} with the following continuation of Example~\ref{exmpl:LKLProblem1}.
    For simplicity, we restrict attention to the univariate case in both the observation and outcome space, that is, $l=k=1$.
    For notational convenience, we further set $\mathcal{S}_{t-1}=1$.

	\begin{figure}[tb]
		\centering
		\begin{subfigure}[b]{0.49\textwidth}
			\centering
			\includegraphics[width=\linewidth]{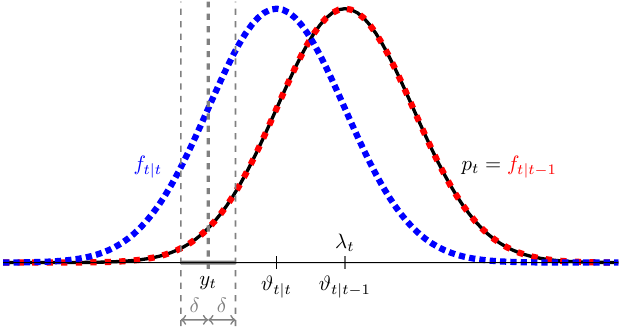}
			\caption{$\lambda_{t} = \vartheta_\tbu >  \vartheta_\tpu$}
			\label{subfig:p=ft}
		\end{subfigure}
		\hfill 
		\begin{subfigure}[b]{0.49\textwidth}
			\centering
			\includegraphics[width=\linewidth]{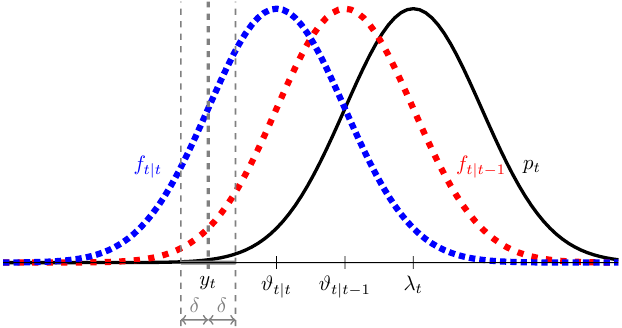} 
			\caption{$\lambda_{t} >  \vartheta_\tbu  > \vartheta_\tpu$}
			\label{subfig:p>ft}
		\end{subfigure}
		\\[1cm]
		\begin{subfigure}[b]{0.49\textwidth}
			\centering
			\includegraphics[width=\linewidth]{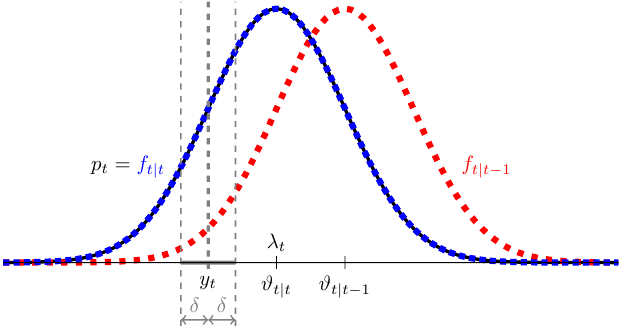}
			\caption{$\lambda_{t} =  \vartheta_\tpu \color{black} <  \vartheta_\tbu$}
			\label{subfig:p=ft+1}
		\end{subfigure}
		\hfill 
		\begin{subfigure}[b]{0.49\textwidth}
			\centering
			\includegraphics[width=\linewidth]{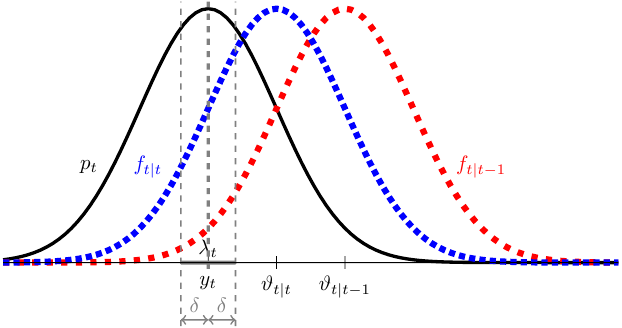} 
			\caption{$\lambda_{t} <  \vartheta_\tpu \color{black} <  \vartheta_\tbu$}
			\label{subfig:p<ft+1}
		\end{subfigure}
		\caption{Illustration of a score-driven conditional mean model with four hypothetical truths in the four panels.} 
		\label{fig:wrongdirectionmean}
	\end{figure}
	
	\begin{aexample}
		\label{exmpl:LKLProblem2}
		Following the setup in Example~\ref{exmpl:LKLProblem1}, consider the true distribution $Y_{t}  \sim p_t=\Nm(\lambda_t,1)$ with unknown time-varying mean $\lambda_t$. 
		We use the (correctly specified) model $f_\tbu = \Nm(\vartheta_\tbu,1)$ based on the time-varying model parameter $\vartheta_\tbu$.
		The Gaussian model likelihood implies the score-driven update $\Delta \phi(y_t, \vartheta_\tbu) = \alpha (y_t - \vartheta_\tbu)$ with learning rate $\alpha > 0$.
		Hence, the score-driven update drives the conditional mean up or down depending on the prediction error $y_t-\vartheta_\tbu$. 
		
		Each panel in Figure~\ref{fig:wrongdirectionmean} shows the starting density $f_\tbu(\cdot) \equiv f(\cdot\vert \vartheta_\tbu)$ in red together with the updated density $f_\tpu(\cdot)  \equiv f(\cdot\vert \vartheta_\tpu)$ in blue based on the (same) realization $y_t$ and associated ball $\ymdel \equiv \ymdellong$ around $y_t$.
		The true (and unknown) density $p_t(\cdot)  \equiv p(\cdot|\lambda_t)$ in black varies across the panels to illustrate four possible scenarios: to the left or right of the predicted density, or coinciding with the predicted or updated density.

		In each panel, the score-driven update guarantees $f_\tpu > f_\tbu$ on $\ymdel$, such that \eqref{eqn:TKL_Problem_Intro} implies a direct improvement in the trimmed KL measure. Therefore, the update should always be beneficial \citep{blasques2015information}.
		While the updated density is indeed shifted towards the true density in panels \subref{subfig:p=ft+1} and \subref{subfig:p<ft+1}, it has moved away from the true density 
		in panels \subref{subfig:p=ft} and \subref{subfig:p>ft}.
		Surprisingly, in panel \subref{subfig:p=ft}, adjustment would be deemed beneficial even though the predicted density was perfect at $f_\tbu = p_t$. 
		(Note that the localization in the state space, which is related to the (small) learning rate $\alpha>0$ in \citet{blasques2015information}, does not affect the conclusion; $\vartheta_\tpu$ being arbitrarily close to $\vartheta_\tbu$ leads to the same problems.) 
	\end{aexample}

    By relying on the closely related topic of localizing scoring rules, the perhaps puzzling conclusion from \eqref{eqn:TKL_Problem_Intro} can be explained by an incorrect localization of the underlying logarithmic scoring rule in the TKL measure.
    In spirit similar to \eqref{eq:LKLDefinition}, \cite{amisano2007comparing} use \emph{trimming} to introduce the weighted likelihood scoring rule, $\mathrm{wl}_R(f,y):=\1_R(y)\log f(y)$, where $R \subseteq \Ym$ is some region of interest and $\1_R(y)$ equals unity if $y\in R$ and zero otherwise. 
	This scoring rule is localizing in that densities that coincide on $R$ achieve the same score for all $y\in \Ym $. 
    
	As pointed out by \cite{diks2011likelihood} and \cite{gneiting_comparing_2011}, however, it fails to be proper, as the conclusion below follows directly from its definition:
	\begin{align}
		\label{eqn:wlProblem}
		f_1(y) > f_2(y), \quad  y\in R 
		\qquad \Longrightarrow \qquad 
		\mathrm{wl}_R(f_1,y) > \mathrm{wl}_R(f_2,y), \quad y\in R.
	\end{align}
	As \eqref{eqn:wlProblem} holds for all $y\in R$, it also holds in expectation, and even if $f_2$ is the true density.
	This deficiency echoes that of the TKL measure in \cite{blasques2015information} mentioned in \eqref{eqn:TKL_Problem_Intro} as we can write 
    $\rkl_{B}(p_{t} \Vert f_\tpu) = \E_{p_{t}} \left[ \wlB \big(p_{t}, X_{t} \big)  - \wlB \big(f_\tpu, X_{t} \big) \right]$.
	For the same reason, the TKL measure can become negative such that it is not a proper divergence measure.

	As an alternative, \cite{diks2011likelihood} propose the proper \emph{censored} likelihood scoring rule $\mathrm{csl}_R(f,y):=\1_R(y)\log f(y) + \1_{R^c}(y)\log F(R^c)$, where $F(R^c)$ denotes the probability the distribution $F$ assigns to $R^c$, the complement of $R$.
	Here, censoring accounts for the omitted set $R^c$ in an aggregated form.
	\cite{depunder2023localzing} recently generalize the concept of censored scoring rules, and for any two densities $p,f \in \Pm$ with $Y \sim p$ and a \emph{localization set} $R \subseteq \Ym$, they define a \emph{local divergence} $\D_R(p\Vert f)$ to satisfy:
    \begin{align}
        \label{eqn:LocDivergence}
		\mathrm{(i)} \;\; \D_R(p\Vert f)\geq 0 
		\quad \text{and} \quad 
		\mathrm{(ii)} \;\;  \D_R(p\Vert f) = 0 \; \iff \; p(Y)\1_{R}(Y) = f(Y)\1_{R}(Y), 
	\end{align}  
	The TKL measure of \cite{blasques2015information} in \eqref{eq:LKLDefinition} does not satisfy the requirements in \eqref{eqn:LocDivergence} and is as such not a proper divergence measure.

This indicates that localizing by \emph{censoring} instead of \emph{trimming} might allow for localized performance guarantees of SD models.
	Hence, we consider the \emph{censored} KL (CKL) divergence
	\begin{align}
		\begin{aligned}
			\label{eq:CKLDefinition}
			\ckl_{B}(p_{t} \Vert f_\tpu) 
			&:= \int_{\ymdel} \log\left(\frac{p_t(x)}{{f}(x|\vartheta_\tpu)}\right) p_t(x) \dd x  
			+ \log \left(\frac{ \int_{\ymdelc} p_t(x) \dd  x}{\int_{\ymdelc}  {f}(x|\vartheta_\tpu) \dd x} \right) \int_{B^c} p_t(x) \dd x,
		\end{aligned}
	\end{align}
	where, as before,  $\ymdel \equiv \ymdellong =\{x \in \Ym: \vert x-y_t \vert \leq \delta\}$ is a ball with a small radius $\delta>0$ around the realization $y_t$. 
    The first part in~\eqref{eq:CKLDefinition} is equivalent to the TKL measure proposed by \cite{blasques2015information} in \eqref{eq:LKLDefinition}, whereas the second part adds a correction term for the ignored (but, in its aggregated form, necessary) information on $B^c$. \citet{depunder2023localzing} show that this CKL divergence, which is based on the censored likelihood score in \citet{diks2011likelihood}, is a \emph{local divergence measure} in the sense of \eqref{eqn:LocDivergence}, for which the inclusion of the second term in \eqref{eq:CKLDefinition} is crucial. 

	On an intuitive level, localizing by \emph{trimming} on $B$ disregards the entire behaviour on $B^c$, with the abovementioned adverse consequences. 
	In contrast, localizing by \emph{censoring} disregards the behaviour on $B^c$ as far as possible while still accounting for aggregate behaviour. By doing so, the CKL divergence inherits the attractive theoretical properties of the (global) KL divergence.

    We now analyze whether censoring can provide sensible outcome-space localized performance guarantees for SD models that hold almost surely, instead of in expectation as our EKL measure.

	\begin{adefinition}[CKL difference] 
		\label{def:CKL}
		For any $(y_t,\vartheta_\tbu,p_t) \in \Ym \times \Theta \times \Pm$ and $\delta>0$, the \emph{CKL difference} for an updating rule $\phi$ is
		\begin{align*}
			\Delta_{\delta}^{\mathsf{CKL}}(\phi)
			\equiv \Delta_{\delta}^{\mathsf{CKL}}(\phi |y_t,\vartheta_\tbu,p_{t}) 
			:&=  \mathsf{CKL}_\ymdellong(p_{t} \Vert f_\tpu) - \mathsf{CKL}_\ymdellong(p_{t} \Vert f_\tbu). 
		\end{align*}
    We say that an update $\phi$ is CKL reducing w.r.t.\ the class $\Pm$ if for all $(y_t,\vartheta_\tbu,p_t) \in \Ym \times \Theta \times \Pm$ such that $\big(p_t(y_t)-f(y_t|\vartheta_{\tbu})\big)  \score(y_t, \vartheta_{\tbu}) \Delta \phi(y_t,\vartheta_{\tbu})\neq 0$ there exists $\delta>0$ that guarantees $\Delta_{\delta}^{\mathsf{CKL}}<0$. 
	\end{adefinition}

Similar to the excluded orthogonality condition $\E_{p_t}\!\left[ \Delta\phi (Y_t,\vartheta_\tbu) \right]^\top 
	\E_{p_t}\!\left[ \score(X_t,\vartheta_\tbu) \right] \not= 0$ in Definition~\ref{def:EKL}, the definition of CKL reducing updates  requires that $$\big(p_t(y_t)-f(y_t\vert\vartheta_{t\vert t-1})\big) s(y_t,\vartheta_{t\vert t-1})\Delta \phi(y_t,\vartheta_{\tbu}) \neq 0,$$
reducing to $\big(p_t(y_t)-f(y_t\vert\vartheta_{t\vert t-1})\big) s(y_t,\vartheta_{t\vert t-1}) \neq 0$ for score-driven updates.    
If
$p_t(y_t)=f(y_t\vert\vartheta_{t\vert t-1})$, the postulated model matches the truth locally at the realized observation, so no first order CKL improvement is available. If $
s(y_t,\vartheta_{t\vert t-1})=0$, the likelihood is locally flat at $y_t$, implying that infinitesimal parameter changes do not affect the local fit. In both cases, CKL changes are governed by higher order terms, whose signs are unknown.
   
\begin{aassumption} 
	\label{ass:thmifflocal}
    \begin{enumerate}[label=(\roman*), itemsep=0.2em, topsep=0.5em]
        \item 
        $\Ym$ and $\Theta$ are open and convex, $\vartheta_\tbu \in \Theta$ and $\vartheta_\tpu(y_t) \in \Theta$.

        \item 
        $\log f(x|\vartheta)$ is twice continuously differentiable in $\vartheta$ and Lipschitz continuous in $x$,  $\forall \vartheta \in \Theta$ and  $x \in \Ym$.

        \item
        $p_t(x)$ is Lipschitz continuous in $x$ and bounded $\forall x \in \Ym$, $\forall p_t\in \Pm$.
    \end{enumerate}
\end{aassumption}
\smallskip

The next theorem indicates the condition under which SD updates imply an improvement in the CKL divergence: $\Delta_{\delta}^{\mathsf{CKL}}(\phi) < 0$. 
For clarity, we focus on the SD case here and generalize this result to general (i.e., possibly non-SD) updating rules in Appendix~\ref{sec:CKLGeneral}.

\begin{atheorem} 
	\label{thm:CKLequivalenceSD}
    Consider a class $\Pm$ and let Assumption~\ref{ass:thmifflocal} hold.
    Then, for any $y_t\in \Ym$, $\vartheta_\tbu \in \Theta$ and $p_t \in\Pm$ such that $\big(p_t(y_t)-f(y_t|\vartheta_{\tbu})\big)  \score(y_t, \vartheta_{\tbu}) \neq 0$, there exists an $\alpha > 0$ and corresponding $\delta \equiv \delta_\alpha > 0$ such that
	\begin{align*}
		\Delta_{\delta}^{\mathsf{CKL}}(\phi_\mathrm{SD} |y_t,\vartheta_\tbu,p_{t}) <0 
		\quad \iff \quad 
		p_t(y_t)>f(y_t \vert \vartheta_\tbu).
	\end{align*}
\end{atheorem}
\color{black}

Similar to \citet[Prop.~1]{blasques2015information}, Theorem~\ref{thm:CKLequivalenceSD} considers two localizations: (i) an incremental SD update $\phi_{\mathrm{SD}}(y_t,\vartheta_\tbu) = \vartheta_\tbu + \alpha \score(y_t,\vartheta_\tbu)$ through a sufficiently small learning rate $\alpha > 0$, and (ii) a sufficiently small  $\delta \equiv \delta_\alpha > 0$, given the (typically small) choice of $\alpha$, that focuses interest around $y_t$.
Theorem \ref{thm:CKLequivalenceSD} shows that SD models are \emph{not} guaranteed to be CKL reducing, in contrast with the assertion in \citet[p.\ 330]{blasques2015information} that ``every [$\ldots$] score update is locally realized KL optimal [$\ldots$] for any true density $p_t$''. 
Our use of the CKL divergence demonstrates that an improvement hinges on the (practically unverifiable) condition $p_t(y_t) > f(y_t \vert \vartheta_\tbu)$. 
This dependence on the true density is consistent with our intuition based on Figure \ref{fig:wrongdirectionmean}, where adjusting the researcher's density upwards at $y_t$ is beneficial whenever $p_t(y_t) > f(y_t \vert \vartheta_\tbu)$ (panels~\subref{subfig:p=ft+1} and~\subref{subfig:p<ft+1}), but detrimental to the model fit whenever the converse of this condition holds (panels \subref{subfig:p=ft} and~\subref{subfig:p>ft}).

The correction note by \citet{blasques2018information} acknowledges the importance of the additional condition $p_t(y_t) > f(y_t|\vartheta_\tbu)$.
However, instead of remedying the cause of the problem (i.e., localization by trimming), they 
address the impropriety issue of the weighted likelihood by adding the restriction that $p_{t} > f_{\tbu}$ on $B$ to ensure that the localized KL divergence is positive. 
This implies that localization of the
integral in \eqref{eq:LKLDefinition} now concerns the adjusted set
\begin{align}
	\label{eq:CorrectedSolution}
	\widetilde{B} := \ymdel \cap \big\{y \in \Ym: p_t(y) > f(y|\vartheta_\tbu) \big\}.
\end{align}
This gives rise to several adverse consequences. 
As the set $\widetilde{B}$ depends on the unknown true density $p_t$, it can never be verified in practice.
Moreover, $\widetilde{B}$ is the empty set when the predicted density $f_\tbu$ dominates the true density $p_t$ on $B$.
This case occurs with positive probability (panel~\subref{subfig:p>ft}, Figure~\ref{fig:wrongdirectionmean}) and means the strict improvement is no longer guaranteed.\footnote{
An alternative interpretation  of \cite{blasques2018information} is that $p_t(y_t)>f_{t|t-1}(y_t)$ is simply imposed as an additional condition throughout the original article. However, the intersection of assumptions may then be empty for the same reason that the ball $\smash{\widetilde{B}}$ may be empty. Indeed, the condition $p_t(y_t)>f_{t|t-1}(y_t)$ will fail for some $y_t$, meaning no conclusion can apparently be drawn about score-driven updates.}

In contrast, the condition $p_t(y_t)>f_{t|t-1}(y_t)$ emerges naturally as an \emph{output} of our Theorem~\ref{thm:CKLequivalenceSD}, which states that SD updates improve the CKL measure if the condition $p_t(y_t)>f_{t|t-1}(y_t)$ holds, while yielding a deterioration  if $p_t(y_t)<f_{t|t-1}(y_t)$. To date, the possibility of a deterioration of fit has not been recognized in the literature. As we show, both cases occur with positive probability, although in practice we never know which is which.

Returning to Figure \ref{fig:wrongdirectionmean}, the condition $p_t(y_t) > f(y_t \vert \vartheta_\tbu)$ holds in panels \subref{subfig:p=ft+1} and \subref{subfig:p<ft+1}, but not in panels \subref{subfig:p=ft} and \subref{subfig:p>ft}. 
Our interpretation, that local improvements in model fit cannot be obtained in all cases, is reflected in the formal results of Theorem~\ref{thm:CKLequivalenceSD}.
However, based on the observation $y_t$, the densities $p_t$ in panels \subref{subfig:p=ft+1} or \subref{subfig:p<ft+1} are more likely to correspond to the truth than the densities in panels~\subref{subfig:p=ft} and \subref{subfig:p>ft}.
This illustrates again our main result that the model fit is guaranteed to improve in the EKL sense, hence, \emph{in expectation}.

\begin{aremark}
	\label{rem:CKLForwardNotion}
	\citet[Appendix 1]{blasques2015information} discuss a \emph{forward-looking notion} that analyses whether the updated density $f_\tpu$ represents an improvement over the predicted density $f_\tbu$ in fitting the one-step-ahead true density $p_{t+1}$. 
	Theorem~\ref{thm:CKLequivalenceSD} can easily be extended to this case by replacing $p_t$ by $p_{t+1}$ to yield
	\begin{align*}
		\Delta_{\delta}^{\mathsf{CKL}}(\phi_\mathrm{SD} |y_t,\vartheta_\tbu,p_{t+1}) <0 
		\quad \iff \quad 
		p_{t+1}(y_t)>f(y_t \vert \vartheta_\tbu).
	\end{align*}
	However, as $y_t$ is derived from $p_t$, we have no information at time $t$ regarding~$p_{t+1}$, unless something is known about the dynamics of the true densities. For example, for a parametric true density $p_{t+1}(\cdot) \equiv p (\cdot\vert \lambda_{t+1})$, if the time-varying parameter process $\{\lambda_t\}_{t=1}^T$ were mean-reverting, for some values of $\omega$ and $\beta$ the predicted density $f(\cdot|\vartheta_{t+1|t})$ with $\vartheta_{t+1|t}=\omega+\beta \vartheta_{\tpu}$ may outperform $f(\cdot|\vartheta_{\tpu})$ in approximating  $p(\cdot|\lambda_{t+1})$. This analysis concerns the prediction step, however, as distinct from the updating step.
\end{aremark}

\setcounter{equation}{0}
\section{Censored KL reductions of general updating rules}
\label{sec:CKLGeneral}

The key result in \citet{blasques2015information} is that updating rules are TKL reducing if and only if they are score-equivalent in the sense of always moving in the same direction as a SD update.
Formally, \citet{blasques2015information} call an updating rule $\phi$ (almost surely) \emph{score equivalent} if $\operatorname{sign} \big( \Delta \phi(y_t, \vartheta_{\tbu}) \big) = \operatorname{sign} \big( \score(y_t, \vartheta_\tbu) \big)$ for all $(y_t, \vartheta_\tbu) \in \Ym \times \Theta$, where $\operatorname{sign}(0) := 0$.
The almost sure score equivalence does not imply score equivalence in expectations in \eqref{eqn:ExpectedScoreEquivalence}.

Here we generalize Theorem~\ref{thm:CKLequivalenceSD} and provide a similar characterization of CKL-reducing updates.
Following Appendix~\ref{sec:CKL}, we retain the restrictions $l=k=1$ and $\mathcal{S}_{t-1}=1$. 

\begin{atheorem} 
	\label{thm:CKLequivalence}
    Consider a class $\Pm$ and let Assumption~\ref{ass:thmifflocal} hold.
    Then, for any $y_t\in \Ym$, $\vartheta_\tbu \in \Theta$ and $p_t \in\Pm$ such that $\big(p_t(y_t)-f(y_t|\vartheta_{\tbu})\big)  \score(y_t, \vartheta_{\tbu}) \Delta \phi(y_t,\vartheta_{\tbu})\neq 0$, there exists a $\bar{\kappa} > 0$ and $\delta \equiv {\delta}_{\bar{\kappa}} > 0$  such that for all $\phi$ with $\vert \Delta \phi(y_t, \vartheta_\tbu)\vert < \bar{\kappa}$,
	we have
	\begin{align*}
		\Delta_{\delta}^{\mathsf{CKL}}(\phi |y_t,\vartheta_\tbu,p_t) <0 
		\quad \iff \quad 
		\score(y_t,\vartheta_\tbu) \Delta \phi(y_t,\vartheta_\tbu) \big(p_t(y_t)-f(y_t \vert \vartheta_\tbu) \big) > 0.
	\end{align*}
\end{atheorem}

The right-hand side of the equivalence in Theorem~\ref{thm:CKLequivalence} shows that whether an updating scheme $\phi$ is CKL reducing depends on the signs of three individual components: the score $\score(y_t,\vartheta_\tbu)$, the update direction $\Delta \phi(y_t,\vartheta_\tbu)$ and, as already present in Theorem~\ref{thm:CKLequivalenceSD} above, the truth through $\big(p_t(y_t)-f(y_t \vert \vartheta_\tbu) \big)$.
For score-equivalent updates, it holds that $\score(y_t,\vartheta_\tbu) \Delta \phi(y_t,\vartheta_\tbu) > 0$ such that Theorem~\ref{thm:CKLequivalenceSD} arises as a special case of Theorem~\ref{thm:CKLequivalence}.
{As in Theorem~\ref{thm:CKLequivalenceSD}, the outcome space locality parameter $\delta \equiv {\delta}_{\bar{\kappa}}$ has to be chosen much smaller than the 
	state-space locality parameter $\bar{\kappa}$; see the arguments at the end of the proof of Theorem~\ref{thm:CKLequivalence} for details.}
The upper bound $\bar{\kappa}>0$ on $\vert \Delta \phi(y_t, \vartheta_{t\vert t-1})\vert$ serves the same purpose as the linear downscaling of the random update $\Delta \phi(Y_t, \vartheta_{t\vert t-1})$ in Theorem~\ref{thm:EKLequivalenceIndepCopy}. 
In the CKL setting, this direct bound is sufficient because the update $\Delta \phi(y_t,\vartheta_{t\vert t-1})$ is deterministic conditional on $(y_t, \vartheta_{t\vert t-1})$.
Theorem~\ref{thm:CKLequivalence} establishes that if $\big(p_t(y_t) - f(y_t \vert \vartheta_\tbu) \big) > 0$ holds, then an arbitrary updating scheme is locally CKL reducing if and only if it is score equivalent to the SD update, i.e., $\score(y_t, \vartheta_\tbu) \Delta \phi(y_t, \vartheta_\tbu) > 0$ at the current $(y_t, \vartheta_\tbu)$.
However, similar to Theorem~\ref{thm:CKLequivalenceSD}, this characterization hinges on the practically unverifiable condition $\big(p_t(y_t) - f(y_t \vert \vartheta_\tbu) \big) > 0$ and is hence of limited use in practice.

In light of Remark \ref{rem:CKLForwardNotion}, a forward notion of Theorem \ref{thm:CKLequivalence} holds equivalently when considering the distribution $P_{t+1}$ with density $p_{t+1}$ by simply using the factor $\big(p_{t+1}(y_t) - f(y_t \vert \vartheta_\tbu) \big)$ on the right-hand side of the logical equivalence statement.
Hence, a separate treatment of these cases as in \citet[Appendix 1]{blasques2015information} is not required here.

\setcounter{equation}{0}
\section{Proofs of the results in Appendices~\ref{sec:CKL} and \ref{sec:CKLGeneral}}
\label{sec:ProofsSupplement}

\begin{proof}[Proof of Theorem \ref{thm:CKLequivalenceSD}]
    \sloppy
	The proof follows directly from the proof of Theorem~\ref{thm:CKLequivalence} given below by using $\phi = \phi_\mathrm{SD}$ with $\vartheta_\tpu=\phi_{\mathrm{SD}}(y_t, \vartheta_\tbu) := \vartheta_\tbu + \alpha \score(y_t,\vartheta_\tbu)$ and $\bar \kappa = \alpha$.
	Notice that for $\phi_\mathrm{SD}$ and any $\alpha > 0$ and  $\score(y_t,\vartheta_\tbu) \neq 0$, it trivially holds that $\Delta \phi_\mathrm{SD}(y_t,\vartheta_\tbu)  \score(y_t,\vartheta_\tbu) =  \alpha \score(y_t,\vartheta_\tbu)^2 > 0$. 
\end{proof}

\newcommand{\Om}{\mathcal{O}}
We require the following lemma for the proof of Theorem~\ref{thm:CKLequivalence}.

\begin{alemma}
	\label{lemma:approximationerrorIB}
Let $\Ym$ be open and convex. For a given $y_t \in \mathcal{Y}$ and $\delta > 0$, consider a function $g:\ymdellong~\to~\R$, that is Lipschitz-continuous on $\ymdellong$ and which is allowed to depend on (time-varying) parameters. 
	Then,
	\begin{align*}
		\int_\ymdellong g(y) \dd y = 2\delta g(y_t) + \Om(\delta^2),
	\end{align*}
	as $\delta \downarrow 0$.  
\end{alemma}

\begin{proof}[Proof of Lemma \ref{lemma:approximationerrorIB}]
	Using $\epsilon(y) := g(y) - g(y_t)$, we have that 
	\begin{align*}
		e:= \left\vert \int_{y_t-\delta}^{y_t+\delta} g(y) \dd y - 2 \delta g(y_t) \right\vert  &=  \left \vert \int_{y_t-\delta}^{y_t}  \epsilon(y) \dd y + \int_{y_t}^{y_t+\delta} \epsilon(y) \dd y \right \vert.
	\end{align*}
	Since $g(y)$ is Lipschitz-continuous on $\ymdellong$, there exists a constant $L>0$ such that $\vert\epsilon(y)\vert = \vert g(y)-g(y_t)\vert \leq L \vert y-y_t \vert$. Consequently, 
	\begin{align*}
		&\left \vert \int_{y_t-\delta}^{y_t}  \epsilon(y) \dd y \right \vert \leq L  \int_{y_t-\delta}^{y_t}  (y_t - y) \dd y=\frac{L}{2} \delta^2, \qquad \text{and} \\
		&\left \vert \int_{y_t}^{y_t+\delta}  \epsilon(y) \dd y \right \vert \leq L   \int_{y_t}^{y_t+\delta}   (y - y_t) \dd y =\frac{L}{2} \delta^2,
	\end{align*}
	and hence $e\leq L\delta^2=\Om(\delta^2)$, as $\delta \downarrow 0$.  
\end{proof}

\begin{proof}[Proof of Theorem \ref{thm:CKLequivalence}]	
	Fix $y_t \in \Ym$ and let $\ymdel \equiv \ymdellong = \{x \in \Ym: \vert x - y_t \vert \leq \delta\}$ be the ball of radius $\delta>0$ around the realization $y_t$ for some (small) $\delta > 0$.
    Given $\vartheta_\tbu \in \Theta$, consider an arbitrary updating rule $\phi$ such that $|\Delta \phi(y_t, \vartheta_\tbu)|\leq \bar{\kappa}$ for some sufficiently small $\bar{\kappa}>0$.
	In the following, we derive the explicit orders in which $\Delta_{\delta}^{\mathsf{CKL}}(\phi)$ shrinks to zero as the two \emph{locality parameters} $\bar{\kappa}$ and $\delta$ tend to zero. 
	The desired result then follows by taking $\bar{\kappa}>0$ and $\delta>0$ small enough. 
	
	We start by (implicitly) defining $I_B$ and $I_{B^c}$ through
	\begin{align*}
		\Delta_{\delta}^{\mathsf{CKL}}(\phi)
		&=  \mathsf{CKL}_\ymdellong(p_{t} \Vert f_\tpu) - \mathsf{CKL}_\ymdellong(p_{t} \Vert f_\tbu) \\
		&= \int_{\ymdel} \log\left(\frac{p_t(x)}{{f}(x|\vartheta_\tpu)}\right) p_t(x) \dd x 
		\;+\; \log \left(\frac{ \int_{\ymdelc} p_t(x) \dd  x}{\int_{\ymdelc}  {f}(x|\vartheta_\tpu) \dd x} \right) \int_{B^c} p_t(x) \dd x \\
		&\qquad-\int_{\ymdel} \log\left(\frac{p_t(x)}{{f}(x|\vartheta_\tbu)}\right) p_t(x) \dd x 
		\;-\; \log \left(\frac{ \int_{\ymdelc} p_t(x) \dd  x}{\int_{\ymdelc}  {f}(x|\vartheta_\tbu) \dd x} \right) \int_{B^c} p_t(x) \dd x \\
		&= - \underbrace{\int_{\ymdel} \log\left(\frac{{f}(x|\vartheta_\tpu)}{{f}(x|\vartheta_\tbu)}\right) p_t(x) \dd x}_{=:~I_B}
		\;-\;  \underbrace{\log \left(\frac{\int_{\ymdelc}  {f}(x|\vartheta_\tpu) \dd x}{\int_{\ymdelc}  {f}(x|\vartheta_\tbu) \dd x} \right) \int_{B^c} p_t(x) \dd x}_{=:~I_{B^c}},
	\end{align*}
	as the negative of the components of the censored KL divergence difference that focus on $B$ and $B^c$, respectively, such that $\Delta_{\delta}^{\mathsf{CKL}}(\phi) = - I_{B}-  I_{B^c}$.
	
	For $I_B$, Lemma~\ref{lemma:approximationerrorIB} yields 
	\begin{align} 
		\label{eq:IBfirststep}
		I_B = 2 \delta \log\left(\frac{f(y_t|\vartheta_{\tpu})}{f(y_t|\vartheta_{\tbu})}\right) p_t(y_t) + \Om(\delta^2),
	\end{align}
	as $\delta \downarrow 0$.
    Since the likelihood ratio in \eqref{eq:IBfirststep} tends to one as $\Delta \phi\equiv \Delta \phi(y_t, \vartheta_{\tbu}) $ tends to zero, we expand $\log(z)$ around $z=1$:
	\begin{align} 
		\label{eq:logexpansionB}
		\log \left(\frac{{f}(y_t|\vartheta_{\tpu})}{ f(y_t|\vartheta_{\tbu})} \right) 
		&= \frac{{f}(y_t|\vartheta_{\tpu})-f(y_t|\vartheta_\tbu)}{ f(y_t|\vartheta_{\tbu})}  + \Om\left(\bar{\kappa}^2\right),
	\end{align}
	as $\bar{\kappa}\downarrow 0$, 
	where the order of the remainder term in \eqref{eq:logexpansionB}
    follows directly from the expansion $f(y_t|\vartheta_{\tpu})-f(y_t|\vartheta_\tbu)=f(y_t|\vartheta_\tbu)\score(y_t,\vartheta_{\tbu}) \Delta \phi + \Om\left((\Delta \phi)^2\right)=\Om(\bar{\kappa})$. 

	The other elements of $I_B$ in \eqref{eq:IBfirststep} are independent of $\vartheta_\tbu$ and $\vartheta_\tpu$. 
	Hence, by substituting \eqref{eq:logexpansionB} into \eqref{eq:IBfirststep}, we find that
	\begin{align} \label{eq:IBsecond}
		I_B = 2\delta \big({f}(y_t|\vartheta_{\tpu})-f(y_t|\vartheta_{\tbu})\big) \frac{p_t(y_t)}{f(y_t|\vartheta_{\tbu})}+ \Om(\delta^2) +  \Om(\delta\bar{\kappa}^2), 
	\end{align}
	as $\delta \downarrow 0$ and $\bar{\kappa}\downarrow 0$.
	
	For $I_{B^c}$, write $\bar{F}_{\tbu} := \int_{B^c} f(x|\vartheta_{\tbu}) \dd x$, $\bar{F}_{\tpu} := \int_{B^c} f(x|\vartheta_{\tpu}) \dd x$ and  $\bar{P}_{t} := \int_{B^c} p_t(x) \dd x$. Then 
		$I_{B^c} = \log\left(\frac{\bar{F}_{\tpu}}{\bar{F}_{\tbu}}\right)  \bar{P}_{t}$.
	As $\delta \downarrow 0$, we have $\bar F_{t|t-1},\bar F_{t|t},\bar P_t \uparrow 1$. 
	Therefore, the ratio $\bar{F}_{\tpu}/\bar{F}_{\tbu}$ tends to one, motivating a Taylor expansion of $\log(z)$ around $1$.
    Specifically,
    	\begin{align}
		\log\left(\frac{\bar{F}_{\tpu}}{\bar{F}_{\tbu}}\right)
		&=  u_\delta + \Om(u_\delta^2), 
        \end{align}
    where $u_\delta := (\bar{F}_{\tpu}-\bar{F}_{\tbu})/\bar{F}_{\tbu}$. 
    By Lemma~\ref{lemma:approximationerrorIB},
    \begin{align*}
        \bar{F}_{\tpu}-\bar{F}_{\tbu} =- \int_B \big(f(x|\vartheta_{\tpu}) - f(x|\vartheta_{\tbu}) \big) \dd x  
        = - 2\delta \big(f(y_t|\vartheta_{\tpu})-f(y_t|\vartheta_{\tbu})\big) + \Om(\delta^2),
    \end{align*}
 and hence
 \begin{align*}
     u_\delta  = - \frac{2\delta \big(f(y_t|\vartheta_{\tpu})-f(y_t|\vartheta_{\tbu})\big) + \Om(\delta^2)}{\bar{F}_{\tbu}}= \Om(\delta \bar{\kappa}) + \Om(\delta^2),
 \end{align*}
recalling $f(y_t|\vartheta_{\tpu})-f(y_t|\vartheta_{\tbu}) =\Om(\bar{\kappa})$.
Moreover, 
$u_\delta^2 = \Om(\delta^2 \bar{\kappa}^2)  + \Om(\delta^3 \bar{\kappa}) + \Om(\delta^4)$.
Consequently,  
	\begin{align*} 
		I_{B^c} 
		&= \Big( 2 \delta\big({f}(y_t|\vartheta_{\tbu}) - {f}(y_t|\vartheta_{\tpu})\big) + \Om(\delta^2) \Big) \frac{ \bar{P}_t }{\bar{F}_{\tbu}} + \Om\left(\delta^2 \bar{\kappa}^2 \right) + \Om(\delta^3 \bar{\kappa}) + \Om(\delta^4), 
	\end{align*}
	as $\delta \downarrow 0$ and $\bar{\kappa}\downarrow 0$,
    in which
    \begin{align*}
		\frac{\bar{P}_t }{\bar{F}_{\tbu}} =  \frac{1- 2\delta p_t(y_t) -  \Om(\delta^2) }{1- 2\delta f(y_t \vert \vartheta_{\tbu}) -  \Om(\delta^2)}
		&=1+ \Om(\delta),
	\end{align*}
	as $\delta \downarrow 0$, by Lemma~\ref{lemma:approximationerrorIB}. Hence,
    \begin{align} \label{eq:IBc}
        I_{B^c} 
		&= - 2 \delta\big( {f}(y_t|\vartheta_{\tpu}) - {f}(y_t|\vartheta_{\tbu}) \big)\frac{ \bar{P}_t }{\bar{F}_{\tbu}} + \Om(\delta^2), 
    \end{align}
    as $\delta \downarrow 0$ and $\bar{\kappa}\downarrow 0$.
	Combining the expressions for $I_B$ and $I_{B^c}$ in \eqref{eq:IBsecond} and \eqref{eq:IBc}, respectively, yields
	\begin{align*}
		\Delta_{\delta}^{\mathsf{CKL}}(\phi) 
		&= - (I_B + I_{B^c}) \\
		&= -  2\delta \big({f}(y_t|\vartheta_{\tpu})-f(y_t|\vartheta_{\tbu})\big)\left(\frac{p_t(y_t)}{f(y_t|\vartheta_{\tbu})} - \frac{\bar{P}_t }{\bar{F}_{\tbu}}\right) +  \Om(\delta^2) +  \Om(\delta\bar{\kappa}^2) \\
		&= -2\delta \big({f}(y_t|\vartheta_{\tpu})-f(y_t|\vartheta_{\tbu})\big)\left(\frac{p_t(y_t)}{f(y_t|\vartheta_{\tbu})} - 1 + \Om(\delta) \right) +  \Om(\delta^2) +  \Om(\delta\bar{\kappa}^2) \\
		&= -2\delta \frac{ f(y_t|\vartheta_{\tpu}) - f(y_t|\vartheta_{\tbu})}{f(y_t\vert \vartheta_{\tbu})} \big( p_t(y_t)-f(y_t|\vartheta_{\tbu}) \big) +  \Om(\delta^2) +  \Om(\delta\bar{\kappa}^2),
	\end{align*}
	as $\delta \downarrow 0$ and $\bar{\kappa} \downarrow 0$.
	
	Using the expansion in \eqref{eq:logexpansionB} for the fraction $\big({f}(y_t|\vartheta_{\tpu}) - f(y_t|\vartheta_{\tbu})\big) \big/ f(y_t\vert \vartheta_{\tbu})$, we~get
	\begin{align*}
		\Delta_{\delta}^{\mathsf{CKL}}(\phi)
		&= - 2\delta  \Big(\log f(y_t|\vartheta_{\tpu})-\log f(y_t|\vartheta_{\tbu})-\Om(\bar{\kappa}^2) \Big) \big(p_t(y_t)-f(y_t|\vartheta_{\tbu})\big)  \\
		& \qquad +  \Om(\delta^2) +  \Om(\delta\bar{\kappa}^2) \\
		&= - 2\delta  \big(\log f(y_t|\vartheta_{\tpu})-\log f(y_t|\vartheta_{\tbu})\big) \big({p}_t(y_t)-f(y_t|\vartheta_{\tbu})\big) +  \Om(\delta^2) +  \Om(\delta\bar{\kappa}^2),\\
		&= -2\delta \big({p}_t(y_t)-f(y_t|\vartheta_{\tbu})\big) \Big(\score(y_t, \vartheta_{\tbu}) \Delta \phi(y_t, \vartheta_{\tbu}) + \Om(\bar{\kappa}^2) \Big) +  \Om(\delta^2) +  \Om(\delta\bar{\kappa}^2) \nonumber \\
		&= -2\delta \big({p}_t(y_t)-f(y_t|\vartheta_{\tbu})\big) \score(y_t, \vartheta_{\tbu}) \Delta \phi(y_t, \vartheta_{\tbu})  +  \Om(\delta^2) +  \Om(\delta\bar{\kappa}^2),
	\end{align*}
	as $\delta \downarrow 0$ and $\bar{\kappa} \downarrow 0$.
	
	We continue to define the term
	\begin{align}
		\label{eqn:DefinitionR}
		r(y_t, \vartheta_{\tbu}, p_t) := \big(p_t(y_t)-f(y_t|\vartheta_{\tbu})\big)  \score(y_t, \vartheta_{\tbu}) \Delta \phi(y_t,\vartheta_{\tbu}),
	\end{align}
	which is of order $\Om(\bar \kappa)$ for $\bar \kappa \downarrow 0$ as $|\Delta \phi(y_t,\vartheta_{\tbu})| = |\vartheta_{\tpu} - \vartheta_{\tbu}| \le \bar \kappa$, whereas the terms $\big(p_t(y_t)-f(y_t|\vartheta_{\tbu})\big)$ and $ \score(y_t, \vartheta_{\tbu})$ are independent of $\delta$ and $\bar \kappa$.
	Taking the following product then yields that
	\begin{align*}
		&\Delta_\delta^\mathsf{CKL}(\phi\vert y_t,\vartheta_\tbu,p_t) \, r(y_t, \vartheta_{\tbu}, p_t) \\
		&\qquad= -2 \delta r^2(y_t, \vartheta_{\tbu}, p_t) + r(y_t, \vartheta_{\tbu}, p_t) \big( \Om(\delta^2) + \Om(\delta \bar \kappa^2) \big)\\
		&\qquad= -2 \delta r^2(y_t, \vartheta_{\tbu}, p_t) + \Om(\delta^2 \bar \kappa) + \Om( \delta \bar \kappa^3),
	\end{align*}
	where we notice that the first term $-2 \delta r^2(y_t, \vartheta_{\tbu}, p_t) =  \Om(\delta \bar \kappa^2)$ is strictly negative and given that we can choose $\delta = \delta_{\bar \kappa}$ much smaller than $\bar \kappa$, the $\Om(\delta \bar \kappa^2)$ term is of lower order than (i.e., it does not vanish as fast as) the following $\Om(\delta^2 \bar \kappa) + \Om(\delta \bar \kappa^3)$ terms.
	
	Hence, there exists $\bar{\kappa} > 0$ and $\delta  = \delta_{\bar \kappa} > 0$ (where the choice of $\delta = \delta_{\bar \kappa}$ depends on the choice of $\bar \kappa$) such that for all updates $\phi$ with $|\Delta \phi(y_t, \vartheta_\tbu)|\leq \bar{\kappa}$,
	it holds that $\Delta_\delta^\mathsf{CKL}(\phi\vert y_t,\vartheta_\tbu,p_t) r(y_t, \vartheta_{\tbu}, p_t) < 0$, implying that $\Delta_\delta^\mathsf{CKL}(\phi|y_t,\vartheta_{\tbu},p_t)$ and $r(y_t, \vartheta_{\tbu}, p_t)  =  \big(p_t(y_t)-f(y_t|\vartheta_{\tbu})\big) \score(y_t, \vartheta_{\tbu}) \Delta \phi(y_t,\vartheta_{\tbu})$ have opposite signs, which is exactly the statement that had to be shown.\footnote{A notable difference between the characterization results of \citet[Prop.~2]{blasques2015information} and our Theorem~\ref{thm:CKLequivalence} shown here is the \emph{locality} of the required score-equivalence.
		In Theorem~\ref{thm:CKLequivalence}, we impose \emph{local} score-equivalence that only has to hold at the pair $(y_t, \vartheta_\tbu)$ whereas \citet[Proposition 2]{blasques2015information} claim a \emph{global} score equivalence for all $(y, \vartheta) \in \mathcal{Y} \times \Theta$.
		This difference is not caused by the different localization methods (i.e., censoring opposed to trimming), but by an inaccuracy in the proof of \citet[Prop.~2]{blasques2015information}:
		Using the notation of \citet{blasques2015information}, in the ``only if'' direction on page 340, the neighborhoods $FY$ and $Y_{\delta_Y}(y_t)$ are not necessarily overlapping, which could be fixed by imposing a \emph{local} score-equivalence opposed to their \emph{global} notion.}
	   (A separate proof for the two directions in the ``$\Longleftrightarrow$'' statement could be carried out as at the end of the proof of Theorem~\ref{thm:EKLequivalenceIndepCopy}.) 
\end{proof}

\singlespacing
\spacingset{1} 
\setlength{\bibsep}{2pt} 
\putbib[Bibliography-LSPS-v2]
\end{bibunit}


\begin{thebibliography}{}

\bibitem[\protect\citeauthoryear{Amisano and Giacomini}{Amisano and
  Giacomini}{2007}]{amisano2007comparing}
Amisano, G. and R.~Giacomini (2007).
\newblock Comparing density forecasts via weighted likelihood ratio tests.
\newblock {\em Journal of Business \& Economic Statistics\/}~{\em 25\/}(2),
  177--190.

\bibitem[\protect\citeauthoryear{Artemova, Blasques, van Brummelen, and
  Koopman}{Artemova et~al.}{2022a}]{artemova2022score1}
Artemova, M., F.~Blasques, J.~van Brummelen, and S.~J. Koopman (2022a).
\newblock Score-driven models: {M}ethodology and theory.
\newblock In {\em Oxford Research Encyclopedia of Economics and Finance}.

\bibitem[\protect\citeauthoryear{Artemova, Blasques, van Brummelen, and
  Koopman}{Artemova et~al.}{2022b}]{artemova2022score2}
Artemova, M., F.~Blasques, J.~van Brummelen, and S.~J. Koopman (2022b).
\newblock Score-driven models: {M}ethods and applications.
\newblock In {\em Oxford Research Encyclopedia of Economics and Finance}.

\bibitem[\protect\citeauthoryear{Beutner, Lin, and Lucas}{Beutner
  et~al.}{2023}]{beutner2023consistency}
Beutner, E.~A., Y.~Lin, and A.~Lucas (2023).
\newblock Consistency, distributional convergence, and optimality of
  score-driven filters.
\newblock {\em Tinbergen Institute Discussion Paper No. 2023-051/III\/}.
\newblock
  \href{https://papers.tinbergen.nl/23051.pdf}{https://papers.tinbergen.nl/23051.pdf}.

\bibitem[\protect\citeauthoryear{Blasques, Francq, and Laurent}{Blasques
  et~al.}{2023}]{blasques2023quasi}
Blasques, F., C.~Francq, and S.~Laurent (2023).
\newblock Quasi score-driven models.
\newblock {\em Journal of Econometrics\/}~{\em 234\/}(1), 251--275.

\bibitem[\protect\citeauthoryear{Blasques, Koopman, and Lucas}{Blasques
  et~al.}{2015}]{blasques2015information}
Blasques, F., S.~J. Koopman, and A.~Lucas (2015).
\newblock Information-theoretic optimality of observation-driven time series
  models for continuous responses.
\newblock {\em Biometrika\/}~{\em 102\/}(2), 325--343.

\bibitem[\protect\citeauthoryear{Blasques, Koopman, and Lucas}{Blasques
  et~al.}{2018}]{blasques2018information}
Blasques, F., S.~J. Koopman, and A.~Lucas (2018).
\newblock Amendments and corrections: Information-theoretic optimality of
  observation-driven time series models for continuous responses.
\newblock {\em Biometrika\/}~{\em 105\/}(3), 753.

\bibitem[\protect\citeauthoryear{Blasques, Lucas, and van Vlodrop}{Blasques
  et~al.}{2021}]{blasques2021finite}
Blasques, F., A.~Lucas, and A.~C. van Vlodrop (2021).
\newblock Finite sample optimality of score-driven volatility models: {S}ome
  {M}onte {C}arlo evidence.
\newblock {\em Econometrics and Statistics\/}~{\em 19}, 47--57.

\bibitem[\protect\citeauthoryear{Boyd and Vandenberghe}{Boyd and
  Vandenberghe}{2004}]{boyd2004convex}
Boyd, S.~P. and L.~Vandenberghe (2004).
\newblock {\em Convex optimization}.
\newblock CUP.

\bibitem[\protect\citeauthoryear{Catania, D’Innocenzo, and Luati}{Catania
  et~al.}{2026}]{catania2026unobserved}
Catania, L., E.~D’Innocenzo, and A.~Luati (2026).
\newblock Unobserved component models, approximate filters and dynamic adaptive
  mixture models.
\newblock {\em Journal of Econometrics\/}~{\em 253}, 106155.

\bibitem[\protect\citeauthoryear{Creal, Koopman, and Lucas}{Creal
  et~al.}{2013}]{creal2013generalized}
Creal, D., S.~J. Koopman, and A.~Lucas (2013).
\newblock Generalized autoregressive score models with applications.
\newblock {\em Journal of Applied Econometrics\/}~{\em 28\/}(5), 777--795.

\bibitem[\protect\citeauthoryear{Creal, Koopman, Lucas, and Zamojski}{Creal
  et~al.}{2024}]{Creal2024GMM}
Creal, D., S.~J. Koopman, A.~Lucas, and M.~Zamojski (2024).
\newblock Observation-driven filtering of time-varying parameters using moment
  conditions.
\newblock {\em Journal of Econometrics\/}~{\em 238\/}(2), 105635.

\bibitem[\protect\citeauthoryear{De~Punder, Diks, Laeven, and van
  Dijk}{De~Punder et~al.}{2026}]{depunder2023localzing}
De~Punder, R., C.~Diks, R.~Laeven, and D.~J. van Dijk (2026).
\newblock Localizing strictly proper scoring rules.
\newblock {\em Journal of the American Statistical Association\/}.
\newblock
  \href{https://doi.org/10.1080/01621459.2025.2576189}{https://doi.org/10.1080/01621459.2025.2576189}.

\bibitem[\protect\citeauthoryear{Delle~Monache, De~Polis, and
  Petrella}{Delle~Monache et~al.}{2023}]{delle2023modeling}
Delle~Monache, D., A.~De~Polis, and I.~Petrella (2023).
\newblock Modeling and forecasting macroeconomic downside risk.
\newblock {\em Journal of Business \& Economic Statistics\/}~{\em 42\/}(3),
  1010–1025.

\bibitem[\protect\citeauthoryear{Diks, Panchenko, and van Dijk}{Diks
  et~al.}{2011}]{diks2011likelihood}
Diks, C., V.~Panchenko, and D.~J. van Dijk (2011).
\newblock Likelihood-based scoring rules for comparing density forecasts in
  tails.
\newblock {\em Journal of Econometrics\/}~{\em 163\/}(2), 215--230.

\bibitem[\protect\citeauthoryear{Durbin and Koopman}{Durbin and
  Koopman}{2012}]{durbin2012time}
Durbin, J. and S.~J. Koopman (2012).
\newblock {\em Time series analysis by state space methods}.
\newblock OUP.

\bibitem[\protect\citeauthoryear{D’Innocenzo, Lucas, Schwaab, and
  Zhang}{D’Innocenzo et~al.}{2024}]{d2024modeling}
D’Innocenzo, E., A.~Lucas, B.~Schwaab, and X.~Zhang (2024).
\newblock Modeling extreme events: Time-varying extreme tail shape.
\newblock {\em Journal of Business \& Economic Statistics\/}~{\em 42\/}(3),
  903--917.

\bibitem[\protect\citeauthoryear{Gasperoni, Luati, Paci, and
  D’Innocenzo}{Gasperoni et~al.}{2023}]{gasperoni_score-driven_2023}
Gasperoni, F., A.~Luati, L.~Paci, and E.~D’Innocenzo (2023).
\newblock Score-driven modeling of spatio-temporal data.
\newblock {\em Journal of the American Statistical Association\/}~{\em
  118\/}(542), 1066--1077.

\bibitem[\protect\citeauthoryear{Gneiting and Raftery}{Gneiting and
  Raftery}{2007}]{gneiting2007strictly}
Gneiting, T. and A.~Raftery (2007).
\newblock Strictly proper scoring rules, prediction, and estimation.
\newblock {\em Journal of the American Statistical Association\/}~{\em
  102\/}(477), 359--378.

\bibitem[\protect\citeauthoryear{Gneiting and Ranjan}{Gneiting and
  Ranjan}{2011}]{gneiting_comparing_2011}
Gneiting, T. and R.~Ranjan (2011).
\newblock Comparing density forecasts using threshold- and quantile-weighted
  scoring rules.
\newblock {\em Journal of Business \& Economic Statistics\/}~{\em 29\/}(3),
  411--422.

\bibitem[\protect\citeauthoryear{Goodfellow, Bengio, Courville, and
  Bengio}{Goodfellow et~al.}{2016}]{goodfellow2016deep}
Goodfellow, I., Y.~Bengio, A.~Courville, and Y.~Bengio (2016).
\newblock {\em Deep learning}.
\newblock MIT press.

\bibitem[\protect\citeauthoryear{Gorgi, Lauria, and Luati}{Gorgi
  et~al.}{2024}]{Gorgi2023}
Gorgi, P., C.~S. Lauria, and A.~Luati (2024).
\newblock On the optimality of score-driven models.
\newblock {\em Biometrika\/}~{\em 111\/}(3), 865--880.

\bibitem[\protect\citeauthoryear{Harvey}{Harvey}{2013}]{harvey2013dynamic}
Harvey, A.~C. (2013).
\newblock {\em Dynamic models for volatility and heavy tails: {W}ith
  applications to financial and economic time series}.
\newblock CUP.

\bibitem[\protect\citeauthoryear{Harvey}{Harvey}{2022}]{harvey2022score}
Harvey, A.~C. (2022).
\newblock Score-driven time series models.
\newblock {\em Annual Review of Statistics and Its Application\/}~{\em 9},
  321--342.

\bibitem[\protect\citeauthoryear{Harvey and Luati}{Harvey and
  Luati}{2014}]{harvey2014filtering}
Harvey, A.~C. and A.~Luati (2014).
\newblock Filtering with heavy tails.
\newblock {\em Journal of the American Statistical Association\/}~{\em
  109\/}(507), 1112--1122.

\bibitem[\protect\citeauthoryear{Hol{\`y} and Tomanov{\'a}}{Hol{\`y} and
  Tomanov{\'a}}{2022}]{holy2022modeling}
Hol{\`y}, V. and P.~Tomanov{\'a} (2022).
\newblock Modeling price clustering in high-frequency prices.
\newblock {\em Quantitative Finance\/}~{\em 22\/}(9), 1649--1663.

\bibitem[\protect\citeauthoryear{Kingma and Ba}{Kingma and
  Ba}{2014}]{kingma2014adam}
Kingma, D.~P. and J.~Ba (2014).
\newblock Adam: A method for stochastic optimization.
\newblock {\em Preprint\/}.
\newblock
  \href{https://arxiv.org/abs/1412.6980}{https://arxiv.org/abs/1412.6980}.

\bibitem[\protect\citeauthoryear{Koopman, Lucas, and Scharth}{Koopman
  et~al.}{2016}]{koopman2016predicting}
Koopman, S.~J., A.~Lucas, and M.~Scharth (2016).
\newblock Predicting time-varying parameters with parameter-driven and
  observation-driven models.
\newblock {\em Review of Economics and Statistics\/}~{\em 98\/}(1), 97--110.

\bibitem[\protect\citeauthoryear{Lange, van Os, and van Dijk}{Lange
  et~al.}{2025}]{lange2022robust}
Lange, R.-J., B.~van Os, and D.~J. van Dijk (2025).
\newblock Implicit score-driven filters for time-varying parameter models.
\newblock {\em Preprint\/}.
\newblock
  \href{https://arxiv.org/abs/2512.02744}{https://arxiv.org/abs/2512.02744}.

\bibitem[\protect\citeauthoryear{Mai and Johansson}{Mai and
  Johansson}{2021}]{mai2021stability}
Mai, V.~V. and M.~Johansson (2021).
\newblock Stability and convergence of stochastic gradient clipping: Beyond
  {L}ipschitz continuity and smoothness.
\newblock In {\em International Conference on Machine Learning}, pp.\
  7325--7335. PMLR.

\bibitem[\protect\citeauthoryear{Nesterov}{Nesterov}{2018}]{nesterov2018lectures}
Nesterov, Y. (2018).
\newblock {\em Lectures on convex optimization}.
\newblock Springer.

\bibitem[\protect\citeauthoryear{Nicholson}{Nicholson}{1979}]{nicholson1979eigenvalue}
Nicholson, D.~W. (1979).
\newblock Eigenvalue bounds for {$AB+ BA$}, with {$A, B$} positive definite
  matrices.
\newblock {\em Linear Algebra and its Applications\/}~{\em 24}, 173--184.

\bibitem[\protect\citeauthoryear{Van~der Vaart}{Van~der
  Vaart}{1998}]{van1998asymptotic}
Van~der Vaart, A.~W. (1998).
\newblock {\em Asymptotic statistics}.
\newblock CUP.

\end{thebibliography}


\begin{thebibliography}{}

\bibitem[\protect\citeauthoryear{Amisano and Giacomini}{Amisano and
  Giacomini}{2007}]{amisano2007comparing}
Amisano, G. and R.~Giacomini (2007).
\newblock Comparing density forecasts via weighted likelihood ratio tests.
\newblock {\em Journal of Business \& Economic Statistics\/}~{\em 25\/}(2),
  177--190.

\bibitem[\protect\citeauthoryear{Bertsekas}{Bertsekas}{2016}]{bertsekas2016nonlinear}
Bertsekas, D.~P. (2016).
\newblock {\em Nonlinear programming}.
\newblock Athena Scientific.

\bibitem[\protect\citeauthoryear{Blasques, Francq, and Laurent}{Blasques
  et~al.}{2023}]{blasques2023quasi}
Blasques, F., C.~Francq, and S.~Laurent (2023).
\newblock Quasi score-driven models.
\newblock {\em Journal of Econometrics\/}~{\em 234\/}(1), 251--275.

\bibitem[\protect\citeauthoryear{Blasques, Koopman, and Lucas}{Blasques
  et~al.}{2015}]{blasques2015information}
Blasques, F., S.~J. Koopman, and A.~Lucas (2015).
\newblock Information-theoretic optimality of observation-driven time series
  models for continuous responses.
\newblock {\em Biometrika\/}~{\em 102\/}(2), 325--343.

\bibitem[\protect\citeauthoryear{Blasques, Koopman, and Lucas}{Blasques
  et~al.}{2018}]{blasques2018information}
Blasques, F., S.~J. Koopman, and A.~Lucas (2018).
\newblock Amendments and corrections: Information-theoretic optimality of
  observation-driven time series models for continuous responses.
\newblock {\em Biometrika\/}~{\em 105\/}(3), 753.

\bibitem[\protect\citeauthoryear{Bollerslev}{Bollerslev}{1987}]{bollerslev1987conditionally}
Bollerslev, T. (1987).
\newblock A conditionally heteroskedastic time series model for speculative
  prices and rates of return.
\newblock {\em The Review of Economics and Statistics\/}~{\em 69\/}(3),
  542--547.

\bibitem[\protect\citeauthoryear{Creal, Koopman, Lucas, and Zamojski}{Creal
  et~al.}{2024}]{Creal2024GMM}
Creal, D., S.~J. Koopman, A.~Lucas, and M.~Zamojski (2024).
\newblock Observation-driven filtering of time-varying parameters using moment
  conditions.
\newblock {\em Journal of Econometrics\/}~{\em 238\/}(2), 105635.

\bibitem[\protect\citeauthoryear{De~Punder, Diks, Laeven, and van
  Dijk}{De~Punder et~al.}{2026}]{depunder2023localzing}
De~Punder, R., C.~Diks, R.~Laeven, and D.~J. van Dijk (2026).
\newblock Localizing strictly proper scoring rules.
\newblock {\em Journal of the American Statistical Association\/}.
\newblock
  \href{https://doi.org/10.1080/01621459.2025.2576189}{https://doi.org/10.1080/01621459.2025.2576189}.

\bibitem[\protect\citeauthoryear{Diks, Panchenko, and van Dijk}{Diks
  et~al.}{2011}]{diks2011likelihood}
Diks, C., V.~Panchenko, and D.~J. van Dijk (2011).
\newblock Likelihood-based scoring rules for comparing density forecasts in
  tails.
\newblock {\em Journal of Econometrics\/}~{\em 163\/}(2), 215--230.

\bibitem[\protect\citeauthoryear{Donker~van Heel, Lange, van Os, and van
  Dijk}{Donker~van Heel et~al.}{2025}]{donker2025stability}
Donker~van Heel, S.~W., R.-J. Lange, B.~van Os, and D.~J. van Dijk (2025).
\newblock Stability and performance guarantees for misspecified multivariate
  score-driven filters.
\newblock {\em Preprint\/}.
\newblock
  \href{https://arxiv.org/abs/2502.05021}{https://arxiv.org/abs/2502.05021}.

\bibitem[\protect\citeauthoryear{Durbin and Koopman}{Durbin and
  Koopman}{2012}]{durbin2012time}
Durbin, J. and S.~J. Koopman (2012).
\newblock {\em Time series analysis by state space methods}.
\newblock OUP.

\bibitem[\protect\citeauthoryear{Gneiting and Ranjan}{Gneiting and
  Ranjan}{2011}]{gneiting_comparing_2011}
Gneiting, T. and R.~Ranjan (2011).
\newblock Comparing density forecasts using threshold- and quantile-weighted
  scoring rules.
\newblock {\em Journal of Business \& Economic Statistics\/}~{\em 29\/}(3),
  411--422.

\bibitem[\protect\citeauthoryear{Gorgi, Lauria, and Luati}{Gorgi
  et~al.}{2024}]{Gorgi2023}
Gorgi, P., C.~S. Lauria, and A.~Luati (2024).
\newblock On the optimality of score-driven models.
\newblock {\em Biometrika\/}~{\em 111\/}(3), 865--880.

\bibitem[\protect\citeauthoryear{Kolda and Bader}{Kolda and
  Bader}{2009}]{kolda_tensor_2009}
Kolda, T.~G. and B.~W. Bader (2009).
\newblock Tensor {decompositions} and {applications}.
\newblock {\em SIAM Review\/}~{\em 51\/}(3), 455--500.

\bibitem[\protect\citeauthoryear{Koopman, Lucas, and Scharth}{Koopman
  et~al.}{2016}]{koopman2016predicting}
Koopman, S.~J., A.~Lucas, and M.~Scharth (2016).
\newblock Predicting time-varying parameters with parameter-driven and
  observation-driven models.
\newblock {\em Review of Economics and Statistics\/}~{\em 98\/}(1), 97--110.

\bibitem[\protect\citeauthoryear{Lange, van Os, and van Dijk}{Lange
  et~al.}{2025}]{lange2022robust}
Lange, R.-J., B.~van Os, and D.~J. van Dijk (2025).
\newblock Implicit score-driven filters for time-varying parameter models.
\newblock {\em Preprint\/}.
\newblock
  \href{https://arxiv.org/abs/2512.02744}{https://arxiv.org/abs/2512.02744}.

\end{thebibliography}
\end{document}